\documentclass[a4paper,11pt]{article}
\usepackage{dsfont}
\usepackage{amssymb}
\usepackage{latexsym}
\usepackage{amsmath}
\usepackage{color}
\usepackage{comment}
\usepackage{amsthm}
\usepackage[dvips]{graphicx}
\usepackage{here}
\usepackage{layout} 
\usepackage{ulem}
\usepackage{enumerate}
\usepackage{bm}
\usepackage{bbm} 
\usepackage[bbgreekl]{mathbbol} 
\usepackage{authblk}
\usepackage{setspace} 
\newif\ifrs
\rstrue
\ifrs \usepackage{mathrsfs} \fi  
\newif\ifcol
\coltrue 

\newtheorem{theorem}{Theorem}[section]

\newtheorem{lemma}[theorem]{Lemma}

\newtheorem{proposition}[theorem]{Proposition}

\newtheorem{remark}[theorem]{Remark}

\numberwithin{equation}{section}
\newtheorem{theorem*}{Theorem}
\newtheorem{ass*}[theorem*]{Assumption}
\newtheorem{note*}[theorem*]{Note}
\newtheorem{lemma*}[theorem*]{Lemma}
\newtheorem{definition*}[theorem*]{Definition}
\newtheorem{proposition*}[theorem*]{Proposition}
\newtheorem{corollary*}[theorem*]{Corollary}
\newtheorem{remark*}[theorem*]{Remark}
\newtheorem{example*}[theorem*]{Example}
\numberwithin{equation}{section}
\newif\ifcol
\colfalse
\ifcol
\newcommand{\colorr}{\color[rgb]{0.8,0,0}}
\newcommand{\colorg}{\color[rgb]{0,0.5,0}}
\newcommand{\colorb}{\color[rgb]{0,0,0.8}}

\newcommand{\colorn}{\color[rgb]{1,1,1}}

\newcommand{\colory}{\color{yellow}}
\newcommand{\coloro}{\color[rgb]{1,0.4,0}}

\else
\newcommand{\colorb}{\color{black}}

\newcommand{\colorr}{\color{black}}
\newcommand{\colorg}{\color{black}}

\newcommand{\colorn}{\color{black}}
\newcommand{\colory}{\color{black}}
\newcommand{\coloro}{\color{black}}

\fi
\newif\ifcol
\colfalse
\ifcol
\newcommand{\cred}{\color[rgb]{0.8,0,0}}

\newcommand{\cblue}{\color[rgb]{0,0,0.8}}
\else
\newcommand{\cred}{\color{black}}
\newcommand{\cblue}{\color{black}}
\fi
\newif\ifcol
\colfalse
\ifcol

\newcommand{\tred}{\color[rgb]{0.8,0,0}}
\newcommand{\tgreen}{\color[rgb]{0,0.5,0}}
\newcommand{\tblue}{\color[rgb]{0,0,0.8}}
\else
\newcommand{\tred}{\color{black}}
\newcommand{\tgreen}{\color{black}}
\newcommand{\tblue}{\color{black}}
\fi
\newif\ifcol
\colfalse
\ifcol

\newcommand{\fred}{\color[rgb]{0.8,0,0}}

\else
\newcommand{\fred}{\color{black}}
\fi
\newif\ifcol
\colfalse
\ifcol

\newcommand{\vred}{\color{black}}
\newcommand{\vgreen}{\color[rgb]{0,0.5,0}}
\newcommand{\vblue}{\color[rgb]{0,0,0.8}}
\newcommand{\vredd}{\color[rgb]{0.8,0,0}}
\else
\newcommand{\vred}{\color{black}}
\newcommand{\vgreen}{\color{black}}
\newcommand{\vblue}{\color{black}}
\newcommand{\vredd}{\color{black}}
\fi
\newif\ifcol
\colfalse
\ifcol

\newcommand{\xred}{\color[rgb]{0.8,0,0}}

\newcommand{\xblue}{\color[rgb]{0,0,0.8}}

\else
\newcommand{\xred}{\color{black}}
\newcommand{\xblue}{\color{black}}
\fi
\excludecomment{en-text}
\includecomment{jp-text}
\includecomment{comment}
\def\infm{{\infty\text{--}}}
\def\inftym{\infm}
\def\koko{{\coloroy{koko}}}
\def\bd{\begin{description}}
\def\ed{\end{description}}

\def\D2{\bbD_{2,\infty-}}

\def\tj{{t_j}}
\def\tjm{{t_{j-1}}}

\def\D{{\bf D}}

\def\F{{\bf F}}

\def\H{{\bf H}}

\def\Y{{\bf Y}}

\def\cala{{\cal A}}

\def\cald{{\cal D}}

\def\calf{{\cal F}}

\def\calh{{\cal H}}

\def\call{{\cal L}}

\def\caln{{\cal N}}

\def\calx{{\cal X}}
\def\caly{{\cal Y}}

\def\ds{\displaystyle}
\def\yeq{\>=\>}

\def\sfd{{\sf d}}
\def\sfp{{\sf p}}
\def\sfr{{\sf r}}

\def\ep{\epsilon}
\def\half{\frac{1}{2}}

\def\down{\downarrow}

\def\y{\vspace*{3mm}\\}
\def\halflineskip{\vspace*{3mm}}
\def\nn{\nonumber}
\def\be{\begin{equation}}
\def\ee{\end{equation}}
\def\bea{\begin{eqnarray}}
\def\eea{\end{eqnarray}}
\def\beas{\begin{eqnarray*}}
\def\eeas{\end{eqnarray*}}
\def\bi{\begin{itemize}}
\def\ei{\end{itemize}}
\def\im{\item}
\def\bd{\begin{description}}
\def\ed{\end{description}}
\newcommand{\bbD}{{\mathbb D}}
\newcommand{\bbE}{{\mathbb E}}
\newcommand{\bbF}{{\mathbb F}}

\newcommand{\bbH}{{\mathbb H}}
\newcommand{\bbI}{{\mathbb I}}

\newcommand{\bbL}{{\mathbb L}}

\newcommand{\bbN}{{\mathbb N}}

\newcommand{\bbR}{{\mathbb R}}

\newcommand{\bbT}{{\mathbb T}}
\newcommand{\bbU}{{\mathbb U}}

\newcommand{\bbY}{{\mathbb Y}}

\newcommand{\lone}{5.1 }
\newcommand{\ltwo}{5.2 }
\newcommand{\lthree}{5.3 }
\newcommand{\lfour}{5.4 }
\newcommand{\lfive}{5.5 }
\newcommand{\lsix}{5.6 }

\newcommand{\bH}{{\bf H}}
\newcommand{\ol}{\overline}
\newcommand{\sfA}{{\sf A}}
\renewcommand{\colory}{\color{yellow}}
\renewcommand{\koko}{{\colory koko}}
\newcommand{\wh}{\widehat}
\newcommand{\wt}{\widetilde}
\newcommand{\bbi}{{\mathbb i}}
\newcommand{\yfbox}{\fbox}

\begin{document}

\title{{\vblue Quasi-likelihood analysis for adaptive estimation of a degenerate diffusion process}
\footnote{
This work was in part supported by 
Japan Science and Technology Agency CREST JPMJCR2115; 
Japan Society for the Promotion of Science Grants-in-Aid for Scientific Research 
No. 23H03354 (Scientific Research);  
and by a Cooperative Research Program of the Institute of Statistical Mathematics. 
}
}
\author[1]{Arnaud Gloter}
\author[2,3]{Nakahiro Yoshida}
\affil[1]{Laboratoire de Math\'ematiques et Mod\'elisation d'Evry, Universit\'e d'Evry
\footnote{
Laboratoire de Math\'ematiques et Mod\'elisation d'Evry, CNRS, Univ Evry, 
Universit\'e Paris-Saclay, 91037, Evry, France. e-mail: arnaud.gloter@univ-evry.fr}
        }
\affil[2]{Graduate School of Mathematical Sciences, University of Tokyo
\footnote{Graduate School of Mathematical Sciences, University of Tokyo: 3-8-1 Komaba, Meguro-ku, Tokyo 153-8914, Japan. e-mail: nakahiro@ms.u-tokyo.ac.jp}
        }
\affil[3]{Japan Science and Technology Agency CREST
        }
\maketitle
\ \\
{\it Summary} 
{\vblue 
The adaptive quasi-likelihood analysis is developed for a degenerate diffusion process. 
Asymptotic normality and moment convergence are proved for the quasi-maximum likelihood estimators and 
quasi-Bayesian estimators, in the adaptive scheme. 
}

\section{Introduction}\label{202002201325}
{\vblue 
We consider parametric estimation of a degenerate diffusion process 
satisfying the stochastic differential equation \vspace{-3mm}
\bea\label{201905262213}
\left\{\begin{array}{ccl}
dX_t &=& A(Z_t,\theta_2)dt+B(Z_t,\theta_1)dw_t
\y
dY_t&=&H(Z_t,\theta_3)dt,
\end{array}\right.
\eea
where $A:\bbR^{\sfd_Z}\times\overline{\Theta}_2\to\bbR^{\sfd_X}$, 
$B:\bbR^{\sfd_Z}\times\overline{\Theta}_1\to\bbR^{\sfd_X}\otimes\bbR^\sfr$ and 
$H:\bbR^{\sfd_Z}\times\overline{\Theta}_3\to\bbR^{\sfd_Y}$. 
The process  $Z_t=(X_t,Y_t)_{t\in\bbR_+}$, $\bbR_+=[0,\infty)$,  is assumed to be 
driven by 
an $\sfr$-dimensional $\F$-Wiener process $w=(w_t)_{t\in\bbR_+}$ 
on a stochastic basis $(\Omega,\calf,\F,P)$ equipped with a right-continuous filtration $\F=(\calf_t)_{t\in\bbR_+}$.   

The estimation theory for the degenerate diffusion process (\ref{201905262213}) is motivated by 
the growing interest in hypo-elliptic diffusions, e.g., 
the harmonic oscillator, the Van der Pol oscillator and 
the FitzHugh-Nagumo neuronal model. 
Estimation is studied in 
Gloter \cite{gloter2006parameter}, 
Samson and Thieullen \cite{samson2012contrast}, 
Comte et al. \cite{comte2017adaptive}, 
Ditlevsen and Samson \cite{ditlevsen2019hypoelliptic}, 
Melnykova \cite{melnykova2019parametric}, 
Delattre et al. \cite{delattre2020rate}, 
Le\'on and Samson \cite{leon2018hypoelliptic}, 
Gloter and Yoshida \cite{gloter2020adaptive,gloter2021adaptive,gloter2024nonadaptive} and 
Iguchi et al. \cite{iguchi2022numerical}. 

Bayesian methods are playing a key role in the modern statistics, synchronizing with the developments of 
computational environments. 
Based on the likelihood ratio, 
Ibragimov and Khasminskii 
\cite{IbragimovHascprimeminskiui1972,IbragimovHascprimeminskiui1973,IbragimovHascprimeminskiui1981} 
gave a new perspective on the asymptotic decision theory. 
Their theory features a large deviation inequality for the likelihood ratio random field, and thanks to it, 
they proved the convergence of moments and asymptotic behavior of the Bayesian estimator. 
The developments in this direction toward semimartingales were by 
Kutoyants
\cite{kutoyants1984parameter,Kutoyants1994,Kutoyants2004,kutoyants2012statistical}. 
The quasi-likelihood analysis (QLA, Yoshida \cite{yoshida2011polynomial, yoshida2021simplified}) 
embodied their program for sampled nonlinear stochastic processes, with a polynomial type large deviation inequality proved for the asymptotically quadratic quasi-likelihood ratio random field. 
The QLA makes it possible to deduce the limit distribution and the moment convergence of 
the quasi-maximum likelihood estimator (QMLE) and 
the quasi-Bayesian estimator (QBE), by simple $L^p$-estimates of the functionals associated with the quasi-likelihood random field. 
Due to the universality, 
the QLA has many applications, e.g.,  
Yoshida \cite{yoshida2011polynomial}, 
Uchida and Yoshida \cite{uchida2012adaptive, uchida2014adaptive, UchidaYoshida2013}, 
Ogihara and Yoshida \cite{ogihara2014quasi, OgiharaYoshida2011}, 
Masuda \cite{masuda2013convergence}, 
Eguchi and Masuda \cite{eguchi2018schwarz}, 
Umezu et al. \cite{umezu2019aic}, 
Masuda and Shimizu \cite{masuda2017moment}, 
Suzuki and Yoshida \cite{suzuki2020penalized}, 
Yoshida \cite{yoshida2018partial}, 
Clinet and Yoshida \cite{clinet2017statistical}, 
Inatsugu and Yoshida \cite{InatsuguYoshida2020Accepted}), 
Muni Toke and Yoshida \cite{muni2017modelling,muni2019analyzing,muni2022marked}, 
among many others. 
Kamatani and Uchida \cite{KamataniUchida2014} proposed a hybrid multi-step estimators for sampled diffusion processes, 
and suggested that the use of a quasi-Bayesian estimator as the initial estimator yields stable estimation avoiding 
the quasi-maximum likelihood estimator apt to be trapped by local maxima. 
The importance of the moment convergence of the estimator is established because 
it is indispensable for making the theory in various fields in statistics, e.g., 
asymptotic decision theory, information criteria for model selection, 
prediction, asymptotic expansion, higher-order inference, bootstrap and resampling plans.

The aim of this paper is to construct the adaptive QLA for the degenerate diffusion model (\ref{201905262213}). 
This enables us to formulate various (at least 16 patterns of) adaptive schemes, especially, the schemes incorporating the quasi-Bayesian estimators 
into them. 
With the QLA, we show the moment convergence of the QLA estimators (i.e., QMLE and QBE) 
as well as the asymptotic normality. 

This paper is organized as follows. 
In Section \ref{202306160554}, we explain how to construct the adaptive schemes for the degenerate diffusion 
model (\ref{201905262213}), and give a summary of the main result. 
Precise treatment of the problem starts from Section \ref{202001141611}. 
In Section \ref{202007281646}, the adaptive QMLE and QBE for $\theta_3$ are investigated. 
These estimators assume initial estimators for $\theta_1$ and $\theta_2$. 
Estimation of $\theta_1$ and $\theta_2$ is possible based on the non-degenerate component of the model (\ref{201905262213}). 
Since the mode of convergence of estimators is stronger than the classical theory, 
we apply the QLA theory in Section \ref{202007281718} to the initial estimators for $\theta_1$ and $\theta_2$. 
Section \ref{2206241735} shows that the initial estimator for $\theta_1$ can be improved 
by using the degenerate component of the model (\ref{201905262213}) to get more information from the data. 
The constructed estimators above are integrated into a joint estimator in Section \ref{202306160251}. 
{\xred 
Section \ref{202402141230} provides the results of numerical simulations. 
Sections \ref{202101221241}-\ref{202001141618} are devoted to the proofs of the main results. 
}

\begin{en-text}
The spaces $\Theta_i$ ($i=1,2,3$) are the unknown parameter spaces of 
{\cred the components of} 
$\theta=(\theta_1,\theta_2,\theta_3)$ to be estimated 
from the data $(Z_{\tj})_{j=0,1,...,n}$, {\cblue where} $\tj=jh$, $h=h_n$ satisfying 
$h\to0$, $nh\to\infty$ and $nh^2\to0$ as $n\to\infty$.

Estimation theory 
has been well developed for {\cred diffusion processes}. 
Even focusing {\cblue on} parametric estimation for ergodic diffusions, there is huge amount of studies: 
Kutoyants \cite{Kutoyants1984,Kutoyants2004,Kutoyants1997}, 
Prakasa Rao \cite{PrakasaRao1983, PrakasaRao1988},
Yoshida \cite{Yoshida1990,Yoshida1992b,yoshida2011polynomial},
Bibby and S{\o}rensen \cite{BibbySoerensen1995},
Kessler \cite{Kessler1997},
K{\"u}chler and S{\o}rensen \cite{kuchler1997exponential},
{\cblue Genon--Catalot et al. \cite{genoncatalot1999parameter},}
Gloter \cite{gloter2000discrete,gloter2001parameter,gloter2007efficient}, 
Sakamoto and Yoshida \cite{SakamotoYoshida2009}, 
Uchida \cite{uchida2010contrast}, 
Uchida and Yoshida \cite{UchidaYoshida2001,uchida2011estimation, uchida2012adaptive}, 
Kamatani and Uchida \cite{KamataniUchida2014}, 
De Gregorio and Iacus \cite{de2012adaptive}, 
{\cblue Genon--Catalot and Lar\'edo \cite{genoncatalot2016estimation},}
Suzuki and Yoshida \cite{suzuki2020penalized} 
among many others. 
{\cred Nakakita and Uchida \cite{nakakita2019inference} and 
Nakakita et al. \cite{nakakita2020quasi} studied estimation under measurement error; 
related are 
Gloter and Jacod \cite{gloter2001diffusions1,gloter2001diffusions}. 
}
\noindent
{\cblue Non parametric estimation for the coefficients of an ergodic diffusion has also been widely studied : Dalayan and Kutoyants \cite{dalalyan2002asymptotically}, 
Kutoyants \cite{Kutoyants2004}, 
Dalalyan \cite{dalalyan2005sharp}, 
Dalalyan and Rei{ss} \cite{dalalyan2006asymptotic, dalalyan2007asymptotic}, 
Comte and Genon--Catalot \cite{comte2006penalized}, 
Comte et al. \cite{comte2007penalized}, 
Schmisser \cite{schmisser2013penalized} to name a few.
Historically attentions were paid to inference for non-degenerate cases.  

Recently there is 

}	
	

In this paper, we will present several estimation schemes. 
Since we assume discrete-time observations of $Z{\cred=(Z_t)_{t\in\bbR+}}$, 
quasi-likelihood estimation for $\theta_1$ and $\theta_2$ {\cred is} known; 
only difference from the standard diffusion case is the existence of 
the covariate $Y{\cred=(Y_t)_{t\in\bbR+}}$ 
in the equation of $X{\cred=(X_t)_{t\in\bbR+}}$ but it causes no theoretical difficulty. 
\end{en-text}
\begin{en-text}
We will give an exposition for construction of those standard estimators 
in Sections \ref{201905291607} and \ref{201906041938} for selfcontainedness. 
\end{en-text}
\begin{en-text}
Thus, our first approach in Section {\tred \ref{202007281646}} 
is toward estimation of $\theta_3$ with initial estimators 
for $\theta_1$ and $\theta_2$. 
The idea for construction of the quasi-likelihood function in the elliptic case was 
based on the local Gaussian approximation of the transition density. 
Then it is natural to approximate the distribution of the increments of $Y$ 
by {\cred that of} the principal Gaussian variable in the expansion of {\cred the increment}. 
However, this method causes deficiency, as we will observe there; {\cred see} 
{\tred Section} \ref{201906181743}. 
We present a more efficient method by incorporating an additional Gaussian part from $X$. 
The {\tred error rate} 
attained by the estimator for $\theta_3$ is $n^{-1/2}h^{1/2}$ and it is 
much faster than the rate $(nh)^{-1/2}$ for $\theta_2$ and $n^{-1/2}$ for $\theta_1$. 
Section \ref{202001141623} treats some adaptive estimators 
using suitable initial estimators for $(\theta_1,\theta_2,\theta_3)$, and 
shows joint {\cred asymptotic normality}. 
Then it should be remarked that the asymptotic variance of our estimator $\wh{\theta}_1$ 
for $\theta_1$ has improved {\cred that of} the ordinary volatility parameter estimator, 
e.g. $\wh{\theta}_1^0$ recalled in Section {\tred\ref{202007281718}} 
that would be 
asymptotically optimal if the system consisted only of $X$. 
\end{en-text}
\begin{en-text}
In Section \ref{202001141632}, 
we consider a non-adaptive joint quasi-maximum likelihood estimator. 
This method does not require initial estimators. 
From computational point of view, adaptive methods often have merits 
by reducing dimension of parameters, 
but 
the non-adaptive method is still theoretically interesting. 
\end{en-text}
\begin{en-text}
{\cred
Section \ref{202001141611} collects the assumptions under which we will work. 
Section \ref{202001141618} offers several basic estimates to the increments of $Z$.  
}

To investigate efficiency of the presented estimators, {\cred we need} 
the LAN property of the exact likelihood function of the hypo-elliptic diffusion. 
{\tred Another important and natural question the reader must have is the 
asymptotic behavior of the joint quasi-maximum likelihood estimator 
based on a quasi-likelihood random field for the full parameter $\theta$; 
an expression of the random field 
has already appeared in (\ref{201906060545}) essentially. 
In the present situation, the three parameters have different convergence rates and in particular 
the handling of $\theta_3$ is not straightforward because 
for estimation of $\theta_3$, the parameters $(\theta_2,\theta_3)$ become 
nuisance, but any estimator of them has very large error compared with $\theta_3$. 
The user could get some estimated value with the joint quasi-likelihood random field, 
however, there is no theoretical backing for such a scheme. 
Though somewhat sophisticated treatments are necessary, 
we can validate the joint quasi-maximum likelihood estimator and 
can show that the same asymptotic variance is attained, 
up to the first order, as the one-step quasi-likelihood estimator 
provided in this article. 
We will discuss these problems elsewhere, while 
we recommend the reader to see Gloter and Yoshida \cite{gloter2020adaptive} 
for more complete exposition including the non-adaptive approach 
and additional information. 
 }
\end{en-text}

{\vblue 
\section{Adaptive schemes and their asymptotic properties}\label{202306160554}
{\vredd 
In our adaptive method, first we make an estimator $\wh{\theta}_3^{\sfA_3}$ for $\theta_3$ by using some initial estimators $\wh{\theta}_1^{0,\sfA_0}$ and $\wh{\theta}_2^{0,\sfA_2}$ for $\theta_1$ and $\theta_2$, respectively, and next an estimator $\wh{\theta}_1^{\sfA_1}$ to improve 
the initial estimator $\wh{\theta}_1^{0,\sfA_0}$ 
with the help of $\wh{\theta}_3^{\sfA_3}$. 
Here $\sfA_i$ ($i=0,1,2,3$) denote $M$ (QMLE) or $B$ (QBE). 
The initial estimators $\wh{\theta}_1^{0,\sfA_0}$ and $\wh{\theta}_2^{0,\sfA_2}$  
are constructed by using the non-degenerate component of the model (\ref{201905262213}).
}

More precisely, 
the estimator $\wh{\vartheta}=\big(\wh{\theta}_1^{A_1},\wh{\theta}_2^{A_2},\wh{\theta}_3^{A_3}\big)$ for the vector parameter $\theta=(\theta_1,\theta_2,\theta_3)$ is adaptively constructed as follows. 
\begin{enumerate}[(Step 1)]
\im $\wh{\theta}_1^{0,\sfA_0}$ for $\sfA_0=M/B$. \\ 
Construct an initial estimator $\wh{\theta}_1^{0,\sfA_0}$ for $\theta_1$ 
based on a quasi-log likelihood function $\H_n^{(1)}(\theta_1)$. 
Here $A_0$ is either $M$ or $B$.
$\wh{\theta}_1^{0,M}$ and $\wh{\theta}_1^{0,B}$ denote 
the QMLE and QBE, respectively, 
with respect to $\H_n^{(1)}(\theta_1)$. $[$Section \ref{202007281718}$]$
\label{202306160300}

\im $\wh{\theta}_2^{0,\sfA_2}$ for $A_2=M/B$. \\
Based on a quasi-log likelihood function $\H_n^{(2)}(\theta_2)$, 
construct an estimator $\wh{\theta}_2^{0,\sfA_2}$ by using $\wh{\theta}_1^{0,\sfA_0}$ in Step \ref{202306160300}, 
where $\sfA_2=M$ (QMLE) or $\sfA_2=B$ (QBE). $[$Section \ref{202007281718}$]$
\label{202306160235}

\im $\wh{\theta}_3^{\sfA_3}$ for $A_3=M/B$. \\
For $\theta_3$, make either QMLE $\wh{\theta}_3^M$ or QBE $\wh{\theta}_3^B$ 
with a quasi-log likelihood function $\bbH_n^{(3)}(\theta_3)$ 
defined by using $\big(\wh{\theta}_1^{0,\sfA_0},\wh{\theta}_2^{0,\sfA_2}\big)$ obtained in Steps \ref{202306160300} and \ref{202306160235}. 
$[$Section \ref{202007281646}$]$
\label{202306160335}

\im $\wh{\theta}_1^{A_1}$ for $A_1=M/B$.\\
For the second-stage estimation for $\theta_1$, 
make either QMLE $\wh{\theta}_1^M$ or QBE $\wh{\theta}_1^B$ 
with a quasi-log likelihood function $\bbH_n^{(1)}(\theta_1)$ 
defined by using $\big(\wh{\theta}_2^{0,\sfA_2},\wh{\theta}_3^{\sfA_3}\big)$. 
Remark that $\sfA_1$ may be different from $\sfA_0$. 
The improved estimator $\wh{\theta}_1^{A_1}$ exploits the estimator $\wh{\theta}_3^{\sfA_3}$ 
as the initial estimator $\wh{\theta}_3^{0}$ for $\theta_3$. 
$[$Section \ref{2206241735}$]$ 
\label{202306160356}

\im $\wh{\theta}_2^{\sfA_2}=\wh{\theta}_2^{0,\sfA_2}$ for $A_2=M/B$. \\
Set $\wh{\theta}_2^{0,\sfA_2}$ in Step \ref{202306160235} as the final estimator $\wh{\theta}_2^{\sfA_2}$ of $\theta_2$. 
\label{202306160423}
\end{enumerate}

The main theorem (Theorem \ref{2206290542}) in Section \ref{202306160251} establishes 
joint asymptotic normality and moment convergence of 
$\wh{\vartheta}=\big(\wh{\theta}_1^{A_1},\wh{\theta}_2^{A_2},\wh{\theta}_3^{A_3}\big)$, that is, 
\beas 
E\big[f\big(a_n^{-1}(\wh{\vartheta}-\theta^*)\big)\big]
&\to&
\bbE[f\big({\sf Z}_1,{\sf Z}_2,{\sf Z}_3\big)]
\eeas
as $n\to\infty$, for every $f\in C_p(\bbR^\sfp)$, where 
$\big({\sf Z}_1,{\sf Z}_2,{\sf Z}_3\big)\sim N_{\sfp_1}(0,\Gamma_{11}^{-1})\otimes N_{\sfp_2}(0,\Gamma_{22}^{-1})\otimes N_{\sfp_3}(0,\Gamma_{33}^{-1})$, and 
\bea\label{202306160306}
a_n &=& \left( \begin{array}{ccc}
 n^{-1/2}{\tt I}_{\sfp_1}&0&0\\
 0&{\colorr n^{-1/2}}h^{{\colorr-1/2}}{\tt I}_{\sfp_2}&0\\
 0&0&{\colorr n^{-1/2}}h^{{\colorr1/2}}{\tt I}_{\sfp_3}
\end{array} \right).
\eea
Here, for $\sfp=\sfp_1+\sfp_2+\sfp_3$, $C_p(\bbR^\sfp)$ denotes the set of $\bbR$-valued continuous functions on $\bbR^\sfp$ 
of at most polynomial growth, and 
each ${\tt I}_{\sfp_i}$ denotes the identity matrix of dimension $\sfp_i$, $i\in\{1,2,3\}$. 
The matrices $\Gamma_{ii}$ ($i=1,2,3$) will be specified in the later sections.

The estimator $\wh{\theta}_1^{0,\sfA_0}$ $[$ resp. $\wh{\theta}_2^{0,\sfA_2}$ $]$ 
in Step \ref{202306160300} $[$ resp. \ref{202306160235} $]$ serves as an initial estimator for $\theta_1$ 
$[$ resp. $\theta_2$ $]$ when one constructs the estimator $\wh{\theta}_3^{\sfA_3}$ in Step \ref{202306160335}. 
Though they are specified in Steps \ref{202306160300} and \ref{202306160235}, 
we can use more general initial estimators to construct $\wh{\theta}_3^{\sfA_3}$ 
in Step \ref{202306160335}. 
Indeed, Section \ref{202007281646} adopts general initial estimators 
$\wh{\theta}_1^{0,\sfA_0}$ and $\wh{\theta}_2^{0,\sfA_2}$, and only 
requires them to satisfy a rate of convergence in $L^p$. 
Similarly, Section \ref{2206241735} works with a pair of general initial estimators $\big(\wh{\theta}_2^0,\wh{\theta}_3^0)$ for $(\theta_2,\theta_3)$, 
including the choice of the initial estimators in Step \ref{202306160356}. 
The estimator $\wh{\theta}_1^{A_1}$ in Step \ref{202306160356} improves the asymptotic variance of 
$\wh{\theta}_1^{0,\sfA_0}$ in Step \ref{202306160300}, as proved in Section \ref{2206241735}. 
As a matter of fact, 
the estimators $\wh{\theta}_1^{A_1}$ and $\wh{\theta}_3^{A_3}$ do not need to share the same initial estimator for $\theta_2$. 
For the final estimator for $\theta_2$ in Step \ref{202306160423}, it is also possible to take $\wh{\theta}_2^{0,\sfA_4}$ for $\sfA_4$ 
such that $\sfA_4\in\{M,B\}\setminus\{\sfA_2\}$, if additional computation does not matter. 
Section \ref{202007281718} proves the asymptotic properties of 
$\wh{\theta}_1^{0,\sfA_0}$ and $\wh{\theta}_2^{0,\sfA_2}$ 
in Steps \ref{202306160300} and \ref{202306160235}, respectively. 
The QLA approach is necessary there since we need $L^p$-estimates of the error to go to the further steps. 
Finally, Theorem \ref{2206290542} in Section \ref{202306160251} integrates 
the asymptotic properties of the estimators in Steps \ref{202306160300}-\ref{202306160423}. 
Let us detail the setting of our investigation. 

}

\section{Setting and assumptions}\label{202001141611}
{\vblue 
In this section, we state the assumptions we will work with to prove asymptotic properties of the estimators. 
The sets 
$\Theta_i$ ($i=1,2,3$) are bounded open domain in $\bbR^{\sfp_i}$, respectively. 
The total parameter space $\Theta=\prod_{i=1}^3\Theta_i$ is suppose to have 
a good boundary so that Sobolev's embedding inequality holds valid: 
there exists a positive constant $C_{{\colorg\Theta}}$ such that 
\bea\label{202001131246}
\sup_{\theta\in\Theta}|{\tred f(\theta)}|&\leq& C_{\Theta}\sum_{k=0}^1 
\|\partial_\theta^kf\|_{L^p(\Theta)}
\eea
for all $f\in C^1(\Theta)$ and $p>\sum_{i=1}^3\sfp_i$, 
where $L^p(\Theta)$ is the $L^p$-space with respect to the Lebesgue measure on $\Theta$. 
If $\Theta$ has a Lipschitz boundary, then this condition is satisfied ({\tred cf. Adams \cite{Adams1975}}).
Obviously, the embedding inequality (\ref{202001131246}) is valid for functions depending only on a part of components of $\theta$.
The use of Sobolev's embedding is an easy way to asses the maximum of a random field on $\Theta$ though we could adopt 
other embedding inequalities such as the GRR inequality to relax the assumptions on 
differentiability of functions. 

For a finite-dimensional real vector space ${\sf E}$, 
we denote by $C^{a,b}_p(\bbR^{\sfd_Z}\times\Theta_i;{\sf E})$ the set of functions $f:\bbR^{\sfd_Z}\times\Theta_i\to{\sf E}$ 
such that $f$ is continuously differentiable $a$ times in $z\in\bbR^{\sfd_Z}$ and $b$ times in $\theta_i\in{\tred \Theta_i}$ 
in any order, and $f$ and all such derivatives are continuously extended to $\bbR^{\sfd_Z}\times\overline{\Theta}_i$, 
moreover, they are of at most polynomial growth in $z\in\bbR^{\sfd_Z}$ uniformly in ${\tred\theta_i\in\Theta_i}$. 

The true value of $\theta=(\theta_1,\theta_2,\theta_3)\in\Theta_1\times\Theta_2\times\Theta_3$ is denoted by $\theta^*=(\theta_1^*,\theta_2^*,\theta_3^*)$, 
that is, the process $(Z_t)_{t\in\bbR_+}$ generating the data satisfies the stochastic differential equation (\ref{201905262213}) for 
$\theta=\theta^*$. 
Denote $C=BB^\star$, where $\star$ denotes the matrix transpose. 
We will denote 
$F_x$ for $\partial_xF$, $F_y$ for $\partial_yF$, $F_z$ for $\partial_zF$ and $F_i$ for $\partial_{\theta_i}F$ for a function $F$ of $(z,\theta)$. 

To state the results, depending on the situation, we will assume some of the following conditions. 
\bd
\im[[A1\!\!]]
{\bf (i)} 
$A\in C^{i_A,j_A}_p(\bbR^{\sfd_Z}\times\Theta_2;\bbR^{\sfd_X})$ and 
$B\in C^{i_B,j_B}_p(\bbR^{\sfd_Z}\times\Theta_1;\bbR^{\sfd_X}\otimes\bbR^\sfr)$. 

\bd\im[(ii)] $H\in C^{i_H,j_H}_p(\bbR^{\sfd_Z}\times\Theta_3;\bbR^{\sfd_Y})$. 
\ed
\ed

\bd
\im[[A2\!\!]] {\bf (i)} 
{\xred$\|Z_0\|_p<\infty$} 
for every $p>1$. 
\bd\im[(ii)] The process $Z$ is a stationary process with an invariant probability measure $\nu$ on $\bbR^{\sfd_Z}$ and satisfies the mixing condition 
that 
\beas 
\alpha_Z(s) &\leq& a^{-1}e^{-as}\quad(s>0)
\eeas
for some positive constant $a$, where 
\beas 
\alpha_Z(s) &=& \sup_{t\in\bbR_+}\sup\bigg\{\big|P[K\cap M]-P[K]P[M]\big|;\>
K\in\sigma[Z_r;r\leq t],\>M\in\sigma[Z_r;r\geq t+s]\bigg\}. 
\eeas

\im[(iii)] The function $\theta_1\mapsto C(Z_t,\theta_1)^{-1}$ is continuous 
on $\overline{\Theta}_1$
a.s., and
\beas
\sup_{\theta_1\in\Theta_1}\sup_{t\in\bbR_+}\|\det C(Z_t,\theta_1)^{-1}\|_p&<&\infty
\eeas 
for every $p>1$. 

\im[(iv)] For the $\bbR^{\sfd_Y}\otimes\bbR^{\sfd_Y}$ valued function 
\bea\label{202008020107}
V(z,\theta_1,\theta_3)=H_x(z,\theta_3)C(z,\theta_1)H_x(z,\theta_3)^\star, 
\eea
the function $(\theta_1,\theta_3)\mapsto V(Z_t,\theta_1,\theta_3)^{-1}$ is continuous 
on $\overline{\Theta}_1\times\overline{\Theta}_3$ a.s., 
and 
\beas 
\sup_{(\theta_1,\theta_3)\in\Theta_1\times\Theta_3}\sup_{t\in\bbR_+}\|\det V(Z_t,\theta_1,\theta_3)^{-1}\|_p&<&\infty
\eeas 
for every $p>1$. 
\ed
\ed
\halflineskip

Condition $[A2]$ (iii) and (iv) implicitly assume 
the existence of $C(Z_{{\tred t}},\theta_1)^{-1}$ and $V(Z_t,\theta_1,\theta_3)^{-1}$. 
%
Since $Z$ is stationary with the invariant measure $\nu$, we have 
\beas 
\int |z|^p\nu(dz)+\sup_{\theta_1\in\overline{\Theta}_1}\int \big(\det C(z,\theta_1)\big)^{-p}\nu(dz)
+\sup_{(\theta_1,\theta_3)\in\overline{\Theta}_1\times\overline{\Theta}_3}
\int \big(\det V(z,\theta_1,\theta_3)\big)^{-p}\nu(dz)<\infty
\eeas 
for any $p>0$. 
%

Exponential ergodicity and boundedness of any order of moment of the process are known. 
We refer the reader to Wu \cite{wu2001large} for damping Hamiltonian systems, 
The Lyapounov function method provides exponential mixing (even in the non-stationary case) and 
estimates of moments of the invariant probability measure up to any order. 
Several examples including the van der Pol model are investigated therein. 
Additional information is in Delattre et al. \cite{delattre2020rate}. 
For nondegenerate diffusions, see e.g. 
Pardoux and Veretennikov \cite{pardoux2001poisson}, 
Meyn and Tweedie \cite{MeynTweedie1993a} and 
Kusuoka and Yoshida \cite{KusuokaYoshida2000} among many others. 
\begin{en-text}
It is well known that the ergodic property 
is deeply related to the regularity of the transition probability and 
positive recurrence of the process. Then the H\"{o}rmander condition and 
the drift condition are of practical use. 
The Malliavin calculus is a common tool for the regularity of the transition probability. 
Smoothness of the transition density is nothing but a smoothness of the 
Watanabe's delta function (i.e., the mapping from the location to 
a generalized Wiener functional dwelling in 
a Sobolev space with a negative differentiability index on the Wiener space). 
Ikeda and Watanabe \cite{ikeda2014stochastic} is a standard textbook. 
\end{en-text}
\begin{en-text}
, which uses 
Watanabe's delta function (i.e. a delta function characterized as a generalized 
Wiener functional 
in a Sobolev space with a negative differentiability index on the Wiener space) 
to show it. 
\end{en-text}
\begin{en-text}
The use of the drift condition is also standard. 
We recommend the reader should consult the textbook or a series of papers 
by Meyn and Tweedie. 
Wu \cite{wu2001large} 
provided exponential convergence for stochastic damping Hamiltonian systems. 
\end{en-text}
\begin{en-text}
Our model is nonlinear and therefore nonlinear functionals (e.g. polynomial) of the data $Z$ appear in statistical inference. 
The collection of appearing random variables should be an algebra, i.e., closed in $L^\inftym=\cap_{p>1}L^p$. 
Usually, results on existence of a solution of the SDE entail $L^\inftym$-estimates. 
\end{en-text}
\begin{en-text}
As it is quite well known especially for nondegenerate diffusions, this condition is generally satisfied, for example, see 
E. Pardoux, A.Y. Veretennikov 
On the Poisson equation and diffusion approximation 1
Ann. Probab., 29 (2001), pp. 1061-1085. 
For degenerate diffusions, Wu \cite{wu2001large} solved this question. 
The new version added a remark in (new) Section 2. 
\end{en-text}
\halflineskip
}


%

\begin{en-text}
\bi
\im $\sup_{\theta\in\Theta}\|R_j(\theta)\|_p=O(h^2)$ for every $p>1$. 

\im Conditionally on $\calf_\tjm$, the distribution of $\Delta_jY$ is approximated by 
the normal distribution\\ $N\big(hF_n(Z_\tjm,\theta_3),h^3K(Z_\tjm,(\theta_1,\theta_3)\big)$. 

\im This suggests the use of the quasi-log likelihood function 
\beas 
\bbH^{[3]}_n(\theta_3)
 &=& 
 -\half\sum_{j=1}^n \bigg\{
 \frac{\big(\Delta_jY-hF_n(Z_\tjm,\theta_3)\big)^2}{h^3K(Z_\tjm,\wh{\theta}_1,\theta_3)}
 +\log\big(h^3K(Z_\tjm,\wh{\theta}_1,\theta_3)\big)\bigg\}
 \eeas
where 
$\wh{\theta}_1=\wh{\theta}_{1,n}$ and $\wh{\theta}_2^0=\wh{\theta}_{2,n}$ implicitly in $F_n$ 
are consistent estimators for $\theta_1$ and $\theta_2$ respectively based on the data $(Z_\tj)_{j=0,1,...,n}$. 
\im $B$ is non-degenerate. 
\im Observe $(Z_\tj)_{j=0,1,...,n}$, $\tj=t_j^n=jh$, $h=h_n=\tj-\tjm$
\im Assume $h\to0$, $nh\to\infty$, $nh^2\to0$. 
\im Approximation to $\Delta_jY=Y_\tj-Y_\tjm$
\ei
\end{en-text}

We need the following notation. 
Denote $U^{\otimes k}$ for $U\otimes \cdots\otimes U$ ($k$-times) for a tensor $U$. 
For tensors  
$S^1=(S^1_{i_{1,1},...,i_{1,d_1};j_{1,1},...,j_{1,k_1}})$, ..., 
$S^m=(S^m_{i_{m,1},...,i_{m,d_m};j_{m,1},...,j_{m,k_m}})$ 
and a tensor $T=(T^{i_{1,1},...,i_{1,d_1},...,i_{m,1},...,i_{m,d_m}})$, we write 
\beas &&
T[S^1,...,S^m]
\yeq
T[S^1\otimes\cdots\otimes S^m]
\\&=&
\bigg(
\sum_{i_{1,1},...,i_{1,d_1},...,i_{m,1},...,i_{m,d_m}}
T^{i_{1,1},...,i_{1,d_1},...,i_{m,1},...,i_{m,d_m}}
S^1_{i_{1,1},...,i_{1,d_1};j_{1,1},...,j_{1,k_1}}
\\&&\hspace{3cm}\cdots S^m_{i_{m,1},...,i_{m,d_m};j_{m,1},...,j_{m,k_m}}\bigg)
_{j_{1,1},...,j_{1,k_1},...,j_{m,1},...,j_{m,k_m}}.
\eeas
This notation will be applied for a tensor-valued tensor $T$ as well. 
For example, $H_{xx}$ is a tensor with $\sfd_Y\times\sfd_X^2$ components $(H_{xx})^i_{j,k}$ ($i=1,...,\sfd_Y;\> j,k=1,...,\sfd_X$), and 
$H_{xx}[C(z,\theta_1)]$ is a $\sfd_Y$-dimensional vector with $i$-th component 
$
\sum_{j,k=1}^{\sfd_X}(H_{xx})^i_{j,k}C(z,\theta_1)^{j,k}
$ 
for $C(z,\theta_1)=\big(C(z,\theta_1)^{j,k}\big)_{j,k=1,...,\sfd_X}$.

{\vblue 
There are two ways to estimate $\theta_1$. 
One is the classical approach applied to the non-degenerate component of the stochastic differential equation (\ref{201905262213}) and 
the other is the new approach that incorporates the information from the degenerate component. 
It will turn out that the latter is more efficient than the former though the classical one serves as an initial estimator of $\theta_1$. 
Thus, there are two identifiability conditions to consider. 

Define the function $\bbY^{(1)}$ by 
\beas 
\bbY^{(1)}(\theta_1) 
&=& 
-\half\int \bigg\{
\text{Tr}\big(C(z,\theta_1)^{-1}C(z,\theta_1^*)\big)-\sfd_X
+\log\frac{\det C(z,\theta_1)}{\det C(z,\theta_1^*)}\bigg\}\nu(dz).
\nn\\&&
\eeas
Since 
$|\log x|\leq x+x^{-1}$ for $x>0$, $\bbY^{(1)}(\theta_1)$ is a continuous function on $\overline{\Theta}_1$ 
well defined 
under {\colorg $[A1]$ and }$[A2]$. 
We remark that the function $\bbY^{(1)}(\theta_1)$ is different from the one denoted by the same symbol in 
Gloter and Yoshida \cite{gloter2020adaptive,gloter2021adaptive}. 
}
{\vblue 
Let 
\beas 
\bbY^{(J,1)}(\theta_1)
&=& 
-\half\int \bigg\{\text{Tr}\big(C(z,\theta_1)^{-1}C(z,\theta_1^*)\big)
+\text{Tr}\big(V(z,\theta_1,\theta_3^*)^{-1}V(z,\theta_1^*,\theta_3^*)\big)
-\sfd_Z
\\&&
+\log \frac{\det C(z,\theta_1)\det V(\theta_1,\theta_3^*)}{\det C(z,\theta_1^*)\det V(\theta_1^*,\theta_3^*)}
\bigg\}
\nu(dz).
\eeas

Let 
\beas
S(z,\theta_1,\theta_3)&=&
\left(\begin{array}{ccc}
C(z,\theta_1)&& 2^{-1}C(z,\theta_1)H_x(z,\theta_3)^\star
\\
2^{-1}H_x(z,\theta_3)C(z,\theta_1)&& 3^{-1}H_x(z,\theta_3)C(z,\theta_1)H_x(z,\theta_3)^\star
\end{array}\right).
\eeas
Then 
\begin{en-text}
\bea\label{201906040958} &&
S(z,\theta_1,\theta_3)^{-1}
\nn\\&=&
\left(\begin{array}{cc}
\left\{\begin{array}{c}
C(z,\theta_1)^{-1}\\
+3H_x(z,\theta_3)^\star V(z,\theta_1,\theta_3)^{-1}H_x(z,\theta_3)
\end{array}\right\}
&-6H_x(z,\theta_3)^\star V(z,\theta_1,\theta_3)^{-1}\y
-6V(z,\theta_1,\theta_3)^{-1}H_x(z,\theta_3)&12V(z,\theta_1,\theta_3)^{-1}
\end{array}\right). 
\nn\\&&
\eea
\end{en-text}
\bea\label{201906040958} 
S(z,\theta_1,\theta_3)^{-1}
&=&
\left(\begin{array}{cc}
S(z,\theta_1,\theta_3)^{1,1}&S(z,\theta_1,\theta_3)^{1,2}\y
S(z,\theta_1,\theta_3)^{2,1}&S(z,\theta_1,\theta_3)^{2,2}
\end{array}\right),
\eea
where 
\beas 
S(z,\theta_1,\theta_3)^{1,1}
&=&
C(z,\theta_1)^{-1}
+3H_x(z,\theta_3)^\star V(z,\theta_1,\theta_3)^{-1}H_x(z,\theta_3),
%
\nn\y
S(z,\theta_1,\theta_3)^{1,2}
&=&
-6H_x(z,\theta_3)^\star V(z,\theta_1,\theta_3)^{-1},
%
\nn\y
S(z,\theta_1,\theta_3)^{2,1}
&=&
-6V(z,\theta_1,\theta_3)^{-1}H_x(z,\theta_3),
\nn\y
S(z,\theta_1,\theta_3)^{2,2}
&=& 
12V(z,\theta_1,\theta_3)^{-1}.
\eeas
Recall (\ref{202008020107}): 
$
V(z,\theta_1,\theta_3) = H_x(z,\theta_3)C(z,\theta_1)H_x(z,\theta_3)^\star.
$
Simple calculus gives 
\bea\label{2206251502}
\bbY^{(J,1)}(\theta_1)
&=& 
-\frac{1}{2}\int\bigg\{ 
\big(S(z\theta_1,\theta_3^*)^{-1}-S(z,\theta_1^*,\theta_3^*)^{-1}\big)\big[S(z,\theta_1^*,\theta_3^*)\big]
\nn\\&&\hspace{50pt}
 +\log\frac{\det S(z,\theta_1,\theta_3^*)}{\det S(z,\theta_1^*,\theta_3^*)}
 \bigg\}\nu(dz).
 \eea

For identifiability of $\theta_2$ and $\theta_3$, we will use the following functions, respectively: 
\bea\label{202002171821}
\bbY^{(2)}(\theta_2)
&=&
-\half \int 
C(z,\theta_1^*)^{{\cred -1}}\big[\big(A(z,\theta_2)-A(z,\theta_2^*)\big)^{\otimes2}\big]\nu(dz),
\eea
and
\beas 
\bbY^{(3)}(\theta_3) &=&
-\int 6V(z,\theta_1^*,\theta_3)^{-1}
\big[\big(H(z,\theta_3)-H(z,\theta_3^*)\big)^{\otimes2}\big]\nu(dz). 
\eeas
The both functions are well defined under $[A1]$ and $[A2]$. 
\begin{en-text}
{\tred 
Obviously, $\nu$ depends on the value $\theta^*$. We suppress $\theta^*$ from notation since it is fixed in this article, where it is not necessary to change $\theta^*$ 
differently from discussion of the asymptotic minimax bound for example.  
}
\end{en-text}
}
\begin{en-text}
{\coloro Let 
\beas 
\bbY^{(J,3)}(\theta_1,\theta_3) &=&
-\int 6V(z,\theta_1,\theta_3)^{-1}
\big[\big(H(z,\theta_3)-H(z,\theta_3^*)\big)^{\otimes2}\big]\nu(dz). 
\eeas
}
\end{en-text}

{\vblue 
We will assume all or some of the following identifiability conditions. 
Verification of them is usually resulted in just checking difference of each coefficient function in (\ref{201905262213}) evaluated at $\theta_i$ and $\theta_i^*$. 
\bd
\im[{\bf[A3]}] 
\bd
\im[(i)] There exists a positive constant $\chi_1$ such that 
\beas 
\bbY^{(1)}(\theta_1) &\leq& -\chi_1|\theta_1-\theta_1^*|^2\qquad(\theta_1\in\Theta_1). 
\eeas
\im[(i$'$)] There exists a positive constant $\chi_1'$ such that 
\beas 
\bbY^{(J,1)}(\theta_1)&\leq& -\chi_1'|\theta_1-\theta_1^*|^2\qquad(\theta_1\in\Theta_1). 
\eeas

\im[(ii)] There exists a positive constant $\chi_2$ such that 
\beas 
\bbY^{(2)}(\theta_2) &\leq& -\chi_2|\theta_2-\theta_2^*|^2\qquad(\theta_2\in\Theta_2). 
\eeas

\im[(iii)]
There exists a positive constant $\chi_3$ such that 
\beas 
\bbY^{(3)}(\theta_3) &\leq& -\chi_3|\theta_3-\theta_3^*|^2\qquad(\theta_3\in\Theta_3). 
\eeas

\begin{en-text}
{\coloro
\im[(iii$'$)]
There exists a positive constant $\chi_3$ such that 
\beas 
\bbY^{(J,3)}(\theta_1,\theta_3) &\leq& -\chi_3|\theta_3-\theta_3^*|^2
\qquad(\theta_1\in\Theta_1,\>\theta_3\in\Theta_3). 
\eeas
}
\end{en-text}
\ed
\ed

\begin{en-text}
In the hypoelliptic case, as it is the most interesting case, checking these identifiability conditions is usually easy since $\nu$ is equivalent to or at least 
dominated by the Lebesgue measure and 
admits a density that is positive on a non-empty open set. 
Thus, identifiability is a problem of parameterization of the model. 
In particular, it is obvious that this condition causes no difficulty for 
linearly parametrized models often appearing in applications. 
\end{en-text}

Obviously, Condition $[A3]$ (i) implies Condition $[A3]$ (i$'$).

Conditions 
$h\to0$, $nh\to\infty$ and $nh^2\to0$ as $n\to\infty$ will be assumed in this article.  
The condition $nh^2\to0$ is called 
the condition for rapidly increasing experimental design 
(Prakasa Rao \cite{PrakasaRao1988}). 
In the non-degenerate case, 
Yoshida \cite{Yoshida1992b} relaxed this condition to $nh^3\to0$, and 
Kessler \cite{Kessler1997} to $nh^p\to0$ for any positive number $p$. 
In the frame of the the quasi-likelihood analysis, 
Uchida and Yoshida \cite{uchida2012adaptive} gave adaptive methods under the condition $nh^p\to0$. 
These approaches under $nh^p\to0$ 
need more smoothness of the model than our assumptions  
because they 
inevitably involve higher-order expansions of the semigroup of the process. 
In this article, we do not pursue the problem of relaxation of the condition $nh^2\to0$, but 
recently, Iguchi et al. \cite{iguchi2022numerical} studied estimation under $nh^3\to0$. 
\begin{en-text}
In this paper, when estimating the order of a random variable, 
eventually we use either $n\to\infty$, 
$h\to0$, $nh\to\infty$ or $nh^2\to0$, and that's all. 
So, it is easy for the reader to recognize which convergence is used in each case. 
For example, if the reader finds $O_p(\sqrt{n}h)$, then quite likely it will be estimated 
as $o_p(1)$. However, we left traces as many as possible in the proof. 
\end{en-text}
Additionally to the basic conditions, 
we also assume a weak balance condition that 
$nh\geq n^{\ep_0}$ for large $n$ for some positive constant $\ep_0$. 
This condition is used in Yoshida \cite{yoshida2011polynomial} and other studies. 
In summary, we assume 
\beas 
n\to\infty,\quad nh\to\infty,\quad nh^2\to0,\quad nh\geq n^{\ep_0}.
\eeas
}

{\vblue 
In this paper, we will discuss several methods of adaptive estimation using some initial estimators 
$\wh{\theta}_i^0$ for $\theta_i$, $i=1,2,3$. 
The following standard convergence rates will be assumed in part or fully: 
\bd
\im[{\bf[A4]}] 
\bd
\im[(i)]
$\ds \wh{\theta}_1^0-\theta_1^* \yeq O_{{\vblue L^\inftym}}(n^{-1/2})$ as $n\to\infty$
\im[(ii)] 
$\ds \wh{\theta}_2^0-\theta_2^*\yeq O_{{\vblue L^\inftym}}(n^{-1/2}h^{-1/2})$  as $n\to\infty$
\im[(iii)]
$\ds \wh{\theta}_3^0-\theta_3^*\yeq O_{L^\inftym}(n^{-1/2}h^{1/2})$  as $n\to\infty$. 
\ed\ed
Here and after, we write $V_n=O_{L^\inftym}(b_n)$ [resp. $o_{L^\inftym}(b_n)$] 
for a sequence of random variables/vectors/matrices and a sequence of positive numbers $b_n$, to mean that 
$\||V_n|\|_{L^p(P)}=O(b_n)$ [resp. $\||V_n|\|_{L^p(P)}=o(b_n)$] as $n\to\infty$ for all $p>1$. 
The above condition $[A4]$ is requesting the same rates of convergence for the initial estimators 
as Condition {\vredd$[A4^\sharp]$} of Gloter and Yoshida \cite{gloter2020adaptive, gloter2021adaptive}, 
but stronger 
since we will treat a much stronger mode of convergence, i.e., convergence of moments, than those papers. 
The $L^p$-convergences in $[A4]$ are verifiable if the QLA theory is applied to the initial estimators. 
See Sections \ref{202007281718} and \ref{202007281646}. 
}

{\vblue 
\section{Adaptive quasi-likelihood analysis for $\theta_3$}\label{202007281646}
This section clarifies asymptotic properties of the quasi-maximum likelihood estimator and the quasi-Bayesian estimator for $\theta_3$, 
given initial estimators {\xred$\wh{\theta}_1^0$} and {\xred$\wh{\theta}_2^0$}. 

\subsection{Adaptive QMLE and QBE for $\theta_3$}\label{202007290607}
\begin{en-text}
\begin{remark}\rm 
Clearly, this notation has an advantage over the notation by matrix product since 
the elements $S^1, ..., S^m$ quite often have a long expression 
in the inference. The matrix notation repeats $S^i$s twice for the quadratic form,  
thrice for the cubic form, and so on. 
This notation was introduced by 
\cite{yoshida2011polynomial} and 
already adopted by many papers, e.g., 
\cite{uchida2010contrast}, 
\cite{uchida2012adaptive}, 
\cite{yoshida2013martingale}, 
\cite{masuda2013convergence}, 
\cite{KamataniUchida2014}, 
\cite{masuda2017moment}, 
\cite{eguchi2018schwarz}, 
\cite{nualart2019asymptotic}, 
\cite{nakakita2019inference}, 
\cite{nakakita2020quasi}, 
just to name a few. 
\end{remark}
\end{en-text}

Define the $\bbR^{\sfd_Y}$-valued function $G_n(z,\theta_1,\theta_2,\theta_3)$ by 
\bea\label{20200803} 
G_n(z,\theta_1,\theta_2,\theta_3) 
&=& 
H(z,\theta_3)+\frac{h}{2} L_H\big(z,\theta_1,\theta_2,\theta_3\big),
\eea
where 
\beas
L_H(z,\theta_1,\theta_2,\theta_3)
&=&
H_x(z,\theta_3)[A(z,\theta_2)]+\half H_{xx}(z,\theta_3)[C(z,\theta_1)]
\nn\\&&
+H_y(z,\theta_3)[H(z,\theta_3)].
\eeas
The estimating function will be constructed by using the approximately centered increments 
\bea\label{201906030041}
\cald_j(\theta_1,\theta_2,\theta_3) &=&
\left(\begin{array}{c}
h^{-1/2}\big(\Delta_jX-{\colorr h}A(Z_\tjm,\theta_2)\big)\y
h^{-3/2}\big(\Delta_jY-{\colorr h}G_n(Z_\tjm,\theta_1,\theta_2,\theta_3)\big)
\end{array}\right).
\eea
It is possible to make another estimating function with the second component of (\ref{201906030041}). 
However, this approach gives less optimal asymptotic variance. 
See Section 8 of Gloter and Yoshida \cite{gloter2021adaptive} for details. 
}

\begin{en-text}
We will work with some initial estimators $\wh{\theta}_1^0$ for $\theta_1^0$ and $\wh{\theta}_2^0$ for $\theta_2$. 
The expansions (\ref{201906021458}) and (\ref{201905291809}) 
with Lemma {\colorg\lfive} 
suggest 
two approaches for estimating $\theta_3$. 
The first approach is based on the likelihood of 
{\tred$\Delta_jY$ only, without assistance of $\Delta_jX$.}  
The second one uses the likelihood corresponding to $\cald_j(\theta_1,\theta_2,\theta_3)$. 
However, it is possible to show that the first approach gives less optimal asymptotic variance; 
see {\tred Section} \ref{201906181743}. 
So, we will {\fred take} the second approach here. 
\end{en-text}

{\vblue 
Let 
\bea\label{202008020108}
\kappa({\fred z,\theta_1,\theta_3})
&=&
3^{-1/2}H_x({\fred z,\theta_3})B({\fred z,\theta_1}). 
\eea
%
The principal conditionally Gaussian part of $(\Delta_jX,\Delta_jY)$ has the covariance matrix 
\bea\label{2206241549}
S(z,\theta_1,\theta_3)
&=&
\begin{pmatrix}
B(z,\theta_1)B(z,\theta_1)^\star&\frac{\sqrt{3}}{2}B(z,\theta_1)\kappa(z,\theta_1,\theta_3)^\star\\
\frac{\sqrt{3}}{2}\kappa(z,\theta_1,\theta_3)B(z,\theta_1)^\star&
\kappa(z,\theta_1,\theta_3)\kappa(z,\theta_1,\theta_3)^\star
\end{pmatrix},
\eea
if properly scaled and evaluated at $Z_\tjm=z$ and $(\theta_1,\theta_3)=(\theta_1^*,\theta_3^*)$. 
See Lemmas \lfour 
and \lfive 
{\xred of Gloter and Yoshida \cite{gloter2021adaptive}}
for validation.

We will work with some initial estimators $\wh{\theta}_1^0$ for $\theta_1^0$ and $\wh{\theta}_2^0$ for $\theta_2$. 
Sections \ref{202007281718} 
recalls standard construction of estimators for $\theta_1$ and $\theta_2$, respectively.
With
$
\wh{S}(z,\theta_3)=S(z,\wh{\theta}_1^0,\theta_3), 
$
we define a log quasi-likelihood function $\bbH^{(3)}_n(\theta_3)$ for $\theta_3$ by 
\bea\label{201906041919} 
\bbH^{(3)}_n(\theta_3)
 &=& 
 -\half\sum_{j=1}^n \bigg\{
 \wh{S}(Z_\tjm,\theta_3)^{-1}\big[\cald_j(\wh{\theta}_1^0,\wh{\theta}_2^0,\theta_3)^{\otimes2}\big]
 +\log\det \wh{S}(Z_\tjm,\theta_3)\bigg\}.
\eea
Let {\xred$\wh{\theta}_3^M$} be a quasi-maximum likelihood estimator (QMLE, {\vblue adaptive QMLE}) for $\theta_3$
 for $\bbH^{(3)}_n$, that is, ${\vblue\wh{\theta}_3^M}$ is a 
 $\overline{\Theta}_3$-valued  measurable mapping satisfying 
 \beas 
 \bbH^{(3)}_n({\vblue\wh{\theta}_3^M})
 &=&
 \max_{\theta_3\in\overline{\Theta}_3}\bbH^{(3)}_n(\theta_3). 
 \eeas
 The adaptive QMLE 
 ${\vblue\wh{\theta}_3^M}$ for $ \bbH^{(3)}_n$ 
 depends on $n$ as it does on the data $(Z_\tj)_{j=0,1,...,n}$; 
 $\wh{\theta}_1^0$ in the function $\wh{S}$ also depends on $(Z_\tj)_{j=0,1,...,n}$. 
 
 {\vblue 
 {\vredd 
 The advantage of the frame of the QLA is that it gives asymptotic properties of the quasi-Bayesian estimator as well as the $L^\inftym$-boundedness of the QLA estimators. 
 Let 
 \bea\label{202306172011}
 \wt{\theta}_3^B
 &=& 
 \bigg(\int_{\Theta_3}\exp\big(\bbH^{(3)}_n(\theta_3)\big)\varpi^{(3)}(\theta_3)d\theta_3\bigg)^{-1}
 \int_{\Theta_3}\theta_3\exp\big(\bbH^{(3)}_n(\theta_3)\big)\varpi^{(3)}(\theta_3)d\theta_3, 
 \eea
 where $\varpi^{(3)}$ is a continuous prior density for $\theta_3$ that satisfies 
 $0<\inf_{\Theta_3}\varpi^{(3)}(\theta_3)\leq\sup_{\Theta_3}\varpi^{(3)}(\theta_3)<\infty$. 
The quasi-Bayesian estimator $\wt{\theta}_3^B$ takes values in the convex hull of $\Theta_3$. 
The adaptive quasi-Bayesian estimator (QBE, adaptive QBE, adaBayes) $ \wh{\theta}_3^B$ for $\theta_3$ is 
defined by 
 \bea\label{202306172012}
 \wh{\theta}_3^B
 &=& 
 \left\{\begin{array}{ll}
 \wt{\theta}_3^B&\text{if }\  \wt{\theta}_3^B\in\Theta_3\y
 \theta_3^o&\text{otherwise},
 \end{array}\right.
 \eea
 where $\theta_3^o$ is an arbitrarily prescribed value in $\Theta_3$. 
 The reason of such a modification is to make it possible for the Bayesian estimators to serve as an initial estimator, that will be plugged in the model. 

In what follows, when proving asymptotic properties of a quasi-Bayesian estimator, 
we apply the QLA (Yoshida \cite{yoshida2011polynomial, yoshida2021simplified}) to 
the unmodified version of the Bayesian type estimator $\wt{\theta}_i^B$ for $\theta_i$ like $\wt{\theta}_3^B$ for $\theta_3$. 
Then a conclusion is that there exists a sequence $(\varphi_n)$ of positive numbers 
such that $0<\liminf_{n\to\infty}\log\varphi_n/\log n\leq\limsup_{n\to\infty}\log\varphi_n/\log n<\infty$ and that 
$\big\||\wt{\theta}_i^B-\theta_i^*|\big\|_p=O(\varphi_n^{-1})$ as $n\to\infty$ for all $p>1$. 
Then, for any continuous function $f$ of at most polynomial growth and any constant $L$, we obtain 
\beas 
E\big[f\big(\varphi_n(\wh{\theta}_i^B-\theta_i^*)\big)\big]
&=&
E\big[f\big(\varphi_n(\wt{\theta}_i^B-\theta_i^*)\big)\big]+O(n^{-L})\quad(n\to\infty)
\eeas
since $\Theta_i$ is bounded and the value of $f$ can be controlled by a power of $\varphi_n$. 
Consequently, the asymptotic properties we prove for $\wh{\theta}_i^B$ are the same as those of 
$\wt{\theta}_i^B$. 
In this sense, we will identify $\wh{\theta}_i^B$ with $\wt{\theta}_i^B$ without mentioning the procedure of the modification (\ref{202306172012}) of $\wt{\theta}_i^B$, and simply say $\wh{\theta}_i^B$ is defined 
by the expression on the right-hand side of (\ref{202306172011}). 
If $\Theta_i$ is convex, it is not necessary to consider the modification. 
}
}
}

\subsection{{\vblue Quasi-likelihood analysis for $\wh{\theta}_3^M$ and $\wh{\theta}_3^B$}} 
{\vblue The Fisher information matrix of $\theta_3$ at $\theta_3^*$ is given by }
\bea\label{201906190133} 
\Gamma_{33}
&=&
\int S(z,\theta_1^*,\theta_3^*)^{-1}
\left[
\left(\begin{array}{c}0\\
\partial_3H(z,\theta_3^*)
\end{array}\right)^{\otimes2}
\right]
\nu(dz)
\nn\\&=&
\int 12V(z,\theta_1^*,\theta_3^*)^{-1}\big[\big(\partial_3H(z,\theta_3^{{\colorg*}})\big)^{\otimes2}\big]\nu(dz)
\nn\\&=&
\int12\partial_3H(z,\theta_3^*)^\star V(z,\theta_1^*,\theta_3^*)^{-1}\partial_3H(z,\theta_3^*)\nu(dz).
\eea
{\vblue 
Let 
\bea\label{2206290423}
\zeta_j &=& \sqrt{{\colorr{3}}}\int_\tjm^\tj\int_\tjm^t dw_sdt.
\eea
The quasi-score function of $\theta_3$ at $\theta_3^*$ is 
\bea\label{2206241544}
M^{(3)}_n
&=&
n^{-1/2}\sum_{j=1}^nS(Z_\tjm,\theta_1^*,\theta_3^*)^{-1}
\Bigg[
\left(\begin{array}{c}
h^{-1/2}B(Z_\tjm,\theta_2^*)\Delta_jw\\
h^{-3/2}\kappa(Z_\tjm,\theta_1^*,\theta_3^*)\zeta_j
\end{array}\right),\>
\left(\begin{array}{c}0\\
\partial_3H(Z_\tjm,\theta_3^*)
\end{array}\right)
\Bigg]. 
\nn\\&&
\eea

The following theorem provides asymptotic normality of $\wh{\theta}_3^M$ and $\wh{\theta}_3^B$. 
The convergence of {\vblue $\wh{\theta}_3^M$ and $\wh{\theta}_3^B$} is much faster than other components of 
estimators.
The proof 
is in Section \ref{202101221241}.
Let $\bbi=(i_A,j_A,i_B,j_B,i_H,j_H)$. 
}
\begin{theorem}\label{201906031711}
Suppose that $[A1]$ with {\vblue $\bbi=(1,1,2,1,3,3)$}, 
$[A2]$, 
$[A3]$ $(iii)$ and $[A4]$ {\vredd $(i)$} are satisfied. Then 
\bd
\im[(a)] {\vblue For $\sfA\in\{M,B\}$,}
\beas 
n^{1/2}h^{-1/2}\big(\wh{\theta}_3^\sfA-\theta_3^*\big)
- \Gamma_{33}^{-1}M_n^{(3)}
&\to^p& 0
\eeas
as $n\to\infty$. 
\im[(b)] {\vblue For $\sfA\in\{M,B\}$ and $f\in C_p(\bbR^{\sfp_3})$, 
\beas 
E\big[f\big(n^{1/2}h^{-1/2}\big(\wh{\theta}_3^\sfA-\theta_3^*\big)\big)\big]
&\to&
\bbE[f({\sf Z}_3)]
\eeas
as $n\to\infty$, where  
${\sf Z}_3\sim N(0,\Gamma_{33}^{-1})$. }
\ed
\end{theorem}
{\vblue 
\begin{remark}\rm 
The differentiability $\bbi=(1,1,2,1,3,3)$ requested in Theorem \ref{201906031711} is slightly stronger than 
$\bbi=(1,1,2,1,3,2)$ of Theorem 4.5 of Gloter and Yoshida \cite{gloter2021adaptive} 
and Theorem 3.3 of Gloter and Yoshida \cite{gloter2021adaptive}. 
Condition $[A4]$ here includes the rates of convergence in $L^\inftym$, differently from the above papers, 
while the resulting mode of convergence is much stronger, besides the asymptotic normality of the quasi-Bayesian estimator has been proved. 
\end{remark}
}

{\vblue 
We consider a convergence of the joint distribution of the estimators. 
In the following condition, 
$\wh{\theta}_i$ denotes an estimator for $\theta_i$ based on the data $(Z_\tj)_{j=0,...,n}$, $i\in\{1,2\}$, and 
$\bbT^{(1)}_n$, $\bbT^{(2)}_n$ and ${\mathfrak L}$ denote random vectors 
of dimensions $\sfp_1$, $\sfp_2$ and $\sfp$, respectively. 
The estimator $\wh{\theta}_i$ may coincide with $\wh{\theta}_i^0$ but they can differ from each other in general. 
\bd
\im[{\bf[A5]}] 
\bd
\im[(i)] 
$\ds n^{1/2}\big(\wh{\theta}_1-\theta_1^*\big)-\bbT^{(1)}_n\to^p0$, 
and $\ds \wh{\theta}_1-\theta_1^* \yeq O_{{\vblue L^\inftym}}(n^{-1/2})$ as $n\to\infty$.

\im[(ii)] 
$\ds (nh)^{1/2}\big(\wh{\theta}_2-\theta_2^*\big)-\bbT^{(2)}_n\to^p0$, 
and $\ds \wh{\theta}_2-\theta_2^*\yeq O_{{\vblue L^\inftym}}(n^{-1/2}h^{-1/2})$  as $n\to\infty$.
\im[(iii)] 
$\ds \big(\bbT^{(1)}_n,\bbT^{(2)}_n,\Gamma_{33}^{-1}M_n^{(3)}\big)\to^d {\mathfrak L}$ as $n\to\infty$. 
\ed\ed

Let 
$
\wh{\vartheta}^\sfA
=
\big(\wh{\theta}_1,\wh{\theta}_2,\wh{\theta}_3^\sfA\big)
$. 
The matrix $a_n$ is given in (\ref{202306160306}). 
Then, as a corollary to Theorem \ref{201906031711}, 
\begin{theorem}\label{2206250928}
Suppose that $[A1]$ with {\vblue $\bbi=(1,1,2,1,3,3)$}, 
$[A2]$, 
$[A3]$ $(iii)$, $[A4]$ {\vredd $(i)$} and $[A5]$ are satisfied. Then 
\beas 
E\big[f\big(a_n^{-1}\big(\wh{\vartheta}^{\sfA}-\theta^*\big)\big)\big]
&\to&
\bbE[f({\mathfrak L})]
\eeas
as $n\to\infty$ for every $f\in C_p(\bbR^\sfp)$ and $\sfA\in\{M,B\}$. 
\end{theorem}

\begin{en-text}
In Theorem \ref{2206250928}, if $\wh{\theta}_3$ is constructed with $\big(\wh{\theta}_1,\wh{\theta}_2)$ 
for the initial estimator of $(\theta_1,\theta_2)$, then obviously we can erase Condition $[A4]$ from the assumptions 
because it is included in Condition $[A5]$. 
However, i
\end{en-text}
It is necessary to keep this general form of Theorem \ref{2206250928} distinguishing 
$\big(\wh{\theta}_1,\wh{\theta}_2)$ from $\big(\wh{\theta}_1^0,\wh{\theta}_2^0)$
when $\big(\wh{\theta}_1,\wh{\theta}_2)$ is constructed with {\xred$\wh{\theta}_3^M$} based on $\big(\wh{\theta}_1^0,\wh{\theta}_2^0)$ 
different from $\big(\wh{\theta}_1,\wh{\theta}_2)$. 
As seen in Section \ref{202007281718}, 
the quasi-likelihood estimators for $\theta_i$ ($i=1,2$) constructed only 
with the non-degenerate component of the stochastic differential equation (\ref{201905262213}) 
satisfy Condition $[A5]$ {\vredd (i) and (ii)}. So, if these estimators are chosen for $\wh{\theta}_i^0$, then 
we may obtain the joint convergence and $L^\inftym$-boundedness of the scaled error 
of $\big(\wh{\theta}_1^0,\wh{\theta}_2^0,\wh{\theta}_3^\sfA\big)$. 
However, the estimator of $\theta_1$ based only on the non-degenerate component is not efficient. 
That is, it is possible to improve it by incorporating with the information from the degenerate component. 
We will discuss this phenomenon in Section \ref{2206241735}.
}

{\vblue 
\section{Improving estimation of $\theta_1$}\label{2206241735}

Gloter and Yoshida \cite{gloter2020adaptive,gloter2021adaptive} show that 
the initial adaptive QMLE $\wh{\theta}_1^{0,M}$ given in Section \ref{202007281718} is asymptotically normal, 
and that 
its asymptotic variance is improved by a one-step estimator. 
It is possible to prove the $L^\inftym$-boundedness of the error of the one-step estimator, 
if one follows the proof of Gloter and Yoshida \cite{gloter2020adaptive,gloter2021adaptive}. 
However, for that, we need the $L^\inftym$-boundedness of the initial estimator, which 
is discussed in Section \ref{202007281718}. 
Here we will construct the second-stage estimators $\wh{\theta}_1^{\sfA}$ for $\theta_1$ 
in a different way.  In particular, the adaptive QBE can also be addressed, not only the adaptive QMLE. 

Consider a pair of initial estimators $\big(\wh{\theta}_2^0,\wh{\theta}_3^0)$ for $(\theta_2,\theta_3)$.
\begin{en-text}
satisfying {\vredd $[A4]$ (ii)} and 
\bd
\im[{\bf[A4]}] 
{\bf (iii)} 
$\ds \wh{\theta}_3^0-\theta_3^*\yeq O_{L^\inftym}(n^{-1/2}h^{1/2})$  as $n\to\infty$. 
\ed
\end{en-text}
Obviously the estimator $\wh{\theta}_3^\sfA$ in Section \ref{202007281646} can serve as $\wh{\theta}_3^0$, but 
we may consider any possibly less efficient estimator as $\wh{\theta}_3^0$; 
see Gloter and Yoshida \cite{gloter2021adaptive} 
for a simpler estimator that is less efficient but having a sufficient convergence rate.  
Define the random field $\bbH_n^{(1)}(\theta_1)$ by 
\bea\label{201906060545} 
\bbH_n^{(1)}(\theta_1)
 &=& 
 -\half\sum_{j=1}^n \bigg\{
S(Z_\tjm,\theta_1,\wh{\theta}_3^0)^{-1}\big[\cald_j(\theta_1,\wh{\theta}_2^0,\wh{\theta}_3^0)^{\otimes2}\big]
 +\log\det S(Z_\tjm,\theta_1,\wh{\theta}_3^0)\bigg\}.
\eea
Any measurable map $\wh{\theta}_1^M$ is called an adaptive quasi-maximum likelihood estimator for $\theta_1$ if it satisfies 
\beas 
\bbH_n^{(1)}\big(\wh{\theta}_1^M\big) 
&=& 
\sup_{\theta_1\in\ol{\Theta}_1}\bbH_n^{(1)}(\theta_1). 
\eeas
The adaptive quasi-Bayesian estimator $\wh{\theta}_1^B$ for $\theta_1$ is defined by 
\beas 
\wh{\theta}_1^B
&=&
\bigg(\int_{\Theta_1}\exp\big(\bbH_n^{(1)}(\theta_1)\big)\varpi_1(\theta_1)d\theta_1\bigg)^{-1}
\int_{\Theta_1}\theta_1\exp\big(\bbH_n^{(1)}(\theta_1)\big)\varpi_1(\theta_1)d\theta_1,
\eeas
where $\varpi_1$ is a continuous prior density for $\theta_1$ such that 
$0<\inf_{\theta_1\in\Theta_1}\varpi_1(\theta_1)\leq\sup_{\theta_1\in\Theta_1}\varpi_1(\theta_1)<\infty$. 

Let 
\bea\label{2206261226}
\Gamma_{11}
&=&
\half\int  \rm{Tr}\big\{S^{-1}(\partial_1S)S^{-1}\partial_1S(z,\theta_1^*, \theta_3^*)\big\}\nu(dz)
\nn\\&=&
\half\int \bigg[
 \rm{Tr}\big\{\big(C^{-1}(\partial_1C)C^{-1}\partial_1C\big)(z,\theta_1^*)\big\}
\nn \\&&\qquad\quad
 +\rm{Tr}\big\{\big(V^{-1}H_x(\partial_1C)H_x^\star V^{-1}H_x(\partial_1C)H_x^\star\big)(z,\theta_1^*,\theta_3^*)\big\}
 \bigg]\nu(dz).
\eea
{\colorg
If $H_x$ is an invertible (square) matrix, then $\Gamma_{11}$ coincides with 
\beas
\int  \rm{Tr}\big\{\big(C^{-1}(\partial_1C)C^{-1}\partial_1C\big)(z,\theta_1^*)\big\}\nu(dz).
\eeas
Otherwise, it is not always true. 
{\xred 
As a matter of fact, more generally, if 
$\sfd_Y \le \sfd_X$ and $\mathop{rank}(H_x)= \sfd_Y$, then we have
\bea
\label{E: comparison variance}	
		\int \text{Tr}\big\{\big(V^{-1}H_x(\partial_1C)H_x^\star V^{-1}H_x(\partial_1C)H_x^\star\big)(z,\theta_1^*,\theta_3^*) \big\}\nu(dz)
	\nn\\ \le \int \text{Tr}\big\{\big(C^{-1}(\partial_1C)C^{-1}\partial_1C\big)(z,\theta_1^*)\big\} \nu(dz).
\eea
Here the order $\leq$ is the partial order among $\sfp_1\times\sfp_1$ symmetric matrices; 
remark that the matrices on both sides are of $\sfp_1\times\sfp_1$ 
due to the product of two factors $\partial_1C$. 
See Section \ref{202402141102} for the proof of (\ref{E: comparison variance}). 
}
}

We need the following notations: 
\bea\label{201906030126}
\xi^{(\ref{201906030126})}_j
&=&
h^{-1/2}B(Z_\tjm,\theta_1^*)\Delta_jw,
\eea
\bea\label{201906030150}
\xi^{(\ref{201906030150})}_j
&=&
h^{-1/2}(B_xB)(Z_\tjm,\theta_1^*)\int_\tjm^\tj \int_\tjm^tdw_sdw_t,
\eea
\bea\label{201906030334}
\xi_j^{(\ref{201906030334})}
&=&
\kappa(Z_\tjm,\theta_1^*,\theta_3^*)\zeta_j,
\eea
and
\bea\label{201906030345}
\xi_j^{(\ref{201906030345})}
&=&
((H_xB)_xB)(Z_\tjm,\theta_1^*,\theta_3^*)\int_\tjm^\tj\int_\tjm^t\int_\tjm^sdw_rdw_sdt.
\eea

Let 
\beas 
M_n^{(1)} &=& 
\half n^{-1/2}\sum_{j=1}^n
\big(S^{-1}(\partial_1S)S^{-1}\big)(Z_\tjm,\theta_1^*,\theta_3^*)
\big[\widetilde{\cald}_j(\theta_1^*,\theta_2^*,\theta_3^*)^{\otimes2}-S(Z_\tjm,\theta_1^*,\theta_3^*)\big], 
\eeas
where
\bea\label{2206290639}
\widetilde{\cald}_j(\theta_1^*,\theta_2^*,\theta_3^*)
&=&
\left(\begin{array}{c}
{\cred\xi^{(\ref{201906030126})}_j+\xi^{(\ref{201906030150})}_j}\y
h^{-3/2}\big(\xi_j^{(\ref{201906030334})}+\xi_j^{(\ref{201906030345})}\big)
\end{array}\right).
\eea
}

{\vblue  
For the adaptive QMLE and the adaptive QBE for $\theta_1$ 
with respect to $\bbH_n^{(1)}$, we have the following results. 
\begin{theorem}\label{2206261240}
Suppose that $[A1]$ with $\bbi=(1,1,2,3,3,1)$, $[A2]$, $[A3]$ $(i')$ and {\vredd $[A4]$ $(ii), (iii)$}
are satisfied. 
Then
\bd
\im[(a)] For $\sfA\in\{M,B\}$,
\beas 
n^{1/2}\big(\wh{\theta}_1^\sfA-\theta_1^*\big)
- \Gamma_{11}^{-1}M_n^{(1)}
&\to^p& 0
\eeas
as $n\to\infty$. 
\im[(b)] For $\sfA\in\{M,B\}$, 
\beas 
E\big[f\big(n^{1/2}\big(\wh{\theta}_1^\sfA-\theta_1^*\big)\big)\big]
&\to&
\bbE[f({\sf Z}_1)]
\eeas
as $n\to\infty$ for every $f\in C_p(\bbR^{\sfp_1})$, where  
${\sf Z}_1\sim N(0,\Gamma_{11}^{-1})$. 
\ed
\end{theorem}
See Section \ref{2206261259} for the proof. 
}

\section{{\vblue $L^\inftym$-boundedness of the initial estimators for $\theta_1$ and $\theta_2$}}\label{202007281718}
{\vblue 
We will construct initial estimators for $\theta_1$ and $\theta_2$ in this paper's setting, 
and investigate their asymptotic properties by means of the QLA. 
In particular, our approach enables us to treat quasi-Bayesian estimators, which is important 
not only to establish Bayesian inference but also to construct 
adaptive methods and hybrid methods (Uchida and Yoshida \cite{uchida2012adaptive}, 
Kamatani and Uchida \cite{KamataniUchida2014}). 
The validation of the Bayesian method is a problem that has the same root as 
the $L^\inftym$-boundedness of the QMLE. 
The QLA solves this difficulty by a polynomial type large deviation inequality 
generally proved under the asymptotic quadradicity of the quasi-likelihood ratio (Yoshida \cite{yoshida2011polynomial,yoshida2021simplified}). 

For $\theta_1$, we consider the random field
\beas 
\bH^{(1)}_n(\theta_1) 
&=& 
-\half \sum_{j=1}^n \bigg\{C(Z_\tjm,\theta_1)^{-1}\big[h^{-1}(\Delta_jX)^{\otimes2}\big]
+\log\det C(Z_\tjm,\theta_1)\bigg\},
\eeas
where $\Delta_jX=X_\tj-X_\tjm$. 
The random field $\bH^{(1)}_n$ is a continuous function on $\overline{\Theta}_1$ a.s. under Conditions $[A1]$ and $[A2]$ (iii). 
Differently from the random field $\bbH^{(1)}_n(\theta_1)$ of (\ref{201906060545}), 
$\bH^{(1)}_n(\theta_1)$ uses the increments of $X$ only. 

Given the data $(Z_\tj)_{j=0,1,...,n}$, 
the quasi-maximum likelihood estimator (QMLE) for $\theta_1$
is any measurable function $\wh{\theta}_1^{0,M}=\wh{\theta}_{1,n}^{0,M}$ of $(Z_\tj)_{j=0,1,...,n}$ satisfying 
\beas 
\bH^{(1)}_n(\wh{\theta}_1^{0,M}) &=& \max_{\theta_1\in\overline{\Theta}_1}\bH^{(1)}_n(\theta_1)\quad a.s.
\eeas
The quasi-Bayesian estimator (QBE) for $\theta_1$ is defined by 
\bea\label{202402141203}
\wh{\theta}_1^{0,B} 
&=& 
\bigg(\int_{\Theta_1}\exp\big(\bH^{(1)}(\theta_1)\big)\varpi^{(1)}(\theta_1)d\theta_1\bigg)^{-1}
\int_{\Theta_1}\theta_1\exp\big(\bH^{(1)}(\theta_1)\big)\varpi^{(1)}(\theta_1)d\theta_1,
\eea
where $\varpi^{(1)}$ is a continuous prior density for $\theta_1$ satisfying 
$0<\inf_{\theta_1\in\Theta_1}\varpi^{(1)}(\theta_1)\leq\sup_{\theta_1\in\Theta_1}\varpi^{(1)}(\theta_1)<\infty$. 
The QBE $\wh{\theta}_1^{0,B}$ takes values in the convex hull of $\Theta_1$, and needs modification so as to be $\ol{\Theta}_1$-valued, 
when used as an initial estimator of $\theta_1$. 
However, the asymptotic properties of the modified estimator are unchanged since $\theta_1^*$ is in the open set $\Theta_1$. 
See the discussion in the paragraph including (\ref{202306172011}). 

The information matrix $\Gamma^{(1)}$ for $\theta_1$ evaluated at $\theta_1^*$ 
is a $\sfp_1\times\sfp_1$ symmetric matrix defined by 
\bea\label{201905281059}
\Gamma^{(1)}[u_1^{\otimes2}]
&=&
\half\int_{\bbR^{\sfd_Z}}\text{Tr}\big\{C^{-1}(\partial_1C)[u_1]C^{-1}(\partial_1C)[u_1](z,\theta_1^*)\big\}\nu(dz)
\eea
for $u_1\in\bbR^{\sfp_1}$. 
The quasi-likelihood analysis establishes that 
the QMLE $\wh{\theta}_1^{0,M}$ and the QBE $\wh{\theta}_1^{0,B}$ share the same asymptotic normality and moments convergence, 
besides they are asymptotically equivalent. 
\begin{en-text}
Routinely, $n^{1/2}$-consistency and asymptotic normality of $\wh{\theta}_1^0$ can be established. 
We will give a brief for self-containedness and for the later use. 

We will see the existence and positivity of $\Gamma^{(1)}$ in the following theorem. 
{\tred We refer the reader to Gloter and Yoshida \cite{gloter2021adaptive} 
for a proof. }
\end{en-text}
}
{\vblue 
Let 
\beas 
\wh{M}^{(1)}_n 
&=&
\half n^{-1/2}\sum_{j=1}^n \big(C^{-1}(\partial_1C)C^{-1}\big)(Z_\tjm,\theta_1^*)
\big[\big(h^{-1/2}B(Z_\tjm,\theta_1^*)\Delta_jw\big)^{\otimes2}
-C(Z_\tjm,\theta_1^*)\big].
\eeas
\begin{theorem}\label{201905291601}
Suppose that $[A1]$ with $\bbi=(1,0,2,3,0,0)$, 
$[A2]$ $(i)$, $(ii)$, $(iii)$, and 
$[A3]$ $(i)$ are satisfied. 
Then 
\bd\im[(a)] Then $\Gamma^{(1)}$ exists and is positive-definite, and 
\beas
\sqrt{n}\big(\wh{\theta}_1^{0,\sfA}-\theta_1^*\big)
-(\Gamma^{(1)})^{-1}\wh{M}^{(1)}_n&\to^p&0
\eeas
as $n\to\infty$ for $\sfA\in\{M,B\}$. 
Moreover, 
$
{\vredd \wh{M}}^{(1)}_n
\to^d 
N_{\sfp_1}(0,\Gamma^{(1)})
$
as $n\to\infty$. 
\im[(b)] For $\sfA\in\{M,B\}$, 
\beas 
E\big[f\big(\sqrt{n}\big(\wh{\theta}_1^{0,\sfA}-\theta_1^*\big)\big)\big]
&\to&
\bbE[f(\eta_1)]
\eeas
as $n\to\infty$ for evey $f\in C_p(\bbR^{\sfp_1})$, where 
$\eta_1\sim N_{\sfp_1}\big(0,(\Gamma^{(1)})^{-1}\big)$. 
\ed
\end{theorem}
}
\begin{en-text}
\begin{remark}\rm
It is possible to show that the quasi-Bayesian estimator (QBE) also enjoys 
the same asymptotic properties as the QMLE in Theorem \ref{201905291601}, 
if we follows the argument in Yoshida \cite{yoshida2011polynomial}. 
This means we can use both estimators together with the estimator for $\theta_2$ 
e.g. given in Section {\tred\ref{202007281718}}, 
to construct a one-step estimator for $\theta_3$ 
based on the scheme presented in Section {\tred\ref{202007290607}}, 
and consequently we can construct a one-step estimator for 
$\theta=(\theta_1,\theta_2,\theta_3)$ by the method in Section \ref{202001141623}. 
\end{remark}
\end{en-text}

{\vblue 

We shall turn our attention to estimation of $\theta_2$. 
An adaptive scheme is considered. 
Suppose that an estimator $\wh{\theta}_1^0$ based on the data $(Z_\tj)_{j=0,1,...,n}$ 
satisfies Condition {\fred$[A4]$} (i). 
\begin{en-text}
, i.e., 
\beas 
\wh{\theta}_1^0 -\theta_1^* &=& O_p(n^{-1/2})
\eeas
as $n\to\infty$. 
\end{en-text}
It is clear that the above constructed estimators $\wh{\theta}_1^{0,\sfA}$ for $\theta_1$ serve for $\wh{\theta}_1^0$. 
However, one can adopt any estimator satisfying this condition.

The random field $\bbH^{(2)}_n$ on $\overline{\Theta}_2$ to be considered is 
\bea\label{202002171935} 
\bH^{(2)}_n(\theta_2) 
&=&
-\half\sum_{j=1}^n C(Z_\tjm,\wh{\theta}_1^0)^{{\cred -1}}
\big[h^{-1}\big(\Delta_jX-hA(Z_\tjm,\theta_2)\big)^{\otimes2}\big].
\eea
Denote by $\wh{\theta}_2^{0,M}=\wh{\theta}_{2,n}^{0,M}$ any sequence of 
quasi-maximum likelihood estimator for $\bH^{(2)}_n$, that is, 
\beas 
\bH^{(2)}_n(\wh{\theta}_2^{0,M}) &=& \sup_{\theta_2\in\overline{\Theta}_2}\bH^{(2)}_n(\theta_2).
\eeas
The adaptive quasi-Bayesian estimator for $\theta_2$ is defined by 
\beas
\wh{\theta}_2^{0,B}
&=&
\bigg(\int_{\Theta_2}\exp\big(\bH^{(2)}_n(\theta_2)\big)\varpi^{(2)}(\theta_2)d\theta_2\bigg)^{-1}
\int_{\Theta_2}\theta_2\exp\big(\bH^{(2)}_n(\theta_2)\big)\varpi^{(2)}(\theta_2)d\theta_2, 
\eeas
where 
$\varpi^{(2)}$ is a continuous prior density function on $\Theta_2$ such that 
$0<\inf_{\theta_2\in\Theta_2}\varpi^{(2)}(\theta_2)\leq\sup_{\theta_2\in\Theta_2}\varpi^{(2)}(\theta_2)<\infty$. 

The information matrix for $\theta_2$ evaluated at $\theta_2^*$ is given by 
\bea\label{202001281931} 
\Gamma_{22}
&=&
\int S(z,\theta_1^*,\theta_3^*)^{-1}
\left[
\left(\begin{array}{c}\partial_2A(z,\theta_2^*)\\
2^{-1}\partial_2L_H(z,\theta_1^*,\theta_2^*,\theta_3^*)
\end{array}\right)^{\otimes2}
\right]\nu(dz)
\nn\\&=&
\int \partial_2A(z,\theta_2^*)^\star C(z,\theta_1^*)^{-1}\partial_2A(z,\theta_2^*)\nu(dz)
\eea
and the quasi-score function evaluated at $\theta_2^*$ is 
\bea\label{202002171941} 
\wh{M}^{(2)}_n
&=&
T^{-1/2}\sum_{j=1}^n C(Z_\tjm,\theta_1^*)^{-1}
\big[B(Z_\tjm,\theta_1^*)\Delta_jw,\partial_2A(Z_\tjm,\theta_2^*)\big], \quad T=nh.
\eea

See Section \ref{2206281655} for the proof of the following theorem. 
}

{\vblue 
\begin{theorem}\label{202002172010}
Suppose that Conditions $[A1]$ with $\bbi=(1,3,1,1,0,0)$, $[A2]$ $(i)$-$(iii)$, 
$[A3]$ $(ii)$ and $[A4]$ $(i)$. 
Then 
\bd
\im[(a)] For $\sfA\in\{M,B\}$, 
\beas
(nh)^{1/2}\big(\wh{\theta}_2^{0,\sfA}-\theta_2^*\big)
-\Gamma_{22}^{-1}\wh{M}^{(2)}_n
&\to^p&
0
\eeas
and 
$
\wh{M}^{(2)}_n
\to^d
N(0,\Gamma_{22})
$
as $n\to\infty$. 
\im[(b)] For $\sfA\in\{M,B\}$, 
\beas
E\big[f\big((nh)^{1/2}\big(\wh{\theta}_2^{0,\sfA}-\theta_2^*\big)\big)\big]
&\to& 
\bbE[f({\sf Z}_2)]
\eeas
as $n\to\infty$ 
for any $f\in C_p(\bbR^{\sfp_2})$, where ${\sf Z}_2\sim N(0,\Gamma_{22}^{-1})$. 
\ed
\end{theorem}
\halflineskip

{\vredd 
The present $\bbi=(1,3,1,1,0,0)$ is better than 
$\bbi=(1,3,2,1,0,0)$ in the previous work \cite{gloter2020adaptive,gloter2021adaptive}. }

\section{Asymptotic properties of the adaptive estimators}\label{202306160251}
According to the discussion so far, a possible adaptive scheme is to compute the estimators 
through the procedure given in Section \ref{202306160554}. 
The following theorem validates this method. 
\begin{theorem}\label{2206290542}
Suppose that Conditions $[A1]$ with $\bbi=(1,3,2,3,3,3)$, $[A2]$, $[A3]$ $(i')$, $(ii)$, $(iii)$. 
Suppose that 
$\wh{\theta}_1^{\sfA_1}$ $(\sfA_1\in\{M,B\})$, 
$\wh{\theta}_2^{\sfA_2}=\wh{\theta}_2^{0,\sfA_2}$  $(\sfA_2\in\{M,B\})$ and 
$\wh{\theta}_3^{\sfA_3}$  $(\sfA_3\in\{M,B\})$ 
are any estimators given in Section \ref{202306160554}. 
Then the joint estimator $\wh{\vartheta}=\big(\wh{\theta}_1^{A_1},\wh{\theta}_2^{A_2},\wh{\theta}_3^{A_3}\big)$ admits the convergence 
\beas 
E\big[f\big(a_n^{-1}(\wh{\vartheta}-\theta^*)\big)\big]
&\to&
\bbE[f\big({\sf Z}_1,{\sf Z}_2,{\sf Z}_3\big)]
\eeas
as $n\to\infty$, for every $f\in C_p(\bbR^\sfp)$, where 
$\big({\sf Z}_1,{\sf Z}_2,{\sf Z}_3\big)\sim N_{\sfp_1}(0,\Gamma_{11}^{-1})\otimes N_{\sfp_2}(0,\Gamma_{22}^{-1})\otimes N_{\sfp_3}(0,\Gamma_{33}^{-1})$. 
\end{theorem}
\halflineskip
%

The proof of Theorem \ref{2206290542} is in Section \ref{202306160609}. 
{\vredd An asymptotically equivalent representation of $a_n^{-1}(\wh{\vartheta}-\theta^*)$ is also given in (\ref{202306230102}). }

\begin{en-text}
{\vgreen
\begin{enumerate}[(1)]
\im $\wh{\theta}_1^{0,A_0}$ for $A_0=M/B$. 
\im $\wh{\theta}_2^{0,A_2}$ for $A_2=M/B$ by using $\wh{\theta}_1^{0,A_0}$. 
\im $\wh{\theta}_3^{A_3}$ for $A_3=M/B$ by using $\big(\wh{\theta}_1^{0,A_0},\wh{\theta}_2^{0,A_2}\big)$. 
\im $\wh{\theta}_1^{A_1}$ for $A_1=M/B$ by using $\big(\wh{\theta}_2^{0,A_2},\wh{\theta}_3^{A_3}\big)$. 
\im $\wh{\theta}_2^{A_2}=\wh{\theta}_2^{0,A_2}$ for $A_2=M/B$. 
\end{enumerate}
Then $\wh{\vartheta}=\big(\wh{\theta}_1^{A_1},\wh{\theta}_2^{A_2},\wh{\theta}_3^{A_3}\big)$ is jointly asymptotically normal. 
}
\end{en-text}
}

\begin{en-text}
{\tred 
\begin{remark}\label{202008010413}\rm 
Another choice of the initial estimator $\wh{\theta}_2^0$ is 
a simple least squares estimator using the coefficient $A$. 
It is less efficient than the quasi-maximum likelihood estimator for  
the first component of the model. 
Theorem \ref{202007281651} in Section \ref{202001141623} shows the one-step estimator for $\theta_2$ 
recovers efficiency even if such a less efficient estimator is used as 
the initial estimator for $\theta_2$. 
\end{remark}
}
\end{en-text}

{\vblue  
}

{\xred 
\section{Numerical simulations}\label{202402141230}
In this section, we study the behaviour of the quasi-Bayesian estimator (QBE) on finite sample, and compare the quality of estimation with the one of the QMLE as studied in 
{\xblue Gloter and Yoshida \cite{gloter2024nonadaptive}}.
We start with the linear model, $Z=(X,Y)\in \mathbb{R}^2$ solution of 
\begin{equation*}
	\left\{\begin{array}{ccl}
		dX_t &=& (-\theta_{2,1} X_t - \theta_{2,2}Y_t ) dt + \theta_1 dw_t
		\y	dY_t&=&\theta_3 X_t dt 
	\end{array}\right. 
\end{equation*}
where $\theta_1>0$, $\theta_{2,1}>0$, $\theta_{2,2}>0$, $\theta_3 >0$. We assume that the true value of the parameters are $\theta_1=1$, $\theta_2=(1,1)$, $\theta_3=1$.

We first compute {\xblue$\wh{\theta}_1^{0,B}$} 
as defined in Section \ref{202007281718} with the choice of uniform prior on a $[10^{-4},10]$. 
To compute the value of \eqref{202402141203}, we use a numerical integration method provided in Python to approximate the integrals. 
We compare this method {\xblue with} a Metropolis-Hasting (M-H) method where we sample the law of the posterior and derive an approximation of the Bayesian estimator as the empirical mean of the M-H Markov chain. From Table \ref{T: Bayesian linear}, we see that both methods provide a good estimation of the parameter $\theta_1$. The first column of results shows the results when using a numerical method to compute the integrals involved in the definition of
{\xblue$\wh{\theta}_1^{0,B}$} 
 while the second shows the results when the estimator is computed by sampling the mean of a M-H method with length $5000$. We see that the two methods give similar results. In practice, the computation of $\hat{\theta}^B_2$ reveals being infeasible through the numerical integration method. This is certainly due to the very steep shape of the quasi likelihood function near the true value of the parameter and the fact that the integration is here multidimensional. Hence, for $\theta_2^B$, and  the subsequent estimators $\hat{\theta}_3^B$ and $\hat{\theta}_1^B$, we only use the M-H method. We see that all the parameters are well estimated. 
 As expected, the estimator $\hat{\theta}_1^B$ has a reduced variance when compared to the initial estimator $\hat{\theta}_1^{0,B}$.  Comparing the results of Table \ref{T: Bayesian linear} for the Bayesian estimator with the study  based on maximization of the quasi-likelihood function given in \cite{gloter2024nonadaptive}, it appears that the results are very similar for both estimators.
\begin{table}
	\caption{ Estimator $\hat{\theta}^B$, $nh=100$.
	\label{T: Bayesian linear}
}\centering
\begin{tabular}{|c||c|c|c|c|c|c|c|}
	\hline
	& &$\hat{\theta}_1^{0,B,\text{num}}$ & $\hat{\theta}_1^{0,B}$ &$\hat{\theta}_{2,1}^B$ & $\hat{\theta}_{2,2}^B$ &$\hat{\theta}_{3}^B$&$\hat{\theta}_1^{B}$
	\\ \hline
	&True value 	&	1 & 1& 1 &1 &1 & 1
	\\
	\hline
	\hline
	\hline
	$h=1/10$&	Mean 
	& 1.002& 1.002  & 1.003  &  0.983 & 1.000 & 0.981
	\\ \hline
	&(std)
	& (2.5e-2)&(2.5e-2) &(0.18) & (0.18) &  (4.9e-3)& (1.8e-2)
	\\
	\hline
	\hline
$	h=1/5$&	Mean 
	& 1.003& 1.003  & 1.002  &  0.927 & 0.999 & 0.963
	\\ \hline
	&(std)
	& (3.0e-2)&(3.0e-2) &(0.14) & (0.14) &(1.1e-2)  &(2.4.e-2)
	\\
	\hline
\end{tabular}
\end{table}

We also studied the non-linear FitzHugh-Nagumo model,
\beas	
\left\{\begin{array}{ccl}
	dX_t&=&(\gamma Y_t -X_t + \beta) dt + \sigma dw_t
	\y dY_t&=&\frac{1}{\varepsilon}(Y_t-Y_t^3-U_t+s) dt
\end{array}\right.
\eeas
where $\varepsilon>0$, $\beta>0$, $\gamma>0$ and $s\in\mathbb{R}$.
It turns out that the preliminary estimator of the diffusion parameter $\hat{\sigma}^{0,B}$ works well, both by numerical integration or Metropolis-Hasting methods. On the other hand, when estimating the drift parameters $\theta_2=(\gamma,\beta)$ and the parameter $\theta_3=(\varepsilon,s)$ we had to use a very large number of iterations in the M-H method in order to get a precise estimation of these parameters. Results are given in Table \ref{T: Bayesian FHN} for the preliminary estimator of $\sigma$, estimation of $\theta_2=(\gamma,\beta)$ and $\theta_3=(\varepsilon,s)$, and the improved estimator of $\sigma$. The Metropolis Hasting method uses series of length $10^6$. We see that, except for the $\hat{\sigma}^B$ the estimation results appear correct and in line with the results of \cite{gloter2024nonadaptive} where the QLE estimator is used.
Let us stress that the estimation of $\theta_3$, whose rate is the fastest, is very sensitive to the length of the series used in the M-H method. Indeed, for $h=1/30$, we also made simulations using series of length $10^5$ in the M-H method and the variance of the estimator $\theta_3$ was increased by a factor larger than five compared to the one given in Table \ref{T: Bayesian FHN}. It is certainly due to the fact that the variance of the M-H method dominates the intrinsic variance of the estimator when the Markov chain used is too short.

The implementation of the improved estimator $\hat{\sigma}^B$ is more delicate. Indeed, it seems necessary that the sampling step is small enough to reduce the bias and variance below the ones of the initial estimator
$\hat{\sigma}^{0,B}$. This was already observed in \cite{gloter2024nonadaptive} for the QLE. Meanwhile, the estimator $\hat{\sigma}^B$ is sensitive to the quality of estimation of the parameter $\theta_3$, and hence on the precision of the M-H method used to evaluate $\hat{\theta}_3^B$. In our simulations, the improved estimator only outperforms the initial one for {\xblue$h=1/30$} 
and when using series of length $10^6$ in the M-H method. 
However, we should remark that the balance condition $nh^2\to0$ is not fully satisfied because $nh^2=1/3$ even in the case $h=1/30$. 

\begin{table}[H]
	\caption{ Estimator $\hat{\theta}^B=(\hat{\sigma}^{0,B},\hat{\gamma}^B,\hat{\beta}^B,\hat{\varepsilon}^B,\hat{s}^B,\hat{\sigma}^B)$, $nh=10$, FHN model.
		\label{T: Bayesian FHN}
	}\centering
	\begin{tabular}{|c||c|c|c|c|c|c|c|c|}
		\hline
		& &$\hat{\sigma}^{0,B,\text{num}}$&$\hat{\sigma}^{0,B}$ &$\hat{\gamma}^B$ & $\hat{\beta}^B$ &$\hat{\varepsilon}^B$ & $\hat{s}^B$ &$\hat{\sigma}^{B}$
		\\ \hline
		&True value 	&	0.3 & 0.3 & 1.5& 0.8 & 0.1 & 0  &0.3
		\\
		\hline
		\hline
		$h=1/30$&	Mean 
		& 0.326 & 0.326 & 1.56 &  0.85& 0.102 & 4.8e-5 & 0.303
		\\
		\hline
		&	(std)
		& (1.7e-2) & (1.7e-2) & (0.29) & (0.25) & (8.6e-4) & (1.6e-3) & (1.4e-2)
		\\
		\hline
		\hline
		$h=1/20$&	Mean 
		& 0.338 & 0.338  & 1.54 &  0.83 & 0.103 & 3.5e-4 & 0.333
		\\
		\hline
		&	(std)
		& (2.3e-2) & (2.3e-2) & (0.27) & (0.23) & (1.1e-3) & (2.1e-3) & (3.0e-2)
		\\
		\hline
		\hline
		$h=1/10$&	Mean 
		& 0.368 & 0.368 & 1.44 &  0.762 & 0.112  & 2.4e-3 & 0.479
		\\
		\hline
		&	(std)
		& (4.1e-2) & (4.1e-2) & (0.22) & (0.17) & (5.1e-3) & (6.7e-3) & (1.1e-1)
		\\
		\hline
		\hline
		$h=1/5$&	Mean 
		& 0.417 & 0.417  & 1.30 &  0.67& 0.134& 1.2e-2 & 0.744
		\\
		\hline
		&	(std)
		& (6.8e-2) & (6.8e-2) & (0.19) & (0.15) & (1.6e-2) & (2.7e-2) & (2.2e-1)
		\\
		\hline
		\hline
	\end{tabular}
\end{table}
}

%
%
%
%
%
%

{\vred 
\section{Proof of Theorem \ref{201906031711}}\label{202101221241}
\subsection{{\vblue Convergence of the random field $\bbY^{(3)}_n(\theta_3)$}
}
{\vblue We start with the following result. }

\begin{lemma}\label{201906021729}
Suppose that $[A1]$ with $(i_A,j_A,i_B,j_B,i_H,j_H)=(0,0,0,1,1,1)$ and $[A2]$ $(i)$, $(iii)$ and $(iv)$ are satisfied. 
Then 
\beas&&
\sup_{t\in\bbR_+}\bigg\|\sup_{(\theta_1,\theta_3)\in\overline{\Theta}_1\times\overline{\Theta}_3}
\big\{\big|S(Z_t,\theta_1,\theta_3)\big|+
\det S(Z_t,\theta_1,\theta_3)^{-1}
+\big|S(Z_t,\theta_1,\theta_3)^{-1}\big|
\big\}\bigg\|_p
\><\>
\infty
\eeas
for every $p>1$ 
\end{lemma}
{\vblue\noindent 
We refer the reader to Gloter and Yoshida \cite{gloter2021adaptive} for a proof of Lemma \ref{201906021729}. 
In what follows, the estimates in 
Lemma \ref{201906021729} will often be used without explicitly mentioned. 

According to the theory of the quasi-likelihood analysis, we will work with the random field
\beas 
\bbY^{(3)}_n(\theta_3)
&=&
n^{-1}h\big\{\bbH^{(3)}_n(\theta_3)-\bbH^{(3)}_n(\theta_3^*)\big\}.
\eeas
Recall $\bbi=(i_A,j_A,i_B,j_B,i_H,j_H)$. 
\begin{lemma}\label{201906021825}
Suppose that $[A1]$ with $\bbi=(1,{\colorg1},{\colorg2},{\colorr1},3,1)
$, $[A2]$ 
and $[A4]$ $(i)$ 
are satisfied. 
Then, there exists a positive number $\ep_1$ such that 
\bea\label{201906021845}
\sup_{\theta_3\in\overline{\Theta}_3}
\big|\bbY_n^{(3)}(\theta_3)
-\bbY^{(3)}(\theta_3)
\big|
&=&
O_{L^\inftym}(n^{-\ep_1})
\eea
as $n\to\infty$. 
\end{lemma}
}

\begin{en-text}
{\tgreen
Let 
\beas 
\bbY^{(3)}_n(\theta_3)
&=&
n^{-1}h\big\{\bbH^{(3)}_n(\theta_3)-\bbH^{(3)}_n(\theta_3^*)\big\}.
\eeas
\begin{theorem}
Suppose that $[A1]$ with $(i_A,j_A,i_B,j_B,i_H,j_H)=(1,{\colorg1},{\colorg2},{\colorr1},3,1)
$ and $[A2]$ are satisfied. 
Then   
\bea\label{201906021845}
\sup_{\theta_3\in\overline{\Theta}_3}
\big|\bbY_n^{(3)}(\theta_3)
-\bbY^{(3)}(\theta_3)
\big|
&\to^p&
0
\eea
as $n\to\infty$, if $\wh{\theta}_1^0\to^p\theta_1^*$ and $\wh{\theta}_2^0\to^p\theta_2^*$. 
Moreover, $\wh{\theta}_3^0\to^p\theta_3^*$ 
if $[A3]$ $(iii)$ is additionally satisfied. 
\end{theorem}
}
\end{en-text}
\begin{en-text}
{\cred {\bf Modify this paragraph later.} 
Remark that if we assume $[A3]$ $(i)$ in the first part of Theorem \ref{201906021825} 
(in place of the convergence of $\wh{\theta}_1^0$), then 
the convergence $\wh{\theta}_1^0\to^p\theta_1^*$ 
for the QMLE $\wh{\theta}_1$ 
follows from 
Theorem \ref{201905291601} (a). \koko about $\wh{\theta}_2^0$\y
}
\end{en-text}
{\vblue 
\noindent
\proof 
Let 
\beas 
\delta_j(\theta_1,\theta_2,\theta_3) 
&=& 
-\cald_j(\theta_1,\theta_2,\theta_3)+\cald_j(\theta_1,\theta_2,\theta_3^*).
\eeas
Then 
Lemma \lsix of \cite{gloter2021adaptive} gives 
\bea\label{0406241050}
\sup_{(\theta_1,\theta_2,\theta_3)\in\overline{\Theta}_1\times\overline{\Theta}_2\times\ol{\Theta}_3}
\big|\delta_j(\theta_1,\theta_2,\theta_3)\big|
&=&
O_{L^\inftym}(h^{-1/2})
\eea
as $n\to\infty$. 

We have
\beas 
\bbY_n^{(3)}(\theta_3)
&{\tred=}& 
-\frac{1}{2n}\sum_{j=1}^n  \wh{S}(Z_\tjm,\theta_3)^{-1}
\big[\big(h^{1/2}\delta_j(\wh{\theta}_1^0,\wh{\theta}_2^0,\theta_3)\big)^{\otimes2}\big]
+n^{-1}hR^{(\ref{201906021417})}_n(\theta_3)
\eeas
for
\bea\label{201906021417}
R^{(\ref{201906021417})}_n(\theta_3)
&=&
{\tred 
{\tgreen h^{-1/2}}
\sum_{j=1}^n 
\wh{S}(Z_\tjm,\theta_3)^{-1}\big[h^{1/2}\delta_j(\wh{\theta}_1^0,\wh{\theta}_2^0,\theta_3),
\cald_j(\wh{\theta}_1^0,\wh{\theta}_2^0,\theta_3^*)\big]
}
\nn\\&&
-\half\sum_{j=1}^n\big(
\wh{S}(Z_\tjm,\theta_3)^{-1}-\wh{S}(Z_\tjm,\theta_3^*)^{-1}\big)
\big[\cald_j(\wh{\theta}_1^0,\wh{\theta}_2^0,\theta_3^*)^{\otimes2}\big]
\nn\\&&
-\half\sum_{j=1}^n \log\frac{\det \wh{S}(Z_\tjm,\theta_3)}{\det \wh{S}(Z_\tjm,\theta_3^*)}.
\eea
Then, from 
(\ref{0406241100}) of Lemma \ref{2206241411}, (\ref{0406241050}), (\ref{201906021417}) 
and Lemma \ref{201906021729}, we obtain  
\bea\label{201906021619} 
n^{-1}h\sup_{\theta_3\in\ol{\Theta}_3}\big|R^{(\ref{201906021417})}_n(\theta_3)
\big|
&=&
{\tgreen O_{{\vblue L^\inftym}}(h^{1/2})+}
{\colorr O_{{\vblue L^\inftym}}(h)}{\yeq\tred O_{{\vblue L^\inftym}}(h^{1/2})}. 
\eea

Let 
\bea\label{201906021814}
\bbY^{(\ref{201906021814})}_n(\theta_1,\theta_3)
&=&
-\frac{1}{2n}\sum_{j=1}^nS(Z_\tjm,\theta_1,\theta_3)^{-1}
\left[\left(\begin{array}{c}0\\
H(Z_\tjm\,\theta_3)-H(Z_\tjm\,\theta_3^*)
\end{array}\right)^{\otimes2}\right]
\nn\\&=&
-\frac{6}{n}\sum_{j=1}^nV(Z_\tjm,\theta_1,\theta_3)^{-1}
\left[\left(
H(Z_\tjm\,\theta_3)-H(Z_\tjm\,\theta_3^*)
\right)^{\otimes2}\right].
\eea

Since the functions $A(z,\theta_2)$, $H(z,\theta_3)$ and $L_H(z,\theta_1,\theta_2,\theta_3)$ are 
dominated by a polynomial in $z$ uniformly in $\theta$, 
by using the  formula
\beas 
h^{1/2}\delta_j(\wh{\theta}_1^0,\wh{\theta}_2^0,\theta_3) 
&=& 
\left(\begin{array}{c}
0
\y
\left\{\begin{array}{c}
H(Z_\tjm,\theta_3)-H(Z_\tjm,\theta_3^*)
\\
+\frac{h}{2}\big(L_H(Z_\tjm,\wh{\theta}_1^{{\cred0}},\wh{\theta}_2^0,\theta_3)-L_H(Z_\tjm,\wh{\theta}_1^{{\cred0}},\wh{\theta}_2^{{\cred0}},\theta_3^*)\big)
\end{array}\right\}
\end{array}\right), 
\eeas
we see 
\bea\label{201906021842} 
\sup_{\theta_3\in\overline{\Theta}_3}
\left|\bbY^{(3)}_n(\theta_3)
-\bbY^{(\ref{201906021814})}_n(\wh{\theta}_1^0,\theta_3)
\right|
&=&
{\vblue O_{L^\inftym}(h^{1/2})},
\eea
taking (\ref{201906021619}) into account. 
%
Moreover, since the derivative $\partial_1V(z,\theta_1,\theta_3)$ is dominated by a polynomial in $z$ 
uniformly in $\theta$, 
\bea\label{201906021843}
\sup_{\theta_3\in\overline{\Theta}_3}\big|
\bbY^{(\ref{201906021814})}_n(\wh{\theta}_1^0,\theta_3)-\bbY^{(\ref{201906021814})}_n(\theta_1^*,\theta_3)\big|
\yeq O_{{\vblue L^\inftym}}(n^{-1/2})
\eea
with the help of the assumption that {\vredd $\wh{\theta}_1^0-\theta_1^*=O_{L^\inftym}(n^{-1/2})$, though this rate is too much 
for the present use}.

Thanks to the stationarity of $Z$, 
$E\big[\bbY^{(\ref{201906021814})}_n(\theta_1^*,\theta_3)\big]=\bbY^{(3)}(\theta_3)$. 
By using Sobolev's inequality and  Lemma 4 of Yoshida \cite{yoshida2011polynomial} with the mixing property, we find
\bea\label{22060241235}
\bigg\|\sup_{\theta_3\in\ol{\Theta}_3}\big|
\bbY^{(\ref{201906021814})}_n(\theta_1^*,\theta_3)-\bbY^{(3)}(\theta_3)\big|\bigg\|_p
&=&
O((nh)^{-1/2})+O((nh)^{-\frac{p-1}{p}})
\nn\\&=&
O((nh)^{-1/2})
\leq 
O(n^{-\ep_0/2})
\eea
as $n\to\infty$ for every $p\geq2$. 

Now, the estimate (\ref{201906021845}) with some positive constant $\ep_1$ follows from the estimates 
(\ref{22060241235}), (\ref{201906021843}) and (\ref{201906021842}). 
\qed\halflineskip
}
\begin{en-text}
Finally, the estimate (\ref{201905280550}) gives 
\bea\label{201906021844} &&
\sup_{\theta_3\in\overline{\Theta}_3}\Bigg|
\bbY^{(\ref{201906021814})}_n(\theta_1^*,\theta_3)
\nn\\&&\hspace{30pt}
+\frac{1}{2nh}\int_0^{nh}S(Z_t,\theta_1^*,\theta_3)^{-1}
\left[\left(\begin{array}{c}0\\
H(Z_t,\theta_3)-H(Z_t,\theta_3^*)
\end{array}\right)^{\otimes2}\right]dt
\Bigg|
\nn\\&&\hspace{10pt}
{\vblue\yeq O_{{\vblue L^\inftym}}(h^{1/2})+O_{{\vblue L^\inftym}}(n^{-1/2})\yeq O_{{\vblue L^\inftym}}(h^{1/2})\yeq o_{{\vblue L^\inftym}}(n^{-1/4}).}
\eea
Now 
{\tgreen(\ref{201906021845})} follows from (\ref{201906021842}), (\ref{201906021843}), (\ref{201906021844}) 
and $[A2]$ (ii) since ${\colorg\partial_3}^iH(z,\theta_1,\theta_3)$ $(i=0,1)$ are dominated by a polynomial in $z$ 
uniformly in $\theta_3$. 
%
%
{\vblue 
By using Sobolev's inequality and the mixing property with Lemma 4 of Yoshida \cite{yoshida2011polynomial}, 
there exists a positive constant $\ep_1$ such that 
\beas 
\sup_{\theta_3\in\ol{\Theta}_3}
\bigg|\frac{1}{2nh}\int_0^{nh}S(Z_t,\theta_1^*,\theta_3)^{-1}
\left[\left(\begin{array}{c}0\\
H(Z_t,\theta_3)-H(Z_t,\theta_3^*)
\end{array}\right)^{\otimes2}\right]dt
-\bbY^{(3)}(\theta_3)
\Bigg|
&=&
O_{L^\inftym}(n^{-\ep_1})
\eeas
as $n\to\infty$. 
Here we used the conditions that $nh\geq n^{\ep_0}$ for large $n$. 
}
\qed\halflineskip
\end{en-text}

\subsection{Random fields}
{\vblue 
To avoid repeatedly {\xred writing} similar big formulas, it is convenient to introduce several symbols of random fields. 
Let 
\beas 
\widetilde{\cald}_j(\theta_1',\theta_2',\theta_3')
&=&
\cald_j(\theta_1',\theta_2',\theta_3')+\widetilde{\cald}_j(\theta_1^*,\theta_2^*,\theta_3^*)
-\cald_j(\theta_1^*,\theta_2^*,\theta_3^*), 
\eeas
where 
$
\widetilde{\cald}_j(\theta_1^*,\theta_2^*,\theta_3^*)
$
is given in (\ref{2206290639}). 
The following random fields depend on $n$. 
\beas
\Psi_{3,1}(\theta_1,\theta_3,\theta_1',\theta_2',\theta_3')
&=&
\sum_{j=1}^nS(Z_\tjm,\theta_1,\theta_3)^{-1}
\bigg[\cald_j(\theta_1',\theta_2',\theta_3'),\>
\left(\begin{array}{c}0\\
\partial_3 H(Z_\tjm,\theta_3)
\end{array}\right)\bigg],
\eeas
\beas
\widetilde{\Psi}_{3,1}(\theta_1,\theta_3,\theta_1',\theta_2',\theta_3')
&=&
\sum_{j=1}^nS(Z_\tjm,\theta_1,\theta_3)^{-1}
\bigg[\widetilde{\cald}_j(\theta_1',\theta_2',\theta_3'),\>
\left(\begin{array}{c}0\\
\partial_3 H(Z_\tjm,\theta_3)
\end{array}\right)\bigg],
\eeas
\beas
\Psi_{3,2}(\theta_1,\theta_2,\theta_3,\theta_1',\theta_2',\theta_3')
&=&
\sum_{j=1}^nS(Z_\tjm,\theta_1,\theta_3)^{-1}
\bigg[\cald_j(\theta_1',\theta_2',\theta_3'),\>
\left(\begin{array}{c}0\\
2^{-1}h\partial_3 L_H(Z_\tjm,\theta_1,\theta_2,\theta_3)
\end{array}\right)\bigg],
\eeas
\beas
\Psi_{3,3}(\theta_1,\theta_3,\theta_1',\theta_2',\theta_3')
&=&
\half\sum_{j=1}^n\big(S^{-1}(\partial_3S)S^{-1}\big)(Z_\tjm,\theta_1,\theta_3)
\big[\cald_j(\theta_1',\theta_2',\theta_3')^{\otimes2}-S(Z_\tjm,\theta_1,\theta_3)\big],
\eeas
\beas
\Psi_{33,1}(\theta_1,\theta_2,\theta_3)
&=&
-\sum_{j=1}^nS(Z_\tjm,\theta_1,\theta_3)^{-1}
\left[
\left(\begin{array}{c}0\\
\bigg\{\begin{array}{c}
\partial_3H(Z_\tjm,\theta_3)
\\
+2^{-1}h\partial_3L_H(Z_\tjm,\theta_1,\theta_2,\theta_3)
\end{array}\bigg\}
\end{array}\right)^{\otimes2}
\right],
\eeas
\beas
\Psi_{33,2}(\theta_1,\theta_2,\theta_3,\theta_1',\theta_2',\theta_3')
&=&
\sum_{j=1}^nS(Z_\tjm,\theta_1,\theta_3)^{-1}
\bigg[\cald_j(\theta_1',\theta_2',\theta_3')
\\&&\qquad\otimes
\left(\begin{array}{c}0\\
\partial_3^2 H(Z_\tjm,\theta_3)
+2^{-1}h\partial_3^2 L_H(Z_\tjm,\theta_1,\theta_2,\theta_3)
\end{array}\right)\bigg],
\eeas
\beas
\Psi_{33,3}(\theta_1,\theta_3)
&=&
-\half\sum_{j=1}^n
\big\{\big(S^{-1}(\partial_3S)S^{-1}\big)(Z_\tjm,\theta_1,\theta_3)
\big[\partial_3S(Z_\tjm,\theta_1,\theta_3)\big]\big\},
\eeas
\beas
\Psi_{33,4}(\theta_1,\theta_2,\theta_3,\theta_1',\theta_2',\theta_3')
&=&
- {\colorr2}\sum_{j=1}^n
S^{-1}(\partial_3S)S^{-1}(Z_\tjm,\theta_1,\theta_3)
\Bigg[\cald_j(\theta_1',\theta_2',\theta_3')
\\&&\hspace{30pt}\otimes
\left(\begin{array}{c}0\\
\partial_3H(Z_\tjm,\theta_3)+2^{-1}h\partial_3L_H(Z_\tjm,\theta_1,\theta_2,\theta_3)
\end{array}\right)
\Bigg],
\eeas
\beas
\Psi_{33,5}(\theta_1,\theta_3,\theta_1',\theta_2',\theta_3')
&=&
\half \sum_{j=1}^n\partial_3\big\{
\big(S^{-1}(\partial_3S)S^{-1}\big)(Z_\tjm,\theta_1,\theta_3)\big\}
\bigg[\cald_j(\theta_1',\theta_2',\theta_3')^{\otimes2}-S(Z_\tjm,\theta_1,\theta_3)\bigg].
\eeas
} 

\begin{en-text}
{\tgreen
\begin{theorem}
Suppose that $[A1]$ with $(i_A,j_A,i_B,j_B,i_H,j_H)=(1,1,2,1,3,2)$, $[A2]$, 
$[A3]$ $(iii)$ and $[A4]$ are satisfied. Then 
\beas 
n^{1/2}h^{-1/2}\big(\wh{\theta}_3^0-\theta_3^*\big)
- \Gamma_{33}^{-1}M_n^{(3)}
&\to^p& 0
\eeas
as $n\to\infty$. In particular, 
\beas 
n^{1/2}h^{-1/2}\big(\wh{\theta}_3^0-\theta_3^*\big)
&\to^d&
N(0,\Gamma_{33}^{-1})
\eeas
as $n\to\infty$. 
\end{theorem}
}
\end{en-text}
\begin{en-text}
{\tred 
Now Theorem \ref{201906031711} follows from 
Lemmas \ref{201906031714} and \ref{201906031658}. 
}
\end{en-text}

\subsection{{\vblue CLT and $L^\inftym$-boundedness of $n^{-1/2}h^{1/2}\>\partial_3\bbH_3^{(3)}(\theta_3^*)$}}
\begin{en-text}
We will prove 
\beas 
\wh{\theta}_3^0-\theta_3^* &=& \xout{{\colorr O_p(n^{-1/2}h^{1/2})}}{\colorr O_p(n^{-1/2})}
\eeas
as a preliminary estimate of the error. 
\end{en-text}
\begin{en-text}
{\tgreen 
Let 
\bea\label{201906190133} 
\Gamma_{33}
&=&
\int S(z,\theta_1^*,\theta_3^*)^{-1}
\left[
\left(\begin{array}{c}0\\
\partial_3H(z,\theta_3^*)
\end{array}\right)^{\otimes2}
\right]
\nu(dz)
\nn\\&=&
\int 12V(z,\theta_1^*,\theta_3^*)^{-1}\big[\big(\partial_3H(z,\theta_3^{{\colorg*}})\big)^{\otimes2}\big]\nu(dz)
\nn\\&=&
\int12\partial_3H(z,\theta_3^*)^\star V(z,\theta_1^*,\theta_3^*)^{-1}\partial_3H(z,\theta_3^*)\nu(dz).
\eea
}
\end{en-text}

{\vblue 
The following algebraic result is simple but crucial. 
\begin{lemma}\label{2206241251}
\bea\label{201906091302} 
S(Z_\tjm,\theta_1,\theta_3^*)^{-1}
\left[\cald_j(\theta_1,\theta_2,\theta_3^*)-\cald_j(\theta_1,\theta_2^*,\theta_3^*),\>
\left(\begin{array}{c}0\\
\partial_3H(Z_\tjm,\theta_3^*)
\end{array}\right)
\right]
&=&
0. 
\eea
\end{lemma}
\proof
Apply linear algebra, with the equality
\beas\label{2206251356}
\cald_j(\theta_1,\theta_2,\theta_3^*)
-\cald_j(\theta_1,\theta_2^*,\theta_3^*)
&=& 
-h^{1/2}
\left(\begin{array}{c}
A(Z_\tjm,\theta_2)-A(Z_\tjm,\theta_2^*)\y
2^{-1}H_x(Z_\tjm,\theta_3^*)\big[A(Z_\tjm,\theta_2)-A(Z_\tjm,\theta_2^*)\big]
\end{array}\right). 
\eeas
\qed
}

\begin{lemma}\label{201906031714}
Suppose that $[A1]$ with $\bbi=(1,1,2,1,3,1)$, $[A2]$  
{\tred and} $[A4]$ {\vredd $(i)$} are satisfied. 
Then 
{\vblue
\bea\label{2206241431}
n^{-1/2}h^{1/2}\>\partial_3\bbH_n^{(3)}(\theta_3^*)
-M^{(3)}_n
&=&
o_{{\vblue L^\inftym}}(1)
\eea
as $n\to\infty$. 
$M^{(3)}_n=O_{L^\inftym}(1)$ and 
$n^{-1/2}h^{1/2}\>\partial_3\bbH_n^{(3)}(\theta_3^*)=O_{L^\inftym}(1)$ as $n\to\infty$. 
Moreover, 
\bea\label{2206241524}
M_n^{(3)} &\to^d & N_{\sfp_3}(0,\Gamma_{33})
\eea
as $n\to\infty$. 
}
\end{lemma}
\proof 
{\vblue 
Prepare the following random vectors: 
\bea\label{201906030100}
R^{(\ref{201906030100})}_n(\wh{\theta}_1^0,\wh{\theta}_2^0)
&=&
n^{-1/2}\Psi_{3,1}(\wh{\theta}_1^0,\theta_3^*,\wh{\theta}_1^0,\wh{\theta}_2^0,\theta_3^*),
\eea
\bea\label{201906022139}
R^{(\ref{201906022139})}_n(\wh{\theta}_1^0,\wh{\theta}_2^0)
&=&
n^{-1/2}\Psi_{3,2}(\wh{\theta}_1^0,\wh{\theta}_2^0,\theta_3^*,\wh{\theta}_1^0,\wh{\theta}_2^0,\theta_3^*)
\nn\\&&
\eea
and 
\bea\label{201906030102}
R^{(\ref{201906030102})}_n(\wh{\theta}_1^0,\wh{\theta}_2^0)
&=&
n^{-1/2}h^{1/2}\Psi_{3,3}(\wh{\theta}_1^0,\theta_3^*,\wh{\theta}_1^0,\wh{\theta}_2^0,
\theta_3^*).
\eea

Thanks to the formula (\ref{201906091302}) of Lemma \ref{2206241251}, we have 
\bea\label{2206241353}
R^{(\ref{201906030100})}_n(\wh{\theta}_1^0,\wh{\theta}_2^0)
&=&
n^{-1/2}\Psi_{3,1}(\wh{\theta}_1^0,\theta_3^*,\wh{\theta}_1^0,\theta_2^*,\theta_3^*). 
\eea
We apply Lemma \lfive (b) of \cite{gloter2021adaptive} and {\vredd $[A4]$ (i)}, and next  
the results in Lemmas \lfour and \lfive (a) of \cite{gloter2021adaptive} 
to show
\bea\label{2206160514}
n^{-1/2}\Psi_{3,1}(\wh{\theta}_1^0,\theta_3^*,\wh{\theta}_1^0,\theta_2^*,\theta_3^*)
&=&
{\colorg n^{-1/2}}
\Psi_{3,1}(\wh{\theta}_1^0,\theta_3^*,\theta_1^*,\theta_2^*,\theta_3^*)
+O_{{\vblue L^\inftym}}({\colorr h^{{\colorr1/2}}})
\\&=&
{\colorg n^{-1/2}}
\wt{\Psi}_{3,1}(\wh{\theta}_1^0,\theta_3^*,\theta_1^*,\theta_2^*,\theta_3^*)+o_{{\vblue L^\inftym}}(1)
\label{201906031505}
\eea
due to $(nh^2)^{1/2}=o(1)$ by assumption. 
{\vredd Condition $[A4]$ (i) was used loosely in that only the property that $\wh{\theta}_1^0-\theta_1^*=o_{L^\inftym}(n^{-1/2}h^{-1/2})$ was necessary in (\ref{2206160514}), 
in order to finally get (\ref{201906031505}). 
We do not use $[A4]$ (ii). }

To estimate ${\colorg n^{-1/2}}\wt{\Psi}_{3,1}(\wh{\theta}_1^0,\theta_3^*,\theta_1^*,\theta_2^*,\theta_3^*)$, 
we consider the random field 
\bea\label{201906091137}
\Phi_n^{(\ref{201906091137})}(u_1)
&=&
{\colorg n^{-1/2}\big\{}
\widetilde{\Psi}_{3,1}(\theta_1^*+r_nu_1,\theta_3^*,\theta_1^*,\theta_2^*,\theta_3^*)
-\widetilde{\Psi}_{3,1}(\theta_1^*,\theta_3^*,\theta_1^*,\theta_2^*,\theta_3^*)
{\colorg\big\}}
\eea
on $\{u_1\in\bbR^{\sfp_1};\>|u_1|<1\}$ for any sequence of positive numbers $r_n\to0$.  
It is easy to show 
\beas 
\sup_{u_1: |u_1|<1}\big\|\partial_{u_1}^i\Phi_n^{(\ref{201906091137})}(u_1)\big\|_p
&=& 
o(1)
\eeas
for $i=0,1$ and every $p>1$, 
by using the orthogonality. 
Therefore, Sobolev's inequality gives 
\beas 
\sup_{u_1: |u_1|<1}|\Phi_n^{(\ref{201906091137})}(u_1)| 
&=&
o_{{\vblue L^\inftym}}(1)
\eeas
as $n\to\infty$. 
In particular, 
\bea\label{2206241352}
{\colorg n^{-1/2}}\wt{\Psi}_{3,1}(\wh{\theta}_1^0,\theta_3^*,\theta_1^*,\theta_2^*,\theta_3^*)
&=&
{\colorg n^{-1/2}}\wt{\Psi}_{3,1}(\theta_1^*,\theta_3^*,\theta_1^*,\theta_2^*,\theta_3^*)+
o_{{\vblue L^\inftym}}(1)
\eea
as $n\to\infty$. {\vredd Here we used $[A4]$ (i) loosely.} 
Thus, it follows 
from (\ref{2206241353}), (\ref{201906031505}) and (\ref{2206241352}) that 
\beas\label{201906091150}
R^{(\ref{201906030100})}_n(\wh{\theta}_1^0,\wh{\theta}_2^0)
&=&
n^{-1/2}\widetilde{\Psi}_{3,1}(\theta_1^*,\theta_3^*,\theta_1^*,\theta_2^*,\theta_3^*)+o_{{\vblue L^\inftym}}(1),
\eeas
and this easily implies 
\bea\label{2206241404} 
R^{(\ref{201906030100})}_n(\wh{\theta}_1^0,\wh{\theta}_2^0)
&=&
M_n^{(3)}+o_{{\vblue L^\inftym}}(1). 
\eea

%
In a similar fashion, it is possible to verify 
\bea\label{2206241427}
R^{(\ref{201906030102})}_n(\wh{\theta}_1^0,\wh{\theta}_2^0)
&=&
n^{-1/2}h^{\colorg 1/2}\Psi_{3,3}(\theta_1^*,\theta_3^*,\theta_1^*,\theta_2^*,\theta_3^*)
+O_{{\vblue L^\inftym}}(h^{1/2})
\nn\\&=&
O_{{\vblue L^\inftym}}({\colorr h^{1/2}}). 
\eea
On the other hand, %
by (\ref{0406241100}) of Lemma \ref{2206241411}, 
it is easy to show that 
\bea\label{2206241428}
R^{(\ref{201906022139})}_n(\wh{\theta}_1^0,\wh{\theta}_2^0)&=&O_{L^\inftym}(n^{1/2}h). 
\eea

From (\ref{201906041919}) and (\ref{201906030041}), we have 
\bea\label{201906030105}
n^{-1/2}h^{1/2}\>\partial_3\bbH_n^{(3)}(\theta_3^*)
&=&
R^{(\ref{201906030100})}_n(\wh{\theta}_1^0,\wh{\theta}_2^0)+R^{(\ref{201906022139})}_n(\wh{\theta}_1^0,\wh{\theta}_2^0)
+R^{(\ref{201906030102})}_n(\wh{\theta}_1^0,\wh{\theta}_2^0).
\eea
Combining (\ref{2206241404}), (\ref{2206241427}) and (\ref{2206241428}) with the decomposition (\ref{201906030105}), 
we obtain (\ref{2206241431}). 
{\vredd By the Burkholder-Davis-Gundy inequality,} 
we have $M^{(3)}_n=O_{L^\inftym}(1)$, and hence 
$n^{-1/2}h^{1/2}\>\partial_3\bbH_n^{(3)}(\theta_3^*)=O_{L^\inftym}(1)$ as $n\to\infty$. 

In view of the relation (\ref{2206241549}), we can show 
the convergence (\ref{2206241524}) as a consequence of the martingale central limit theorem applied to $M_n^{(3)}$ of (\ref{2206241544}). 
\qed\halflineskip
} 

\begin{en-text}
\beas
n^{-1/2}h\>\partial_3\bbH_n^{(2,3)}(\wh{\theta}_2^0,\theta_3^*)
&=&
O_p(h^{1/2})+O_p(n^{1/2}h\sqrt{h})
\\&=&
O_p(h^{1/2})
\eeas
where 
\bea
R^{(\ref{201906022139})}_n(\wh{\theta}_1,\wh{\theta}_2^0)
&=&
n^{-1/2}h^{1/2}\sum_{j=1}^n\wh{S}(Z_\tjm,\theta_3^*)^{-1}
\left[\cald_j(\wh{\theta}_2^0,\theta_3^*),\>
\left(\begin{array}{c}0\\
2^{-1}h\partial_3L_H(Z_\tjm,\wh{\theta}_1,\wh{\theta}_2^0,\theta_3^*)
\end{array}\right)
\right]
\nn\\&=&
n^{-1/2}h^{1/2}\sum_{j=1}^n\wh{S}(Z_\tjm,\theta_3^*)^{-1}
\left[\cald_j(\theta_2^*,\theta_3^*),\>
\left(\begin{array}{c}0\\
2^{-1}h\partial_3L_H(Z_\tjm,\wh{\theta}_1,\wh{\theta}_2^0,\theta_3^*)
\end{array}\right)
\right]
\nn\\&&
+O_p(n^{1/2}h^2)
\nn\\&=&
n^{-1/2}h^{1/2}\sum_{j=1}^n\wh{S}(Z_\tjm,\theta_3^*)^{-1}
\left[
\left[\begin{array}{c}
h^{-1/2}\int_\tjm^\tj B(Z_t,\theta_2^*)dw_t\y
h^{-3/2}K(Z_\tjm,\theta_1^*,\theta_3^*)\zeta_j
\end{array}
\right]\right.
\nn\\&&
\hspace{5cm}\otimes
\left.
\left[\begin{array}{c}0\\
2^{-1}h\partial_3L_H(Z_\tjm,\wh{\theta}_1,\wh{\theta}_2^0,\theta_3^*)
\end{array}\right]
\right]
+O_p(n^{1/2}h^2)
\nn\\&=&
O_p(n^{1/2}h\times\sqrt{h})\yeq o_p(\sqrt{h})
\eea
\end{en-text}

\subsection{{\vblue Convergence of $n^{-1}h\>\partial_3^2\bbH_n^{(3)}(\theta_3)$}}
{\vblue
We are writing $\bbi=(i_A,j_A,i_B,j_B,i_H,j_H)$.}
\begin{lemma}\label{201906031658}
Suppose that $[A1]$ with $\bbi=(1,1,{\colorg 2},1,3,{\vblue3})
$, $[A2]$ and $[A4]$ {\vredd $(i)$} are satisfied. Then, 
{\vblue for some number $n_0$, 
\bea\label{2206230425} 
 \sup_{n\geq n_0}\bigg\|
 \sup_{u_3\in \bbU_n^{(3)}(1)}\left|
n^{-1}h\>\partial_3^2\bbH_n^{(3)}(\theta_3^*+\delta u_3)
-n^{-1}h\>\partial_3^2\bbH_n^{(3)}(\theta_3^*)
\right|\bigg\|_p
&=&
O(\delta)
\eea
as $\delta\down0$ for every $p>1$, where $\bbU_n^{(3)}(r)= \{u_3\in \bbR^{\sfp_3};
\>\theta_3^*+\delta u_3\in\Theta_3,|u_3|<r\}$. 
}
%
{\vblue Moreover, for some positive constant $\ep_2$, it holds that 
\bea\label{2206230443} 
n^{\ep_2}\left|
n^{-1}h\>\partial_3^2\bbH_n^{(3)}(\theta_3^*)
+\Gamma_{33}
\right|
&=&
O_{L^\inftym}(1)
\eea
as $n\to\infty$. }
\end{lemma}
\proof 
\begin{en-text}
Since 
\beas &&
n^{-1}h\>\partial_3\bbH_n^{(3)}(\theta_3)
\\&=&
n^{-1}h^{1/2}\sum_{j=1}^n\wh{S}(Z_\tjm,\theta_3)^{-1}
\left[\cald_j(\wh{\theta}_1^0,\wh{\theta}_2^0,\theta_3),\>
\left[\begin{array}{c}0\\
\partial_3H(Z_\tjm,\theta_3)+2^{-1}h\partial_3L_H(Z_\tjm,\wh{\theta}_1^0,\wh{\theta}_2^0,\theta_3)
\end{array}\right]
\right]
\\&&
+\half n^{-1}h\sum_{j=1}^n
\big(\wh{S}^{-1}(\partial_3\wh{S})\wh{S}^{-1}\big)(Z_\tjm,\theta_3)
\big[\cald_j(\wh{\theta}_1^0,\wh{\theta}_2^0,\theta_3)^{\otimes2}-\wh{S}(Z_\tjm,\theta_3)\big],
\eeas
\end{en-text}
{\vblue 
We use the following decomposition of $n^{-1}h\>\partial_3^2\bbH_n^{(3)}(\theta_3)$: }
\bea\label{2206230451}
n^{-1}h\>\partial_3^2\bbH_n^{(3)}(\theta_3)
&=&
n^{-1}\Psi_{33,1}(\wh{\theta}_1^0,\wh{\theta}_2^0,\theta_3)
+n^{-1}h^{1/2}\Psi_{33,2}(\wh{\theta}_1^0,\wh{\theta}_2^0,\theta_3,\wh{\theta}_1^0,\wh{\theta}_2^0,\theta_3)
\nn\\&&
+n^{-1}h\Psi_{33,3}(\wh{\theta}_1^0,\theta_3)
+n^{-1}h^{1/2}\Psi_{33,4}(\wh{\theta}_1^0,\wh{\theta}_2^0,\theta_3,\wh{\theta}_1^0,\wh{\theta}_2^0,\theta_3)
\nn\\&&
+n^{-1}h\Psi_{33,5}(\wh{\theta}_1^0,\theta_3,\wh{\theta}_1^0,\wh{\theta}_2^0,\theta_3)
.
\eea
{\vblue 
Fix a sufficiently large integer $n_0$ such that $\bbU_n^{(3)}(\delta)=U(\theta_3^*,\delta)=\{\theta_3\in\bbR^{\sfp_3};|\theta_3-\theta_3^*|<\delta\}$ for all $n\geq n_0$ and $\delta\leq1$. 
Denote by $p$ an arbitrary number bigger than $1$. 
We recall the estimate (\ref{2206230411}) of Lemma \ref{2206230342}.  

By $[A1]$ for $\bbi=(0,0,0,0,2,2)$, we obtain
\bea\label{2206221804}
\sup_{n\geq n_0}\bigg\|\sup_{\theta_3\in U(\theta_3^*,\delta)}
\big|n^{-1}\Psi_{33,1}(\wh{\theta}_1^0,\wh{\theta}_2^0,\theta_3)-n^{-1}\Psi_{33,1}(\wh{\theta}_1^0,\wh{\theta}_2^0,\theta_3^*)\big|\bigg\|_p
&=& 
O(\delta)
\eea
as $\delta\down0$. 
Apply Lemma \lsix under $[A1]$ for $\bbi=(0,0,0,0,2,3)$ and $[A2]$ (i), 
\bea\label{2206221819}
\sup_{n\geq n_0}\bigg\|\sup_{\theta_3\in U(\theta_3^*,\delta)}
n^{-1}h^{1/2}\Psi_{33,2}(\wh{\theta}_1^0,\wh{\theta}_2^0,\theta_3,\wh{\theta}_1^0,\wh{\theta}_2^0,\theta_3)
-n^{-1}h^{1/2}\Psi_{33,2}(\wh{\theta}_1^0,\wh{\theta}_2^0,\theta_3^*,\wh{\theta}_1^0,\wh{\theta}_2^0,\theta_3^*)
\big|\bigg\|_p
&=& 
O(\delta).
\nn\\&&
\eea
Remark that the random field $\Psi_{33,2}$ involves $\partial_3^2\partial_z^iH$ for $i=1,2$. 

With $[A1]$ for $\bbi=(0,0,0,0,1,2)$, it is easy to see 
\bea\label{2206230223}
\sup_{n\geq n_0}\bigg\|\sup_{\theta_3\in U(\theta_3^*,\delta)}
n^{-1}h\Psi_{33,3}(\wh{\theta}_1^0,\theta_3)
-n^{-1}h\Psi_{33,3}(\wh{\theta}_1^0,\theta_3^*)
\big|\bigg\|_p
&=& 
\sup_{n\geq n_0}h\times O(\delta)\yeq O(\delta),
\eea
and 
with $[A1]$ for $\bbi=(0,0,0,0,2,2)$ and Lemma \lsix of \cite{gloter2021adaptive}, 
\bea\label{2206230330}
\sup_{n\geq n_0}\bigg\|\sup_{\theta_3\in U(\theta_3^*,\delta)}
n^{-1}h^{1/2}\Psi_{33,4}(\wh{\theta}_1^0,\wh{\theta}_2^0,\theta_3,\wh{\theta}_1^0,\wh{\theta}_2^0,\theta_3)
-n^{-1}h^{1/2}\Psi_{33,4}(\wh{\theta}_1^0,\wh{\theta}_2^0,\theta_3^*,\wh{\theta}_1^0,\wh{\theta}_2^0,\theta_3^*)
\big|\bigg\|_p
&=& 
O(\delta).
\nn\\&&
\eea
Due to Lemma \lsix of \cite{gloter2021adaptive} and $[A1]$ for $\bbi=(0,0,0,0,2,3)$,  
\bea\label{22062300404}
\sup_{n\geq n_0}\bigg\|\sup_{\theta_3\in U(\theta_3^*,\delta)}
n^{-1}h\Psi_{33,5}(\wh{\theta}_1^0,\theta_3,\wh{\theta}_1^0,\wh{\theta}_2^0,\theta_3)
-n^{-1}h\Psi_{33,5}(\wh{\theta}_1^0,\theta_3^*,\wh{\theta}_1^0,\wh{\theta}_2^0,\theta_3^*)
\big|\bigg\|_p
&=& 
O(\delta).
\nn\\&&
\eea
The random field $\Psi_{33,5}$ initially has the second-order derivative $\partial_3^2H_x$, and due to this, 
we used $\partial_3^3H_x$ to show (\ref{2206230425}). 
Additionally to the above mentioned conditions, we assumed $[A1]$ with $\bbi=(1,1,2,1,3,{\vredd1})$, that validates the estimate (\ref{2206230411}) of Lemma \ref{2206230342}. 
Thus, that $[A1]$ for $\bbi=(1,1,2,1,3,3)$ is required to conclude (\ref{2206230425}) 
from (\ref{2206221804})-
(\ref{22062300404}) and (\ref{2206230451}).

Let us prove (\ref{2206230443}). 
\begin{en-text}
For $\cald_j(\wh{\theta}_1^0,\wh{\theta}_2^0,\theta_3)$ in $\Psi_{33,2}$, $\Psi_{33,4}$ and $\Psi_{33,5}$ of (\ref{2206230451}), 
we use Lemma \lfive (b) to replace $\wh{\theta}_i^0$ by $\theta_i^*$ 
for $i=1,2$, and Lemma \lthree (b). 
\end{en-text}
The true value $\theta_3^*$ plugged into $\theta_3$ of the formula (\ref{2206230443}), 
the last four terms on the right-hand side of (\ref{2206230451}), evaluated at $\theta_3^*$, are of 
$O_{L^\inftym}(h^{1/2})$, $O_{L^\inftym}(h)$, $O_{L^\inftym}(h^{1/2})$ and $O_{L^\inftym}(h)$, respectively, 
since 
$
\sup_{(\theta_1,\theta_2)\in\overline{\Theta}_1\times\overline{\Theta}_2}
\big|\cald_j(\theta_1,\theta_2,\theta_3^*)\big|
=
O_{L^\inftym}(1)
$
by Lemma \ref{2206241411}. 
Together with these evaluations, a simple estimate for $n^{-1}\Psi_{33,1}(\wh{\theta}_1^0,\wh{\theta}_2^0,\theta_3^*)$ gives 
\bea\label{2206241504}
n^{-1}h\>\partial_3^2\bbH_n^{(3)}(\theta_3^*)
&=&
-n^{-1}\sum_{j=1}^nS(Z_\tjm,\theta_1^*,\theta_3^*)^{-1}
\left[
\left(\begin{array}{c}0\\
\partial_3H(Z_\tjm,\theta_3^*)
\end{array}\right)^{\otimes2}
\right]
+r_n^{(\ref{201906031646})}(\theta_3^*)
\nn\\&=&
-n^{-1}\sum_{j=1}^n
12V(Z_\tjm,\theta_1^*,\theta_3^*)^{-1}\big[\big(\partial_3H(Z_\tjm,\theta_3^{{\colorg*}})\big)^{\otimes2}\big]
+r_n^{(\ref{201906031646})}(\theta_3^*)
\eea
for some random matrix $r_n^{(\ref{201906031646})}(\theta_3^*)$ such that 
\bea\label{201906031646}
\big|r_n^{(\ref{201906031646})}(\theta_3^*)\big|
&=&
O_{L^\inftym}(n^{-1/2})+O_{L^\inftym}(h^{1/2})\yeq o_{L^\inftym}(n^{-1/4}). 
\eea
{\vredd Here we used the condition that $\wh{\theta}_1^0-\theta_1^*=O_{L^\inftym}(n^{-1/2})$, though the used convergence rate can be relaxed. 
Condition $[A4]$ (ii) was not used since $\wh{\theta}_2^0$ is wrapped by 
the function of order $h$, in $n^{-1}\Psi_{33,1}(\wh{\theta}_1^0,\wh{\theta}_2^0,\theta_3^*)$. }

It is possible to show that 
\bea\label{2206241505}
n^{-1}\sum_{j=1}^n
12V(Z_\tjm,\theta_1^*,\theta_3^*)^{-1}\big[\big(\partial_3H(Z_\tjm,\theta_3^{{\colorg*}})\big)^{\otimes2}\big]
-\Gamma_{33}
&=&
O_{L^\inftym}(n^{-\ep_0/2}), 
\eea
in a similar manner to (\ref{22060241235}) but simpler. 
Consequently, the estimate (\ref{2206230443}) is verified by (\ref{2206241504}), (\ref{201906031646}) and (\ref{2206241505}). 
\qed\halflineskip
}
%
%

{\vblue 
\subsection{Proof of Theorem \ref{201906031711}}
The simplified QLA will be applied to the random field $\bbH_n^{(3)}(\theta_3)$ as $\bbH_n(\theta)$ in Yoshida \cite{yoshida2021simplified}. 
Condition $[U1]$ of \cite{yoshida2021simplified} is nothing but $[A3]$ (iii). 
Condition $[T2]$ of \cite{yoshida2021simplified} is verified in this situation 
by using Lemmas \ref{201906021825}, \ref{201906031714} and \ref{201906031658}. 
The convergence (3.3) required in \cite{yoshida2021simplified} has already been obtained in Lemma \ref{201906031714}. 
Therefore, we can apply Theorem 3.5 of \cite{yoshida2021simplified} in order to conclude the results. 
\qed
}

{\vblue 
\section{Proof of Theorem \ref{2206261240}}\label{2206261259}
\subsection{Convergence of the random field $\bbY^{(J,1)}_n(\theta_1)$}
To apply the theory of QLA, we set 
\beas 
\bbY^{(J,1)}_n(\theta_1)
&=& 
n^{-1}\big(\bbH^{(1)}_n(\theta_1)-\bbH^{(1)}_n(\theta_1^*)\big). 
\eeas
Recall (\ref{2206251502}): 
\beas
\bbY^{(J,1)}(\theta_1)
&=& 
-\frac{1}{2}\int\bigg\{ 
\big(S(z,\theta_1,\theta_3^*)^{-1}-S(z,\theta_1^*,\theta_3^*)^{-1}\big)\big[S(z,\theta_1^*,\theta_3^*)\big]
 +\log\frac{\det S(z,\theta_1,\theta_3^*)}{\det S(z,\theta_1^*,\theta_3^*)}
 \bigg\}\nu(dz).
 \eeas

We have the following result on uniform convergence of $\bbY^{(J,1)}_n$ with a convergence rate. 

\begin{lemma}\label{22062511258} 
Suppose that $[A1]$ with $\bbi=(1,1,2,1,3,1)$, $[A2]$ $(i)$, {\vredd $[A4]$ $(iii)$}. Then, 
there exists a positive constant $\ep_1$ such that 
\beas
\sup_{\theta_1\in\ol{\Theta}_1}\big|
\bbY^{(J,1)}_n(\theta_1)-\bbY^{{\vredd(J,1)}}(\theta_1)\big|
&=&
O_{L^\inftym}(n^{-\ep_1})
\eeas
as $n\to\infty$. 
\end{lemma}
\proof 
Observe that 
\bea\label{2206251322} &&
\sup_{\theta_1\in\ol{\Theta}_1}\big|\cald_j(\theta_1,\wh{\theta}_2^0,\wh{\theta}_3^0)-\cald_j(\theta_1,\theta_2^*,\theta_3^*)\big|
\nn\\&=&
\sup_{\theta_1\in\ol{\Theta}_1}\big|\cald_j(\theta_1,\wh{\theta}_2^0,\theta_3^*)-\cald_j(\theta_1,\theta_2^*,\theta_3^*)\big|
+O_{L^\inftym}(n^{-1/2})
\quad(\text{Lemma \lsix of \cite{gloter2021adaptive} and $[A4]$ (iii)})
\nn\\&=&
O_{L^\inftym}(h^{1/2})+O_{L^\inftym}(n^{-1/2})
\quad(\text{boundedness of $\Theta_2$})
\nn\\&=&
O_{L^\inftym}(n^{-1/4}).
\eea
With (\ref{2206251322}) and Lemmas \ref{201906021729} and \ref{2206241411}, we obtain 
\bea\label{2206251323}
\sup_{\theta_1\in\ol{\Theta}_1}\big|\bbY^{(J,1)}_n(\theta_1)-\bbY^{(\ref{2206251332})}_n(\theta_1)\big|
&=& 
O_{L^\inftym}(n^{-1/4})
\eea
as $n\to\infty$ for 
\bea\label{2206251332}
\bbY^{(\ref{2206251332})}_n(\theta_1)
&=&
-\frac{1}{2n}\sum_{j=1}^n\bigg\{
S(Z_\tjm,\theta_1,\theta_3^*)^{-1}\big[\cald_j(\theta_1,\theta_2^*,\theta_3^*)^{\otimes2}-\cald_j(\theta_1^*,\theta_2^*,\theta_3^*)^{\otimes2}\big]
\nn\\&&\hspace{50pt}
 +\log\frac{\det S(Z_\tjm,\theta_1,\theta_3^*)}{\det S(Z_\tjm,\theta_1^*,\theta_3^*)}
 \bigg\}. 
\eea
Easily, 
\bea\label{2206251343}
\bbY^{(\ref{2206251332})}_n(\theta_1)
&=&
\bbY^{(\ref{2206251417})}_n(\theta_1)+\bbY^{(\ref{2206251418})}_n(\theta_1)+\bbY^{(\ref{2206251419})}_n(\theta_1),
\eea
where 
\bea\label{2206251417}
\bbY^{(\ref{2206251417})}_n(\theta_1)
&=& 
-\frac{1}{2n}\sum_{j=1}^n\bigg\{ 
\big(S(Z_\tjm,\theta_1,\theta_3^*)^{-1}-S(Z_\tjm,\theta_1^*,\theta_3^*)^{-1}\big)\big[\cald_j(\theta_1^*,\theta_2^*,\theta_3^*)^{\otimes2}\big]
\nn\\&&\hspace{50pt}
 +\log\frac{\det S(Z_\tjm,\theta_1,\theta_3^*)}{\det S(Z_\tjm,\theta_1^*,\theta_3^*)}
 \bigg\},
 \eea
\bea\label{2206251418}
\bbY^{(\ref{2206251418})}_n(\theta_1)
&=&
-\frac{1}{2n}\sum_{j=1}^n
S(Z_\tjm,\theta_1,\theta_3^*)^{-1}\big[\big(\cald_j(\theta_1,\theta_2^*,\theta_3^*)-\cald_j(\theta_1^*,\theta_2^*,\theta_3^*)\big)^{\otimes2}\big]
\eea
and 
\bea\label{2206251419}
\bbY^{(\ref{2206251419})}_n(\theta_1)
&=&
-\frac{1}{n}\sum_{j=1}^n
S(Z_\tjm,\theta_1,\theta_3^*)^{-1}\big[\cald_j(\theta_1,\theta_2^*,\theta_3^*)-\cald_j(\theta_1^*,\theta_2^*,\theta_3^*),\cald_j(\theta_1^*,\theta_2^*,\theta_3^*)\big].
\eea
\begin{en-text}
\bea\label{2206251343}
\bbY^{(\ref{2206251332})}_n(\theta_1)
&=&
-\frac{1}{2n}\sum_{j=1}^n\bigg\{ 
S(Z_\tjm,\theta_1,\theta_3^*)^{-1}\big[\big(\cald_j(\theta_1,\theta_2^*,\theta_3^*)-\cald_j(\theta_1^*,\theta_2^*,\theta_3^*)\big)^{\otimes2}\big]
\nn\\&&\hspace{50pt}
+2S(Z_\tjm,\theta_1,\theta_3^*)^{-1}\big[\cald_j(\theta_1,\theta_2^*,\theta_3^*)-\cald_j(\theta_1^*,\theta_2^*,\theta_3^*),\cald_j(\theta_1^*,\theta_2^*,\theta_3^*)\big]
\nn\\&&\hspace{50pt}
+\big(S(Z_\tjm,\theta_1,\theta_3^*)^{-1}-S(Z_\tjm,\theta_1^*,\theta_3^*)^{-1}\big)\big[\cald_j(\theta_1^*,\theta_2^*,\theta_3^*)^{\otimes2}\big]
\nn\\&&\hspace{50pt}
 +\log\frac{\det S(Z_\tjm,\theta_1,\theta_3^*)}{\det S(Z_\tjm,\theta_1^*,\theta_3^*)}
 \bigg\}
\eea
\end{en-text}

By using the equality 
\beas\label{2206251427}
\cald_j(\theta_1,\theta_2^*,\theta_3^*)
-\cald_j(\theta_1^*,\theta_2^*,\theta_3^*)
&=& 
-2^{-1}h^{1/2}
\left(\begin{array}{c}
0\y
L_H(Z_\tjm,\theta_1,\theta_2^*,\theta_3^*)-A(Z_\tjm,\theta_1^*,\theta_2^*,\theta_3^*)
\end{array}\right)
\eeas
and 
Lemma \lthree (b) of \cite{gloter2021adaptive}, we obtain 
\bea\label{2206251440}
\sup_{\theta_1\in\ol{\Theta}_1}\big|\bbY^{(\ref{2206251418})}_n(\theta_1)\big|
&=&
O_{L^\inftym}(h)
\eea
and 
\bea\label{2206251441}
\sup_{\theta_1\in\ol{\Theta}_1}\big|\bbY^{(\ref{2206251419})}_n(\theta_1)\big|
&=&
O_{L^\inftym}(h^{1/2}). 
\eea

Define $\bbY^{(\ref{2206251512})}_n(\theta_1)$ by 
\bea\label{2206251512}
\bbY^{(\ref{2206251512})}_n(\theta_1)
&=& 
-\frac{1}{2n}\sum_{j=1}^n\bigg\{ 
\big(S(Z_\tjm,\theta_1,\theta_3^*)^{-1}-S(Z_\tjm,\theta_1^*,\theta_3^*)^{-1}\big)\big[S(Z_\tjm,\theta_1^*,\theta_3^*)\big]
\nn\\&&\hspace{50pt}
 +\log\frac{\det S(Z_\tjm,\theta_1,\theta_3^*)}{\det S(Z_\tjm,\theta_1^*,\theta_3^*)}
 \bigg\}
 \eea
 We have 
 \bea\label{2206251514}
 \sup_{\theta_1\in\ol{\Theta}_1}\big|\bbY^{(\ref{2206251417})}_n(\theta_1)-\bbY^{(\ref{2206251512})}_n(\theta_1)\big|
 &=&
 O_{L^\inftym}(n^{-1/2})
 \eea
 as a consequence of Sobolev's inequality for functions in $\theta_1$ and the Burkholder-Davis-Gundy inequality. 
Moreover, in the same fashion as the derivation of (\ref{22060241235}), 
by using Sobolev's inequality with the derivative $\partial_1B$ and the mixing property, we obtain the estimate
\bea\label{2206251456}
\sup_{\theta_1\in\ol{\Theta}_1}\big|
\bbY^{(\ref{2206251512})}_n(\theta_1)-\bbY^{{\vredd(J,1)}}(\theta_1)\big|
&=&
O_{L^\inftym}(n^{-\ep_0/2})
\eea
as $n\to\infty$. 

The proof is completed by 
integrating (\ref{2206251323}), (\ref{2206251343}), (\ref{2206251440}), (\ref{2206251441}), (\ref{2206251514}) and (\ref{2206251456}). 
\qed\halflineskip
}


{\vblue 
\subsection{Convergence of $n^{-1}\>\partial_1^2\bbH_n^{(1)}(\theta_1)$}
The following random fields will be used. 
\beas &&
{\vredd \Psi_{1,1}(\theta_1,\theta_3,\theta_1',\theta_2',\theta_3')}
\\&=&
\sum_{j=1}^nS(Z_\tjm,\theta_1,\theta_3)^{-1}
\left[\cald_j(\theta_1',\theta_2',\theta_3'),\>
\left(\begin{array}{c}0\\
2^{-1}\partial_1L_H(Z_\tjm,\theta_1,\theta_2,\theta_3)
\end{array}\right)
\right]
\\&=&
\sum_{j=1}^nS(Z_\tjm,\theta_1,\theta_3)^{-1}
\Bigg[\cald_j(\theta_1',\theta_2',\theta_3'),\>
\left(\begin{array}{c}0\\
{\fred4^{-1}}H_{xx}(Z_\tjm,\theta_3)[\partial_1C(Z_\tjm,\theta_1)]
\end{array}\right)\Bigg],
\eeas
\beas 
\Psi_{1,2}(\theta_1,\theta_3,\theta_1',\theta_2',\theta_3')
&=&
\half\sum_{j=1}^n
\big(S^{-1}(\partial_1S))S^{-1}\big)(Z_\tjm,\theta_1,\theta_3)
\big[\cald_j(\theta_1',\theta_2',\theta_3')^{\otimes2}-S(Z_\tjm,\theta_1,\theta_3)\big],
\eeas
{\vblue 
\beas 
{\vredd\Psi_{11,1}(\theta_1,\theta_3,\theta_1',\theta_3')}
&=&
\sum_{j=1}^nS(Z_\tjm,\theta_1,\theta_3)^{-1}
\left[
\left(\begin{array}{c}0\\
2^{-1}\partial_1L_H(Z_\tjm,\theta_1',\theta_2',\theta_3')
\end{array}\right)^{\otimes2}
\right]
\\&=&
\sum_{j=1}^nS(Z_\tjm,\theta_1,\theta_3)^{-1}
\left[
\left(\begin{array}{c}0\\
4^{-1}H_{xx}(Z_\tjm,\theta_3')[\partial_1C(Z_\tjm,\theta_1')]
\end{array}\right)^{\otimes2}
\right],
\eeas
}
{\vblue 
\beas &&
{\vredd\Psi_{11,2}(\theta_1,\theta_3,\theta_1',\theta_2',\theta_3',\theta_1'',\theta_3'')}
\\&=&
\sum_{j=1}^nS(Z_\tjm,\theta_1,\theta_3)^{-1}
\bigg[\cald_j(\theta_1',\theta_2',\theta_3'),\>
\left(\begin{array}{c}0\\
2^{-1}\partial_1^2L_H(Z_\tjm,\theta_1'',\theta_2'',\theta_3'')
\end{array}\right)
\bigg]
\\&=&
\sum_{j=1}^nS(Z_\tjm,\theta_1,\theta_3)^{-1}
\bigg[\cald_j(\theta_1',\theta_2',\theta_3'),\>
\left(\begin{array}{c}0\\
4^{-1}H_{xx}(Z_\tjm,\theta_3'')[\partial_1^2C(Z_\tjm,\theta_1'')]
\end{array}\right)
\bigg],
\eeas
}
\beas 
\Psi_{11,3}(\theta_1,\theta_3,\theta_1',\theta_2',\theta_3')
&=&
\sum_{j=1}^n
\partial_1\big\{S^{-1}(\partial_1S)S^{-1}(Z_\tjm,\theta_1,\theta_3)\big\}
\big[\cald_j(\theta_1',\theta_2',\theta_3')^{\otimes2}-S(Z_\tjm,\theta_1',\theta_3')\big],
\eeas
\beas
\Psi_{11,4}(\theta_1,\theta_3)
&=&
\sum_{j=1}^n(S^{-1}(\partial_1S)S^{-1})(Z_\tjm,\theta_1,\theta_3)
\big[\partial_1S(Z_\tjm,\theta_1, \theta_3)\big],
\eeas
{\vblue 
\beas&&
{\vredd\Psi_{11,5}(\theta_1,\theta_3,\theta_1',\theta_2',\theta_3',\theta_1'',\theta_3'')}
\\&=&
{\fred 2}
\sum_{j=1}^n
(S^{-1}(\partial_1S)S^{-1})(Z_\tjm,\theta_1,\theta_3)
\bigg[\cald_j(\theta_1',\theta_2',\theta_3'),
\left(\begin{array}{c}0\\
2^{-1}\partial_1L_H(Z_\tjm,\theta_1'',\theta_2'',\theta_3'')
\end{array}\right)
\bigg]
\\&=&
{\fred 2}
\sum_{j=1}^n
(S^{-1}(\partial_1S)S^{-1})(Z_\tjm,\theta_1,\theta_3)
\bigg[\cald_j(\theta_1',\theta_2',\theta_3'),\>
\left(\begin{array}{c}0\\
{\fred4^{-1}}H_{xx}(Z_\tjm,\theta_3'')[\partial_1C(Z_\tjm,\theta_1'')]
\end{array}\right)
\bigg].
\eeas
}
\begin{en-text}
{\fred We sometimes keep parameters in notation even when some of them do not appear in a specific expression of the formula, if such an expression is not necessary for later use; e.g.. $\Psi_{1,1}$ does not depend on $\theta_2$ in fact. 
}
\end{en-text}

{\coloro
{\vblue 
\begin{lemma}\label{201911250513}
Suppose that $[A1]$ with $\bbi={\vredd(1,1,2,\yfbox{3},3,1)}
$, $[A2]$ and {\vredd$[A4]$ $(iii)$}
are satisfied. 
Then, 
for some number $n_0$, 
\bea\label{201911122311} 
 \sup_{n\geq n_0}\bigg\|
 \sup_{\theta_1\in \bbU_n^{(1)}(\delta)}\left|
n^{-1}\>\partial_1^2\bbH_n^{(1)}(\theta_1)
-n^{-1}\>\partial_1^2\bbH_n^{(1)}(\theta_1^*)
\right|\bigg\|_p
&=&
O(\delta)
\eea
as $\delta\down0$ for every $p>1$, where $\bbU_n^{(1)}(\delta)=\{\theta_1\in\Theta_1;|\theta_1-\theta_1^*|<\delta\}$. 
\end{lemma}
\proof 
Take a sufficiently large integer $n_0$ such that $\bbU_n^{(1)}(\delta)=U(\theta_1^*,\delta)\equiv\{\theta_1\in\bbR^{\sfp_1};|\theta_1-\theta_1^*|<\delta\}$ for all $n\geq n_0$ and $\delta\leq1$. 
Let $p>1$. 
Then, it is easy to show 
\bea\label{2206260936}&&
\sup_{n\geq n_0}
\bigg\|\sup_{\theta_1\in U(\theta_1^*,\delta)}\big|
n^{-1}h\Psi_{11,1}(\theta_1,\wh{\theta}_3^0,\theta_1,\wh{\theta}_3^0)
-n^{-1}h\Psi_{11,1}(\theta_1^*,\wh{\theta}_3^0,\theta_1^*,\wh{\theta}_3^0)\big|\bigg\|_p
\nn\\&=& 
\sup_nh\times O(\delta)\yeq O(\delta)
\eea
as $\delta\down0$. 

By using Lemma \lfive (b) of \cite{gloter2021adaptive} ($\bbi=(1,1,2,1,2,0)$)
and Lemma \ref{2206230342} ($\bbi=(1,1,2,1,3,{\vredd1})$), 
\bea\label{2206260954a}&&
\bigg\|\sup_{\theta_1\in U(\theta_1^*,\delta)}\big|
n^{-1}h^{1/2}\Psi_{11,2}(\theta_1,\wh{\theta}_3^0,\theta_1,\wh{\theta}_2^0,\wh{\theta}_3^0,\theta_1,\wh{\theta}_3^0)
-n^{-1}h^{1/2}\Psi_{11,2}(\theta_1^*,\wh{\theta}_3^0,\theta_1^*,\wh{\theta}_2^0,\wh{\theta}_3^0,\theta_1^*,\wh{\theta}_3^0)
\big|\bigg\|_p
\nn\\&=& 
\bigg\|\sup_{\theta_1\in U(\theta_1^*,\delta)}\big|
n^{-1}h^{1/2}\Psi_{11,2}(\theta_1,\wh{\theta}_3^0,\theta_1,\wh{\theta}_2^0,\wh{\theta}_3^0,\theta_1,\wh{\theta}_3^0)
-n^{-1}h^{1/2}\Psi_{11,2}(\theta_1,\wh{\theta}_3^0,\theta_1^*,\wh{\theta}_2^0,\wh{\theta}_3^0,\theta_1,\wh{\theta}_3^0)
\big|\bigg\|_p
\nn\\&&
+\bigg\|\sup_{\theta_1\in U(\theta_1^*,\delta)}\big|
n^{-1}h^{1/2}\Psi_{11,2}(\theta_1,\wh{\theta}_3^0,\theta_1^*,\wh{\theta}_2^0,\wh{\theta}_3^0,\theta_1,\wh{\theta}_3^0)
-n^{-1}h^{1/2}\Psi_{11,2}(\theta_1^*,\wh{\theta}_3^0,\theta_1^*,\wh{\theta}_2^0,\wh{\theta}_3^0,\theta_1^*,\wh{\theta}_3^0)
\big|\bigg\|_p
\nn\\&=& 
h^{1/2}O(\delta)+O(\delta)\yeq O(\delta)
\eea
as $\delta\down0$ uniformly in $n\geq n_0$. 

We have 
\bea\label{2206251031}
\sup_{\theta_1\in\ol{\Theta}_1}\big|\cald_j(\theta_1,\wh{\theta}_2^0,\wh{\theta}_3^0)\big|
&\leq&
\sup_{\theta_1\in\ol{\Theta}_1}\big|\cald_j(\theta_1,\wh{\theta}_2^0,\wh{\theta}_3^0)-\cald_j(\theta_1,\wh{\theta}_2^0,\theta_3^*)\big| 
+\sup_{\theta_1\in\ol{\Theta}_1}\big|\cald_j(\theta_1,\wh{\theta}_2^0,\theta_3^*)\big|
\nn\\&=&
O_{L^\inftym}(h^{-1/2})\times O_{L^\inftym}(n^{-1/2}h^{1/2})+O_{L^\inftym}(1)
\yeq O_{L^\inftym}(1)
\eea
by Lemma \lsix of \cite{gloter2021adaptive} ($\bbi=(0,0,0,0,2,1)$) and {\vredd$[A4]$ (iii)}, and Lemma \ref{2206241411} ($\bbi=(1,1,2,1,3,0)$). 
Moreover, 
\bea\label{2206251050}
\big|\cald_j(\theta_1,\wh{\theta}_2^0,\wh{\theta}_3^0)^{\otimes2}
-\cald_j(\theta_1^*,\wh{\theta}_2^0,\wh{\theta}_3^0)^{\otimes2}\big|
&\leq&
\big|\cald_j(\theta_1,\wh{\theta}_2^0,\wh{\theta}_3^0)-\cald_j(\theta_1^*,\wh{\theta}_2^0,\wh{\theta}_3^0)\big| 
\nn\\&&\times
\big|\cald_j(\theta_1,\wh{\theta}_2^0,\wh{\theta}_3^0)+\cald_j(\theta_1^*,\wh{\theta}_2^0,\wh{\theta}_3^0)\big|
\nn\\&\leq&
\big|\cald_j(\theta_1,\wh{\theta}_2^0,\wh{\theta}_3^0)-\cald_j(\theta_1^*,\wh{\theta}_2^0,\wh{\theta}_3^0)\big| 
\nn\\&&\times
\sup_{\theta_1\in\ol{\Theta}_1}\big|\cald_j(\theta_1,\wh{\theta}_2^0,\wh{\theta}_3^0)+\cald_j(\theta_1^*,\wh{\theta}_2^0,\wh{\theta}_3^0)\big|
\nn\\&=&
O_{L^\inftym}(h^{1/2})|\theta_1-\theta_1^*|
\eea
by (\ref{2206251031}) and the representation 
\beas
\cald_j(\theta_1,\theta_2,\theta_3) - \cald_j(\theta_1^*,\theta_2,\theta_3) 
&=&
\left(\begin{array}{c}
0\y
-4^{-1}h^{1/2}H_{xx}(Z_\tjm,\theta_3)[C(Z_\tjm,\theta_1)-C(Z_\tjm,\theta_1^*)]
\end{array}\right), 
\eeas
that is from (\ref{201906030041}). 
Thus, we obtain 
\bea\label{2206260954}
\sup_{n\geq n_0}\bigg\|\sup_{\theta_1\in U(\theta_1^*,\delta)}\big|
n^{-1}\Psi_{11,3}(\theta_1,\wh{\theta}_3^0,\theta_1,\wh{\theta}_2^0,\wh{\theta}_3^0)
-n^{-1}\Psi_{11,3}(\theta_1^*,\wh{\theta}_3^0,\theta_1^*,\wh{\theta}_2^0,\wh{\theta}_3^0)
\big|\bigg\|_p
&=& 
O(\delta)
\eea
as $\delta\down0$, 
from (\ref{2206251031}) and (\ref{2206251050}). 
Remark that the estimate (\ref{2206260954}) uses $\partial_1^3B$.

It is easy to see 
\bea\label{2206261057}
\sup_{n\geq n_0}\bigg\|\sup_{\theta_1\in U(\theta_1^*,\delta)}\big|
n^{-1} \Psi_{11,4}(\theta_1,\wh{\theta}_3^0)-n^{-1} \Psi_{11,4}(\theta_1^*,\wh{\theta}_3^0)\big|\bigg\|_p
&=& 
O(\delta)
\eea
as $\delta\down0$.

By using (\ref{2206251031}) and an estimate similar to (\ref{2206251050}), we also obtain 
\bea\label{2206261105}&&
\sup_{n\geq n_0}\bigg\|\sup_{\theta_1\in U(\theta_1^*,\delta)}\big|
n^{-1}h^{1/2}\Psi_{11,5}(\theta_1,\wh{\theta}_3^0,\theta_1,\wh{\theta}_2^0,\wh{\theta}_3^0,\theta_1,\wh{\theta}_3^0)
\nn\\&&\hspace{90pt}
-
n^{-1}h^{1/2}\Psi_{11,5}(\theta_1^*,\wh{\theta}_3^0,\theta_1^*,\wh{\theta}_2^0,\wh{\theta}_3^0,\theta_1,\wh{\theta}_3^0)
\big|\bigg\|_p
\yeq
O(\delta)
\eea
as $\delta\down0$.

Now we obtain (\ref{201911122311}) by 
combining (\ref{2206260936}), (\ref{2206260954a}), (\ref{2206260954}), (\ref{2206261057}) and (\ref{2206261105}) 
with the decomposition of $n^{-1}\partial_1^2\bbH_n^{(1)}(\theta_1)$ as 
\bea\label{2206261203}
n^{-1}\partial_1^2\bbH_n^{(1)}(\theta_1)
&=&
-n^{-1}h\Psi_{11,1}(\theta_1,\wh{\theta}_3^0,\theta_1,\wh{\theta}_3^0)
+n^{-1}h^{1/2}
\Psi_{11,2}(\theta_1,\wh{\theta}_3^0,\theta_1,\wh{\theta}_2^0,\wh{\theta}_3^0,\theta_1,\wh{\theta}_3^0)
\nn\\&&
+\half n^{-1}\Psi_{11,3}(\theta_1,\wh{\theta}_3^0,\theta_1,\wh{\theta}_2^0,\wh{\theta}_3^0)
-\half n^{-1} \Psi_{11,4}(\theta_1,\wh{\theta}_3^0)
\nn\\&&
- n^{-1}h^{1/2}
\Psi_{11,5}(\theta_1,\wh{\theta}_3^0,\theta_1,\wh{\theta}_2^0,\wh{\theta}_3^0,\theta_1,\wh{\theta}_3^0). 
\eea
This completes the proof. 
We remark that the term $-2^{-1}n^{-1} \Psi_{11,4}(\theta_1,\wh{\theta}_3^0)$ is the dominating term 
among the terms on the right-hand side of (\ref{2206261203}), as shown in the proof of Lemma \ref{2206261123}
below. 
\qed\halflineskip

\begin{lemma}\label{2206261123}
Suppose that $[A1]$ with $\bbi=(1,1,2,\yfbox{\vredd 2},3,1)
$, $[A2]$ and {\vredd$[A4]$ $(iii)$} are satisfied. 
Then, 
for some number $\ep_1$, 
\bea\label{2206261233}
n^{-1}\partial_1^2\bbH_n^{(1)}(\theta_1^*)
&=& 
-\Gamma_{11}+O_{L^\inftym}(n^{-\ep_1})
\eea
as $n\to\infty$. 
\end{lemma}
\proof 
We will use the decomposition (\ref{2206261203}) of $n^{-1}\partial_1^2\bbH_n^{(1)}(\theta_1)$ 
evaluated at $\theta_1=\theta_1^*$. 
Since the second derivative of $B$ in $\theta_1$ is involved, we need $[A1]$ with $\bbi=(0,0,0,2,0,0)$ at least. 
Obviously, 
\bea\label{2206261140}
n^{-1}h\Psi_{11,1}(\theta_1^*,\wh{\theta}_3^0,\theta_1^*,\wh{\theta}_3^0)
&=&
O_{L^\inftym}(h), 
\eea
and with the help of (\ref{2206251031}) ($\bbi=(1,1,2,1,3,1)$), 
\bea\label{2206261145}
n^{-1}h^{1/2}
\Psi_{11,2}(\theta_1^*,\wh{\theta}_3^0,\theta_1^*,\wh{\theta}_2^0,\wh{\theta}_3^0,\theta_1^*,\wh{\theta}_3^0)
&=&
O_{L^\inftym}(h^{1/2}). 
\eea

Since, as already showed at (\ref{2206251031}) by 
Lemma \lsix of \cite{gloter2021adaptive} ($\bbi=(0,0,0,0,2,1)$) and {\vredd$[A4]$ (iii)}, 
\beas
\sup_{\theta_1\in\ol{\Theta}_1}\big|\cald_j(\theta_1,\wh{\theta}_2^0,\wh{\theta}_3^0)-\cald_j(\theta_1,\wh{\theta}_2^0,\theta_3^*)\big| 
&=&
O_{L^\inftym}(n^{-1/2}), 
\eeas
we obtain
\bea\label{2206261152}
\cald_j(\theta_1^*,\wh{\theta}_2^0,\wh{\theta}_3^0)-\cald_j(\theta_1^*,\theta_2^*,\theta_3^*)
&=& 
O_{L^\inftym}(n^{-1/2})+O_{L^\inftym}(h^{1/2})\yeq O_{L^\inftym}(h^{1/2})
\eea
by Lemma \lfive (b) ($\bbi=(1,1,2,1,2,0)$). 
Then 
\bea\label{2206261157}
n^{-1}\Psi_{11,3}(\theta_1^*,\wh{\theta}_3^0,\theta_1^*,\wh{\theta}_2^0,\wh{\theta}_3^0)
&=&
n^{-1}\Psi_{11,3}(\theta_1^*,\wh{\theta}_3^0,\theta_1^*,\theta_2^*,\theta_3^*)+O_{L^\inftym}(h^{1/2})
\nn\\&=&
n^{-1}\Psi_{11,3}(\theta_1^*,\theta_3^*,\theta_1^*,\theta_2^*,\theta_3^*)+O_{L^\inftym}(h^{1/2})
\quad(\text{{\vredd$[A4]$ (iii))}}
\nn\\&=&
O_{L^\inftym}(n^{-1/2})+O_{L^\inftym}(h^{1/2})
\quad(\text{orthogonality})
\nn\\&=&
O_{L^\inftym}(h^{1/2}). 
\eea

Moreover, 
\bea\label{2206261204}
n^{-1} \Psi_{11,4}(\theta_1^*,\wh{\theta}_3^0)
&=&
n^{-1} \Psi_{11,4}(\theta_1^*,\theta_3^*)+O_{L^\inftym}(n^{-1/2}h^{1/2})
\eea
by {\vredd $[A4]$ (iii)}. 
By the same argument as (\ref{22060241235}) with the mixing property and 
the equality $E[\Psi_{11,4}(\theta_1^*,\theta_3^*)]=\Gamma_{11}$ of (\ref{2206261226}), consequently we obtain 
\bea\label{2206261216}
n^{-1} \Psi_{11,4}(\theta_1^*,\wh{\theta}_3^0)
&=&
2\Gamma_{11}+O_{L^\inftym}(n^{-\ep_0/2}).
\eea

The estimates (\ref{2206251031}) gives 
\bea\label{2206261206}
n^{-1}h^{1/2}\Psi_{11,5}(\theta_1^*,\wh{\theta}_3^0,\theta_1^*,\wh{\theta}_2^0,\wh{\theta}_3^0,\theta_1^*,\wh{\theta}_3^0).
&=&
O_{L^\inftym}(h^{1/2}). 
\eea

We now obtain (\ref{2206261233}) for some positive constant $\ep_1$, 
from (\ref{2206261140}), (\ref{2206261145}), (\ref{2206261157}), (\ref{2206261216}) and (\ref{2206261206}). 
\qed\halflineskip
}
\begin{en-text}
We may assume that $n$ is sufficiently large and $r_n\geq(nh)^{-1/2}\geq h^{1/2}$. 
\end{en-text}
\begin{en-text}
We will use Condition $[A4^{{\fred\sharp}}]$ {\colorg for $\wh{\theta}_2^0$ and $\wh{\theta}_3^0$, 
and the estimate $|\theta_1-\theta_1^*|<r_n$ for $\theta_1\in U(\theta_1^*,r_n)$.}
Then 
\beas &&
\sup_{\theta_1\in U(\theta_1^*,r_n)}\big|n^{-1}\>\partial_1^2\bbH_n^{(1)}(\theta_1)+\Gamma_{11}\big|
\\&\leq&
O_p(h)
\\&&
+n^{-1}h^{1/2}
\sup_{\theta_1\in U(\theta_1^*,r_n)}\big|\Psi_{11,2}(\theta_1,\wh{\theta}_3^0,\theta_1^*,\theta_2^*,\theta_3^*,\theta_1,\wh{\theta}_2^0,\wh{\theta}_3^0)\big|
\nn\\&&\hspace{70pt}
+h^{1/2}O_p(n^{-1/2}+h^{1/2})
\quad(\text{Lemmas }\lsix\text{ of \cite{gloter2021adaptive} and } \lfive(b)
)
\\&&
+n^{-1}\sup_{\theta_1\in U(\theta_1^*,r_n)}\big|\Psi_{11,3}(\theta_1,\wh{\theta}_3^0,\theta_1^*,\theta_2^*,\theta_3^*)\big|
+O_p(h^{1/2}+{\fred n^{-1/2})+O_p(r_n)}
\\&&
\quad(\text{Lemmas }\lsix {\fred \text{ and }} \lfive(b))
\\&&
+\bigg|-\half n^{-1}{\coloro \Psi_{11,4}(\theta_1^*,\theta_3^*)}+\Gamma_{11}\bigg|+O_p(r_n)
{\fred+O_p(n^{-1/2}h^{1/2})}
\\&&
+n^{-1}h^{1/2}
\sup_{\theta_1\in U(\theta_1^*,r_n)}\big|\Psi_{11,5}(\theta_1,\wh{\theta}_3^0,\theta_1,\theta_2^*,\theta_3^*,\theta_1,\wh{\theta}_2^0,\wh{\theta}_3^0)\big|
\nn\\&&\hspace{90pt}
+{\fred h^{1/2}}O_p(h^{1/2}+n^{-1/2})
\quad(\text{Lemmas }\lsix\text{ and } \lfive(b))
\\&=&
O_p(h)
\\&&
+O_p(h^{1/2})\quad(\text{Lemma }\lthree(b))
\\&&
{\fred +O_p(h^{1/2})}
+O_p(n^{-1/2})+O_p(r_n)
\\&&
\quad(\text{random field argument with orthogonality})
\\&&
+o_p(1)\quad(\text{Lemma }\lone(a))
\\&&
+O_p(h^{1/2})\quad(\text{Lemma }\lfive(b))
\\&=&
o_p(1)
\eeas
\end{en-text}
\begin{en-text}
We remark that the used lemmas and appearing functions here 
require the regularity indices $(i_A,j_A,i_B,j_B,i_H,j_H)$ 
for $[A1]$ as follows: 
$(1,0,1,0,3,0)$ for Lemma \lthree(b); 
$({\colorg1,1,2},1,2,{\colorg0})$ for Lemma \lfive(b); 
$(0,0,0,0,2,1)$ for Lemma \lsix; 
$j_B=3$, $j_H=1$ for random field argument for $\Psi_{11,3}$. 
\end{en-text}
\begin{en-text}
\koko{\coloro If we apply the same machinery as in the proof of Lemma \ref{201906081400}, 
it is easy to obtain the result. 
$[$It is remarked that $\partial_1^2$ appears in $\Psi_{11,2}$ and $\Psi_{11,3}$. 
Uniform-in-$\theta_1$ estimate for $\Psi_{11,2}$ is simple since it has the factor $h^{1/2}$ in front of it.  
On the other hand, we use random field argument for $\Psi_{11,3}$ after making the martingale differences. 
We need $\partial_1^3$ at this stage. $]$
For the second assertion, the argument becomes local by Lemma \ref{201906070327}, then 
Lemma \ref{201906081400} and the convergence (\ref{201911122311}) gives it 
by Taylor's formula. 
}
\end{en-text}
}

{\colorg 

{\vblue 
\subsection{CLT and $L^\inftym$-boundedness of $n^{-1/2}\partial_1\bbH^{(1)}_n(\theta_1^*)$}
\begin{lemma}\label{202001181833}
Suppose that $[A1]$ with $\bbi=(1,1,2,1,3,1)$, $[A2]$ and {\vredd $[A4]$ $(ii), (iii)$}
are satisfied. 
Then
\bea\label{202001181432} 
n^{-1/2}\partial_1\bbH^{(1)}_n(\theta_1^*)&=& O_{L^\inftym}(1)
\eea
as $n\to\infty$.
\end{lemma}
\proof 
First, 
\beas&&
n^{-1/2}\Psi_{1,2}(\theta_1^*,\wh{\theta}_3^0,\theta_1^*,\wh{\theta}_2^0,\wh{\theta}_3^0)
\nn\\&=&
n^{-1/2}\Psi_{1,2}(\theta_1^*,\wh{\theta}_3^0,\theta_1^*,\wh{\theta}_2^0,\theta_3^*)+O_{L^\inftym}(1)
\quad(\text{{\vredd$[A4]$ (iii)}, Lemma \lsix of \cite{gloter2021adaptive} and Lemma \ref{2206241411}})
\\&=&
n^{-1/2}\Psi_{1,2}(\theta_1^*,\wh{\theta}_3^0,\theta_1^*,\theta_2^*,\theta_3^*)+O_{L^\inftym}(1)
\quad(\text{{\vredd$[A4]$ (ii)}, Lemma \lfive (b) of \cite{gloter2021adaptive} and Lemma \ref{2206241411}})
\\&=&
n^{-1/2}\Psi_{1,2}(\theta_1^*,\theta_3^*,\theta_1^*,\theta_2^*,\theta_3^*)+O_{L^\inftym}(h^{1/2})+O_{L^\inftym}(1)
\quad(\text{{\vredd$[A4]$ (iii)}, Lemma \lthree or Lemma \ref{2206241411}})
\\&=&
O_{L^\inftym}(1)
\quad(\text{the Burkholder-Davis-Gundy inequality}).
\eeas

Next, in the same way, 
we use {\vredd$[A4]$ (ii), (iii)} together with Lemmas \lsix and \lfive (b) of \cite{gloter2021adaptive} 
to show 
\beas &&
n^{-1/2}h^{1/2}\Psi_{1,1}(\theta_1^*,\wh{\theta}_3^0,\theta_1^*,\wh{\theta}_2^0,\wh{\theta}_3^0)
\nn\\&=& 
n^{-1/2}h^{1/2}\Psi_{1,1}(\theta_1^*,\wh{\theta}_3^0,\theta_1^*,\wh{\theta}_2^0,\theta_3^*)+O_{L^\inftym}(h^{1/2})
\quad(\text{{\vredd$[A4]$ (iii)}, Lemma \lsix of \cite{gloter2021adaptive} and Lemma \ref{2206241411}})
\nn\\&=& 
n^{-1/2}h^{1/2}\Psi_{1,1}(\theta_1^*,\wh{\theta}_3^0,\theta_1^*,\theta_2^*,\theta_3^*)+O_{L^\inftym}(h^{1/2})
\quad(\text{{\vredd$[A4]$ (ii)}, Lemma \lfive (b) of \cite{gloter2021adaptive} and Lemma \ref{2206241411}})
\nn\\&=& 
n^{-1/2}h^{1/2}\Psi_{1,1}(\theta_1^*,\theta_3^*,\theta_1^*,\theta_2^*,\theta_3^*)+O_{L^\inftym}(h)+O_{L^\inftym}(h^{1/2})
\quad(\text{{\vredd$[A4]$ (iii)}})
\nn\\&=&
O_{L^\inftym}(h^{1/2})
\quad(\text{the Burkholder-Davis-Gundy inequality}).
\eeas

Finally, we use the expression 
\beas 
n^{-1/2}\partial_1\bbH^{(1)}_n(\wh{\theta}_1^0)
&=&
n^{-1/2}h^{1/2}\Psi_{1,1}(\wh{\theta}_1^0,\wh{\theta}_3^0,\wh{\theta}_1^0,\wh{\theta}_2^0,\wh{\theta}_3^0)
+n^{-1/2}\Psi_{1,2}(\wh{\theta}_1^0,\wh{\theta}_3^0,\wh{\theta}_1^0,\wh{\theta}_2^0,\wh{\theta}_3^0). 
\eeas
to verify (\ref{202001181432}). 
\qed\halflineskip
}

{\vblue 
Since our condition $[A4]$ is stronger than $[A4^\sharp]$ of Gloter and Yoshida \cite{gloter2021adaptive}, 
the following lemma is a corollary to Lemma 7.14 of \cite{gloter2021adaptive}. 
\begin{lemma}\label{202001181156}
Suppose that $[A1]$ with $\bbi=\colorg{(1,1,2,1,3,1)}$, $[A2]$ and {\vredd$[A4]$}
are satisfied. 
Then 
\bea\label{202001181417} 
n^{-1/2}\partial_1\bbH^{(1)}_n(\theta_1^*)-M_n^{(1)} &\to^p& 0
\eea
as $n\to\infty$. 
In particular, 
\bea\label{20200118141}
n^{-1/2}\partial_1\bbH^{(1)}_n(\theta_1^*) &\to^d& N\big(0,\Gamma_{11}\big)
\eea
as $n\to\infty$. 
\end{lemma}
}
\begin{en-text}
\proof 
We have 
\beas &&
{\sf E}_j(\theta_2,\theta_3) 
\\&:=&
\cald_j(\theta_1^*,\theta_2,\theta_3) -\cald_j(\theta_1^*,\theta_2^*,\theta_3^*) 
\\&=&
\left(\begin{array}{c}
h^{1/2}\big(A(Z_\tjm,\theta_2^*)-A(Z_\tjm,\theta_2)\big)
\\
\left\{\begin{array}{c}
h^{-1/2}\big(H(Z_\tjm,\theta_3^*)-H(Z_\tjm,\theta_3)\big)
\\
+2^{-1}h^{1/2}\big(L_H(Z_\tjm,\theta_1^*,\theta_2^*,\theta_3^*)
-L_H(Z_\tjm,\theta_1^*,\theta_2,\theta_3)\big)
\end{array}
\right\}
\end{array}\right).
\eeas
Define the random field $\Xi_n(u_2,u_3)$ on $(u_2,u_3)\in U(0,1)^2$ by 
\beas &&
\Xi_n(u_2,u_3)
\\&=&
n^{-1/2}\sum_{j=1}^n\big(S^{-1}(\partial_1S)S^{-1}\big)(Z_\tjm,\theta_1^*,\theta_3^*+r_n^{(3)}u_3)
\nn\\&&\hspace{70pt}\cdot
\bigg[\cald_j(\theta_1^*,\theta_2^*,\theta_3^*)
\otimes{\sf E}_j(\theta_2^*+r_n^{(2)}u_2,\theta_3^*+r_n^{(3)}u_3)\bigg]
\eeas
{\tgreen where} $r_n^{(2)}=(nh)^{-1/2}\log(nh)$ and $r_n^{(3)}=n^{-1/2}h^{1/2}\log(nh)$. 
Then the Burkholder-Davis-Gundy inequality gives 
\beas 
\lim_{n\to\infty}\sup_{(u_2,u_3)\in U(0,1)^2}\sum_{i=0,1}
\big\|\partial_{(u_2,u_3)}^i\Xi_n(u_2,u_3)\big\|_p
&=&
0,
\eeas
which implies 
\beas 
\sup_{(u_2,u_3)\in U(0,1)^2}\big|\Xi_n(u_2,u_3)\big|  &\to^p& 0,
\eeas
and hence {\tgreen under $[A4^\sharp]$,} 
\bea\label{202001181318}
n^{-1/2}\sum_{j=1}^n\big(S^{-1}(\partial_1S)S^{-1}\big)(Z_\tjm,\theta_1^*,\wh{\theta}_3^0)
\bigg[\cald_j(\theta_1^*,\theta_2^*,\theta_3^*)
\otimes{\sf E}_j(\wh{\theta}_2^0,\wh{\theta}_3^0)\bigg]
&\to^p&
0
\nn\\&&
\eea
as $n\to\infty$. 
It is easier to see 
\bea\label{202001181322}
n^{-1/2}\sum_{j=1}^n\big(S^{-1}(\partial_1S)S^{-1}\big)(Z_\tjm,\theta_1^*,\wh{\theta}_3^0)
\big[{\sf E}_j(\wh{\theta}_2^0,\wh{\theta}_3^0)^{\otimes2}\big]
&\to^p&
0
\eea
as $n\to\infty$. 
From (\ref{202001181318}) and (\ref{202001181322}), 
\bea\label{202001181326}&&
n^{-1/2}\Psi_{1,2}(\theta_1^*,\wh{\theta}_3^0,\theta_1^*,\wh{\theta}_2^0,\wh{\theta}_3^0)
\nn\\&=&
n^{-1/2}\Psi_{1,2}(\theta_1^*,\wh{\theta}_3^0,\theta_1^*,\theta_2^*,\theta_3^*)
+o_p(1)
\nn\\&=&
n^{-1/2}\Psi_{1,2}(\theta_1^*,\theta_3^*,\theta_1^*,\theta_2^*,\theta_3^*)
+o_p(1)
\eea
as $n\to\infty$, where the last equality is by $[A4^\sharp]$. 

On the other hand, by $[A4^\sharp]$ and Lemmas \lsix of \cite{gloter2021adaptive} and \lfive (b), 
we obtain 
\bea\label{202001181400}&&
n^{-1/2}h^{1/2}\Psi_{1,1}(\theta_1^*,\wh{\theta}_2^0,\wh{\theta}_3^0,\theta_1^*,\wh{\theta}_2^0,\wh{\theta}_3^0)
\nn\\&=&
n^{-1/2}h^{1/2}\Psi_{1,1}(\theta_1^*,\wh{\theta}_2^0,\wh{\theta}_3^0,\theta_1^*,\theta_2^*,\theta_3^*)
+o_p(1).
\eea
By random field argument applied to the first term on the right-hand side of (\ref{202001181400}),  
\bea\label{202001181405}
n^{-1/2}h^{1/2}\Psi_{1,1}(\theta_1^*,\wh{\theta}_2^0,\wh{\theta}_3^0,\theta_1^*,\wh{\theta}_2^0,\wh{\theta}_3^0)
&=&
o_p(1). 
\eea

Consequently, from (\ref{202001181326}) and (\ref{202001181405}), 
we obtain the convergence (\ref{202001181417}) since 
\beas 
n^{-1/2}\partial_1\bbH^{(1)}_n(\theta_1^*)
&=&
n^{-1/2}h^{1/2}\Psi_{1,1}(\theta_1^*,\wh{\theta}_2^0,\wh{\theta}_3^0,\theta_1^*,\wh{\theta}_2^0,\wh{\theta}_3^0)
\nn\\&&
+n^{-1/2}\Psi_{1,2}(\theta_1^*,\wh{\theta}_3^0,\theta_1^*,\wh{\theta}_2^0,\wh{\theta}_3^0)
\\&=&
n^{-1/2}\Psi_{1,2}(\theta_1^*,\theta_3^*,\theta_1^*,\theta_2^*,\theta_3^*)
+o_p(1)
\\&=&
M_n^{(1)}+o_p(1)
\eeas
by using Lemmas {\fred\lfour} and \lfive (a). 
Convergence (\ref{20200118141}) follows from this fact and 
Lemma \lone with $[A2]$. 
\qed\halflineskip
\end{en-text}
}
\halflineskip

{\vblue 
\subsection{Proof of Theorem \ref{2206261240}}
Theorem \ref{2206261240} is proved by the simplified quasi-likelihood analysis in \cite{yoshida2021simplified}. 
More precisely,  Lemmas \ref{202001181833}, \ref{22062511258}, \ref{201911250513} and \ref{2206261123} 
verify 
Condition $[T2]$ of \cite{yoshida2021simplified}. 
Condition $[U1]$ of \cite{yoshida2021simplified} is $[A3]$ (i$'$). 
Lemma \ref{202001181156} gives the convergence (3.3) in \cite{yoshida2021simplified}. 
In this way, Theorem 3.5 of \cite{yoshida2021simplified} applies to the proof of Theorem \ref{2206261240}. 
\qed
}
\color{black}

{\vblue 
\section{Proofs for Section \ref{202007281718}}
\subsection{Proof of Theorem \ref{201905291601}} 

Define a random field $\Y_n^{(1)}(\theta_1)$ by 
$
\Y_n^{(1)}(\theta_1)
= 
n^{-1}\big(\bH^{(1)}_n(\theta_1) -\Y_n^{(1)}(\theta_1^*)\big).
$
Then 
\beas
\Y_n^{(1)}(\theta_1)
&=& 
-\frac{1}{2n} \sum_{j=1}^n \bigg\{\big(C(Z_\tjm,\theta_1)^{-1}-C(Z_\tjm,\theta_1^*)^{-1}\big)\big[h^{-1}(\Delta_jX)^{\otimes2}\big]
+\log\frac{\det C(Z_\tjm,\theta_1)}{\det C(Z_\tjm,\theta_1^*)}\bigg\}. 
\eeas
Under $[A1]$ with $\bbi=(0,0,1,1,0.0)$ and $[A2]$ (i){\vredd-(iii)}, by using 
Lemma \lone (b) of \cite{gloter2021adaptive}, 
Sobolev's inequality, and the mixing property, 
we obtain
\bea\label{2206281136}
\sup_{\theta_1\in\ol{\Theta}_1}\big|\Y_n^{(1)}(\theta_1)-\bbY^{(1)}(\theta_1)\big|=O_{\vredd L^\inftym}(n^{-\ep_1})
\eea
as $n\to\infty$ for some positive constant $\ep_1$.

We have
\bea\label{2206281222}
n^{-1/2}\partial_1\bH_n^{(1)}(\theta_1^*)
&=&
-2^{-1}n^{-1/2}\sum_{j=1}^n(\partial_1C^{-1})(Z_\tjm,\theta_1^*)\big[h^{-1}(\Delta_jX)^{\otimes2}-C(Z_\tjm,\theta_1^*)\big].
\eea
From 
Lemma \ltwo of \cite{gloter2021adaptive} applied under $[A1]$ with $\bbi=(1,0,0,0,0,0)$ and $[A2]$ (i), 
\bea\label{2206281223} 
h^{-1/2}\Delta_jX
&=& 
h^{-1/2}\int_\tjm^\tj B(Z_t,\theta_1^*)dw_t
+
h^{1/2}A(Z_\tjm,\theta_2^*)
+r^{(\ref{201905280600})}_j
\eea
with 
\bea\label{201905280600}
r^{(\ref{201905280600})}_j
&=&
h^{-1/2}\int_\tjm^\tj \big(A(Z_t,\theta_2^*)-A(Z_\tjm,\theta_2^*)\big)dt, 
\eea
for which 
\bea\label{201905280613}
\sup_n\sup_j\big\|r^{(\ref{201905280600})}_j\big\|_p=O(h)
\eea 
for every $p>1$. 
From (\ref{2206281222}) and (\ref{2206281223}), 
\bea\label{2206281224}&&
n^{-1/2}\partial_1\bH_n^{(1)}(\theta_1^*)
\nn\\&=&
-2^{-1}n^{-1/2}\sum_{j=1}^n(\partial_1C^{-1})(Z_\tjm,\theta_1^*)
\bigg[\bigg(h^{-1/2}\int_\tjm^\tj B(Z_t,\theta_1^*)dw_t
+
h^{1/2}A(Z_\tjm,\theta_2^*)\bigg)^{\otimes2}
\nn\\&&\hspace{200pt}
-C(Z_\tjm,\theta_1^*)\bigg]+O_{L^\inftym}(n^{1/2}h)
\nn\\&=&
-2^{-1}n^{-1/2}\sum_{j=1}^n(\partial_1C^{-1})(Z_\tjm,\theta_1^*)
\bigg[\bigg(h^{-1/2}\int_\tjm^\tj B(Z_t,\theta_1^*)dw_t
)\bigg)^{\otimes2}
-C(Z_\tjm,\theta_1^*)\bigg]
\nn\\&&\hspace{200pt}
+O_{L^\inftym}(h^{1/2})+O_{L^\inftym}(n^{1/2}h), 
\eea
with the aid of the orthogonality. 

Let
\bea\label{201996160909}
L_B(z,\theta_1,\theta_2,\theta_3)
&=&
B_x(z,\theta_1)[A(z,\theta_2)]+\half B_{xx}(z,{\fred\theta_1})[C(z,\theta_1)]
+B_y(z,\theta_3)[H(z,\theta_3)].
\eea
Assuming $[A1]$ with $\bbi=(0,0,2,0,0,0)$, 
\beas
h^{-1/2}\int_\tjm^\tj B(Z_t,\theta_1^*)dw_t
&=& 
h^{-1/2}\int_\tjm^\tj B(Z_\tjm,\theta_1^*)dw_t
+h^{-1/2}\int_\tjm^\tj\int_\tjm^t (B_xB)(Z_s,\theta_1^*)dw_sdw_t
\nn\\&&
+h^{-1/2}\int_\tjm^\tj\int_\tjm^t L_B(Z_s,\theta_1^*,\theta_2^*,\theta_3^*)dsdw_t
\nn\\&=:& 
\bbI_{1,j}+\bbI_{2,j}+\bbI_{3,j},
\eeas
where $dw_sdw_t=[dw_s\otimes dw_t]$. 
The function $L_B$ involves $\partial_x^2B$. 
Since $\sup_{n\in\bbN}\sup_{j=1,...,n}\|\bbI_{3,j}\|_p=O(h)$ for every $p>1$, 
the contribution to (\ref{2206281224}) of 
the terms involving $\bbI_{3,j}$ in the expansion of $\big(\bbI_{1,j}+\bbI_{2,j}+\bbI_{3,j}\big)^{\otimes2}$ 
is of $O_{L^\inftym}(n^{1/2}h)=o_{L^\inftym}(1)$. 
Applying It\^o's formula to $\bbI_{1,j}\otimes \bbI_{2,j}$ and using orthogonality, we see 
the contribution of $\bbI_{1,j}\otimes \bbI_{2,j}$ to (\ref{2206281224}) is of $O_{L^\inftym}(h^{1/2})$. 
Similarly, the contribution of $\bbI_{2,j}^{\otimes2}$ to (\ref{2206281224}) is of $O_{L^\inftym}(h)$. 
Therefore, 
\bea\label{2206281429}
n^{-1/2}\partial_1\bH_n^{(1)}(\theta_1^*)
&=&
-2^{-1}n^{-1/2}\sum_{j=1}^n(\partial_1C^{-1})(Z_\tjm,\theta_1^*)
\big[\big(h^{-1/2}B(Z_\tjm\theta_1^*)\Delta_jw\big)^{\otimes2}
-C(Z_\tjm,\theta_1^*)\big]
\nn\\&&\hspace{200pt}
+o_{L^\inftym}(1)
\eea
Consequently, we obtain 
\bea\label{2206281444a}
n^{-1/2}\partial_1\bH_n^{(1)}(\theta_1^*)
&=& 
O_{L^\inftym}(1),
\eea
and also 
\bea\label{2206281444}
n^{-1/2}\partial_1\bH_n^{(1)}(\theta_1^*)
&\to^d& 
N_{\sfp_1}(0,\Gamma^{(1)})
\eea
as $n\to\infty$, by the martingale central limit theorem. 

Assume $[A1]$ with $\bbi=(0,0,1,3,0,0)$. Then 
\bea\label{2206281452}
n^{-1}\partial_1^2\bH_n^{(1)}(\theta_1)
&=&
-2^{-1}n^{-1}\sum_{j=1}^n(\partial_1^2C^{-1})(Z_\tjm,\theta_1)\big[h^{-1}(\Delta_jX)^{\otimes2}-C(Z_\tjm,\theta_1)\big]
\nn\\&&
+2^{-1}n^{-1}\sum_{j=1}^n(\partial_1C^{-1})(Z_\tjm,\theta_1)\big[\partial_1C(Z_\tjm,\theta_1)\big]
\eea
and it is not difficult to prove 
\bea\label{2206281457}
\sup_{n\geq n_0}\bigg\|
\sup_{\theta_1\in U^{(1)}(1)}\big|
n^{-1}\partial_1^2\bH_n^{(1)}(\theta_1^*+\delta u_1)-n^{-1}\partial_1^2\bH_n^{(1)}(\theta_1^*)\big|\bigg\|_p
&=& 
O(\delta)
\eea
as $\delta\down0$, for every $p>1$, for some $n_0\in\bbN$, 
where $U^{(1)}(r)=\{u_1\in \bbR^{\sfp_1};\>|u_1|<r\}$, 
by $L^p$-estimate of $\partial_1^i\big\{n^{-1}\partial_1^2\bH_n^{(1)}(\theta_1^*+\delta u_1)-n^{-1}\partial_1^2\bH_n^{(1)}(\theta_1^*)\}$ for $i=0,1$. 

Moreover, by rewriting $\Delta_jX$ in (\ref{2206281452}) with (\ref{201905262213}), 
neglecting the drift term, and 
applying the Burkholder-Davis-Gundy inequality to the wighted sum of 
\beas
\bigg(h^{-1/2}\int_\tjm^\tj B(Z_t,\theta_1^*)dw_t\bigg)^{\otimes2}-h^{-1}\int_\tjm^\tj C(Z_t,\theta_1^*)dt
\eeas
on the right-hand side of (\ref{2206281452}), 
estimating the sum involving the gap $h^{-1}\int_\tjm^\tj C(Z_t,\theta_1^*)dt-C(Z_\tjm,\theta_1^*)$ 
with Lemma \lone (a), 
and by the mixing property to the second term on the right-hand side of (\ref{2206281452}), we obtain 
\bea\label{2206281512}
n^{-1}\partial_1^2\bH_n^{(1)}(\theta_1^*)+\Gamma^{(1)}
&=& 
O_{L^\inftym}(h^{1/2})
+O_{L^\inftym}(n^{-1/2})+O_{L^\inftym}(n^{-\ep_0/2})
\nn\\&=&
O_{L^\inftym}(n^{-\ep_2})
\eea
as $n\to\infty$, for some positive constant $\ep_2$. 

Consequently, Theorem 3.5 of Yoshida \cite{yoshida2021simplified} with  
Condition $[A3]$ (i), 
the estimates (\ref{2206281444a}), (\ref{2206281136}), (\ref{2206281457}) and (\ref{2206281512}), and 
the convergence (\ref{2206281444}) induces the results of Theorem \ref{201905291601}. 
\qed\halflineskip

\subsection{Proof of Theorem \ref{202002172010}}\label{2206281655}

Set $\Y^{(2)}_n(\theta_2)=T^{-1}\big(\bH^{(2)}_n(\theta_2)-\bH^{(2)}_n(\theta_2^*)\big)$. 
Then 
\bea\label{2206281842}
\Y^{(2)}_n(\theta_2)
&=&
\Y^{(\ref{2206281808})}_n(\wh{\theta}_1^0,\theta_2)+\Y^{(\ref{2206281809})}_n(\wh{\theta}_1^0,\theta_2),
\eea
where 
\bea\label{2206281808}
\Y^{(\ref{2206281808})}_n(\theta_1,\theta_2)
&=&
\frac{1}{T}\sum_{j=1}^n C(Z_\tjm,\theta_1)^{{\cred -1}}
\big[A(Z_\tjm,\theta_2)-A(Z_\tjm,\theta_2^*),\Delta_jX-hA(Z_\tjm,\theta_2^*)\big]
\nn\\&&
\eea
and 
\bea\label{2206281809}
\Y^{(\ref{2206281809})}_n(\theta_1,\theta_2)
&=&
-\frac{1}{2n}\sum_{j=1}^n C(Z_\tjm,\theta_1)^{{\cred -1}}
\big[\big(A(Z_\tjm,\theta_2)-A(Z_\tjm,\theta_2^*)\big)^{\otimes2}\big].
\eea
By {\vredd(\ref{2206281223})} with Lemma \ltwo of \cite{gloter2021adaptive}, under $[A1]$ with $\bbi=(1,0,0,0,0,0)$ and $[A2]$ (i), 
\bea\label{2206281857}
\sup_{\theta_2\in\ol{\Theta}_2}\big|\Y^{(\ref{2206281808})}_n(\wh{\theta}_1^0,\theta_2)\big|
&=&
\sup_{\theta_2\in\ol{\Theta}_2}\big|\Y^{(\ref{2206281810})}_n(\wh{\theta}_1^0,\theta_2)\big|
+O_{L^\inftym}(h^{1/2})
\eea
with 
\bea\label{2206281810}
\Y^{(\ref{2206281810})}_n(\theta_1,\theta_2)
&=&
\frac{1}{T}\sum_{j=1}^n C(Z_\tjm,\theta_1)^{{\cred -1}}
\bigg[A(Z_\tjm,\theta_2)-A(Z_\tjm,\theta_2^*),\int_\tjm^\tj B(Z_t,\theta_1^*)dw_t\bigg]. 
\nn\\&&
\eea
Since by Sobolev's inequality with the Burkholder-Davis-Gundy inequality, we obtain 
\bea\label{2206281858}
\sup_{(\theta_1,\theta_2)\in\ol{\Theta}_1\times\ol{\Theta}_2}\big|\Y^{(\ref{2206281810})}_n(\theta_1,\theta_2)\big|
&=&
O_{L^\inftym}(T^{-1/2})
\eea
with $[A1]$ for $\bbi=(0,1,0,1,0,0)$, it follows that 
\bea\label{2206281823}
\sup_{\theta_2\in\ol{\Theta}_2}\big|\Y^{(\ref{2206281808})}_n(\wh{\theta}_1^0,\theta_2)\big|
&=&
O_{L^\inftym}(T^{-1/2})+O_{L^\inftym}(h^{1/2})
\nn\\&=&
O_{L^\inftym}(T^{-1/2})\leq O_{L^\inftym}(n^{-\ep_0/2}).
\eea

We use $[A1]$ with $\bbi={\vredd(1,1,1,1,0,0)}$ 
and the mixing property to conclude 
\beas&&
\sup_{(\theta_1,\theta_2)\in\ol{\Theta}_1\times\ol{\Theta}_2}\bigg|
\Y^{(\ref{2206281809})}_n(\theta_1,\theta_2)
+\half\int_{\bbR^{\sfd_Z}}C(z,\theta_1)^{-1}\big[\big(A(z,\theta_2)-A(z,\theta_2^*)\big)^{\otimes2}\big]\nu(dz)\bigg|
\nn\\&=&
O_{L^\inftym}(h^{1/2})+O_{L^\inftym}(T^{-1/2})\yeq
O_{L^\inftym}(T^{-1/2}),
\eeas
and hence 
\bea\label{2206281836}
\sup_{\theta_2)\in\ol{\Theta}_2}\big|
\Y^{(\ref{2206281809})}_n(\wh{\theta}_1^0,\theta_2)
-\bbY^{(2)}(\theta_2)
\big|
&=&
O_{L^\inftym}(T^{-1/2})+O_{L^\inftym}(n^{-1/2})
\nn\\&=&
O_{L^\inftym}(T^{-1/2})\leq O_{L^\inftym}(n^{-\ep_0/2})
\eea
due to $[A4]$ (i). 
From (\ref{2206281842}), (\ref{2206281823}) and (\ref{2206281836}), 
\bea\label{2206281844}
\sup_{\theta_2)\in\ol{\Theta}_2}\big|\Y^{(2)}_n(\theta_2)-\bbY^{(2)}(\theta_2)
\big| &=& O_{L^\inftym}(n^{-\ep_0/2}). 
\eea

\begin{en-text}
\beas
T^{-1/2}\partial_2\bH^{(2)}_n(\theta_2) 
&=&
\sum_{j=1}^n C(Z_\tjm,\wh{\theta}_1^0)^{-1}
\big[\Delta_jX-hA(Z_\tjm,\theta_2),\partial_2A(Z_\tjm,\theta_2)\big].
\eeas
\end{en-text}

For uniformity in $\theta_1$, 
we can apply a method similar to (\ref{2206281857})-(\ref{2206281858}) to 
\bea\label{2206290105}
T^{-1/2}\partial_2\bH^{(2)}_n(\theta_2^*) 
&=&
T^{-1/2}\sum_{j=1}^n C(Z_\tjm,\wh{\theta}_1^0)^{-1}
\big[\Delta_jX-hA(Z_\tjm,\theta_2^*),\partial_2A(Z_\tjm,\theta_2^*)\big],
\eea
to obtain 
\bea\label{2206281900}
\big|T^{-1/2}\partial_2\bH^{(2)}_n(\theta_2^*) \big|
&=&
O_{L^\inftym}(1)
\eea
under $[A1]$ with $\bbi=(1,1,0,1,0,0)$; $\partial_1B$ is used for uniformity. 

We have 
\beas
T^{-1}\partial_2^2\bH^{(2)}_n(\theta_2) 
&=&
G_n^{(\ref{2206282235})}(\wh{\theta}_1^0,\theta_2)+G_n^{(\ref{2206282236})}(\theta_2), 
\eeas
where 
\bea\label{2206282235}
G_n^{(\ref{2206282235})}(\theta_1,\theta_2)
&=& 
T^{-1}\sum_{j=1}^n C(Z_\tjm,\theta_1)^{-1}
\bigg[\int_\tjm^\tj B(Z_t,\theta_2^*)dw_t,\partial_2^2A(Z_\tjm,\theta_2)\bigg]
\eea
and 
\bea\label{2206282236}
G_n^{(\ref{2206282236})}(\theta_2)
&=& 
T^{-1}\sum_{j=1}^n C(Z_\tjm,\wh{\theta}_1^0)^{-1}
\bigg[\int_\tjm^\tj A(Z_t,\theta_2^*)dt-hA(Z_\tjm,\theta_2),\partial_2^2A(Z_\tjm,\theta_2)\bigg]
\nn\\&&
-n^{-1}\sum_{j=1}^n C(Z_\tjm,\wh{\theta}_1^0)^{-1}
\big[\big(\partial_2A(Z_\tjm,\theta_2)\big)^{\otimes2}\big]. 
\eea
\begin{en-text}
\beas &&
T^{-1}\sum_{j=1}^n C(Z_\tjm,\wh{\theta}_1^0)^{-1}
\big[\Delta_jX-hA(Z_\tjm,\theta_2),\partial_2^2A(Z_\tjm,\theta_2)\big]
\nn\\&&
-n^{-1}\sum_{j=1}^n C(Z_\tjm,\wh{\theta}_1^0)^{-1}
\big[\big(\partial_2A(Z_\tjm,\theta_2)\big)^{\otimes2}\big]
\nn\\&=&
T^{-1}\sum_{j=1}^n C(Z_\tjm,\wh{\theta}_1^0)^{-1}
\bigg[\Delta_jX-\int_\tjm^\tj A(Z_t,\theta_2^*)dt,\partial_2^2A(Z_\tjm,\theta_2)\bigg]
\nn\\&&
+T^{-1}\sum_{j=1}^n C(Z_\tjm,\wh{\theta}_1^0)^{-1}
\bigg[\int_\tjm^\tj A(Z_t,\theta_2^*)dt-hA(Z_\tjm,\theta_2),\partial_2^2A(Z_\tjm,\theta_2)\bigg]
\nn\\&&
-n^{-1}\sum_{j=1}^n C(Z_\tjm,\wh{\theta}_1^0)^{-1}
\big[\big(\partial_2A(Z_\tjm,\theta_2)\big)^{\otimes2}\big]
\eeas
\end{en-text}
Let 
\beas 
\wt{G}_{n,\delta}^{(\ref{2206282235})}(\theta_1,u_2)
&=&
G_n^{(\ref{2206282235})}(\theta_1,\theta_2^*+\delta u_2)
\quad(u_2\in U^{(2)}(1))
\eeas
for sufficiently small $\delta>0$. 
Then,  
by the Burkholder-Davis-Gundy inequality,  
\beas 
\sum_{i,k=0}^1\sup_{(\theta_1,u_2)\in\ol{\Theta}_1\times U^{(2)}(1)}\big\|\partial_{u_2}^k
\partial_1^i\big\{\wt{G}_{n,\delta}^{(\ref{2206282235})}(\theta_1,u_2)-\wt{G}_{n,\delta}^{(\ref{2206282235})}(\theta_1,0)\big\}\big|\big\|_p
&=&
O(T^{-1/2})\delta \yeq O(\delta),
\eeas
for $U^{(2)}(1)=\{u_2\in\bbR^{\sfp_2};|u_2|<1\}$. 
Here we got the factor $\delta$ after the Burkholder-Davis-Gundy inequality in case $k=0$, and did it when applying $\partial_{u_2}$ in case $k=1$. 
The used differentiability is $\bbi=(0,3,0,1,0,0)$. 
Therefore, Sobolev's inequality gives 
\beas 
\sup_{n\geq n_0}\sup_{(\theta_1,u_2)\in\ol{\Theta}_1\times U^{(2)}(1)}
\big|\wt{G}_{n,\delta}^{(\ref{2206282235})}(\theta_1,u_2)-\wt{G}_{n,\delta}^{(\ref{2206282235})}(\theta_1,0)\big|=O_{L^\inftym}(\delta)
\eeas
as $\delta\down0$, for some integer $n_0$. 
This means 
\bea\label{2206290023}
\sup_{n\geq n_0}\sup_{(\theta_1,u_2)\in\ol{\Theta}_1\times U^{(2)}(1)}
\big|G_n^{(\ref{2206282235})}(\theta_1,\theta_2^*+\delta u_2)-G_n^{(\ref{2206282235})}(\theta_1,\theta_2^*)\big|=O_{L^\inftym}(\delta)
\eea
as $\delta\down0$. 
The same kind inequality as (\ref{2206290023}) for 
$G_n^{(\ref{2206282236})}(\theta_2^*+\delta u_2)-G_n^{(\ref{2206282236})}(\theta_2^*)$ is easy to show by taking advantage of the representation (\ref{2206282236}). 
As a result, 
\bea\label{2206290031}
\sup_{n\geq n_0}\sup_{u_2\in U^{(2)}(1)}\big|
T^{-1}\partial_2^2\bH^{(2)}_n(\theta_2^*+\delta u_2)-T^{-1}\partial_2^2\bH^{(2)}_n(\theta_2^*)\big|
&=& 
O_{L^\inftym}(\delta)
\eea
as $\delta\down0$. 
We used $[A1]$ with $\bbi=(0,3,0,1,0,0)$ for (\ref{2206290031}).

The repeatedly used argument about the uniform convergence also ensures 
\beas
\sup_{\theta_1\in\ol{\Theta}_1}\big|G_n^{(\ref{2206282235})}(\theta_1,\theta_2^*)\big|
&=& 
O_{L^\inftym}(T^{-1/2})
\eeas
on $[A1]$ with $\bbi={\vredd(0,2,0,1,0,0)}$. 
Furthermore, simply, 
\beas
\big|G_n^{(\ref{2206282236})}(\theta_2^*)+\Gamma_{22}
\big|
&=& 
O_{L^\inftym}(h^{1/2})+
O_{L^\inftym}(n^{-1/2})+O_{L^\inftym}(T^{-1/2})
\eeas
by using $[A1]$ with $\bbi=\vredd{(1,2,0,1,0,0)}$ and $[A4]$ (i). 
So we obtained 
\bea\label{2206290057}
T^{-1}\partial_2^2\bH^{(2)}_n(\theta_2^*) +\Gamma_{22}
&=& 
O_{L^\inftym}(n^{-\ep_2})
\eea
as $n\to\infty$ for some positive constant $\ep_2$, on $[A1]$ with $\bbi={\vredd(1,2,0,1,0,0)}$.

From (\ref{2206290105}), Lemma \lone (a) of \cite{gloter2021adaptive}, 
\bea\label{2206290119}&&
T^{-1/2}\partial_2\bH^{(2)}_n(\theta_2^*) 
\nn\\&=&
T^{-1/2}\sum_{j=1}^n C(Z_\tjm,\wh{\theta}_1^0)^{-1}
\big[\Delta_jX-hA(Z_\tjm,\theta_2^*),\partial_2A(Z_\tjm,\theta_2^*)\big]
\nn\\&=&
T^{-1/2}\sum_{j=1}^n C(Z_\tjm,\wh{\theta}_1^0)^{-1}
\bigg[\int_\tjm^\tj B(Z_t,\theta_1^*)dw_t,\partial_2A(Z_\tjm,\theta_2^*)\bigg]+O_p(n^{1/2}h)
\nn\\&=&
T^{-1/2}\sum_{j=1}^n C(Z_\tjm,\theta_1^*)^{-1}
\bigg[\int_\tjm^\tj B(Z_t,\theta_1^*)dw_t,\partial_2A(Z_\tjm,\theta_2^*)\bigg]+o_p(1)
\nn\\&&\hspace{150pt}
(\text{random field argument for $\theta_1$, and $\wh{\theta}_1^0\to^p\theta_1^*$})
\nn\\&=&
\wh{M}_n^{(2)}+O_p(h^{1/2})+o_p(1)
\nn\\&=&
\wh{M}_n^{(2)}+o_p(1)
\eea
on $[A1]$ with $\bbi={\vredd(1,1,1,1,0,0)}$. 
The martingale central limit theorem provides 
\bea\label{2206290129}
\wh{M}_n^{(2)}\to^d N_{\sfp_2}(0,\Gamma_{22}), 
\eea
where $[A1]$ requests 
$\bbi=(1,1,1,0,0,0)$. 

Finally, we apply Theorem 3.5 of \cite{yoshida2021simplified} with the facts 
(\ref{2206281844}), (\ref{2206281900}), (\ref{2206290031}), 
(\ref{2206290057}), (\ref{2206290119}) and (\ref{2206290129}) 
to complete the proof. 
\qed\halflineskip

\section{Proof of Theorem \ref{2206290542}}\label{202306160609}
By Theorems \ref{2206261240}, \ref{202002172010} and \ref{201906031711}, we have 
\bea\label{202306230102}
a_n^{-1}\big(\wh{\vartheta}-\theta^*\big) - \big(\Gamma_{11}^{-1}M_n^{(1)},\Gamma_{22}^{-1}\wh{M}_n^{(2)},\Gamma_{33}^{-1}M_n^{(3)}\big)
&\to^p& 
0
\eea
as $n\to\infty$. 
The $j$-th increment of the sum $\big({\vredd\wh{M}}_n^{(2)},M_n^{(3)}\big)$ are conditionally in the first-order Wiener chaos though {\vredd the principal part of}
$M_n^{(1)}$ conditionally in the second-order. 
Therefore, their increments are conditionally orthogonal, and hence 
$\big({\sf Z}_2,{\sf Z}_3\big)$ is independent of ${\sf Z}_1$. 
By some calculus with (\ref{202002171941}), (\ref{2206241544}) and (\ref{201906040958}), 
it can be shown that the $j$-th increments of ${\vredd\wh{M}}_n^{(2)}$ and $M_n^{(3)}$ are conditionally orthogonal. 
This implies independency between ${\sf Z}_2$ and ${\sf Z}_3$. 
Consequently, ${\sf Z}_1$, ${\sf Z}_2$ and ${\sf Z}_3$ are independent. 
Since $\big(M_n^{(1)},\wh{M}_n^{(2)},M_n^{(3)}\big)$ jointly converges, 
Theorem \ref{2206250928} gives the desired convergence. 
\halflineskip
}

{\xred 
\section{Proof of the inequality (\ref{E: comparison variance})}\label{202402141102}
The following proposition validates the inequality (\ref{E: comparison variance}). 
\begin{proposition}
Assume that $\sfd_Y \le \sfd_X$ and that $\mathop{rank}(H_x)= \sfd_Y$. Then, we have
\beas{\xblue
		{\rm Tr}\big\{\big(V^{-1}H_x(\partial_1C)H_x^\star V^{-1}H_x(\partial_1C)H_x^\star\big)(z,\theta_1^*,\theta_3^*) \big\}
\le{\rm Tr}\big\{\big(C^{-1}(\partial_1C)C^{-1}\partial_1C\big)(z,\theta_1^*)\big\} .
}%
\eeas
\end{proposition}	
\proof
	First, we recall some notations of linear algebra. For a matrix $W \in \mathbb{R}^{\sfd_Y} \otimes \mathbb{R}^{\sfd_X}$ with rank $\sfd_Y$ we introduce the right-inverse, $W^{-1,r}=W^\star(W W^\star)^{-1}$. It is such that $W W^{-1,r}=\mathbf{I}_{\mid{\mathbb{R}^{\sfd_Y}}}$ and 
	$W^{-1,r}W$ is the orthogonal projection on $(\ker(W))^\perp = \mathop{\text{Im}}(W^\star)$. Under the same assumption, the left-inverse of the injective matrix $W^\star$ is defined by
	$(W^\star)^{-1,l}=(W W^\star)^{-1}W$. It satisfies $(W^\star)^{-1,l}W^\star=\mathbf{I}_{\mid{\mathbb{R}^{\sfd_Y}}}$ and 
	 $W^\star(W^\star)^{-1,l}=  W^\star(W W^\star)^{-1}W=W^{-1,r}W=:\widetilde{P}$ is the orthogonal projection on $(\ker(W))^\perp$.
	
	We write $V=H_x C H_x^\star$ in the form $V=W W^\star$ with $W=H_x C^{1/2}$, using that $C$ is symmetric. From the inversibility of $C$ and the assumption on the rank of $H_x$, we deduce that $W \in \mathbb{R}^{\sfd_Y} \otimes \mathbb{R}^{\sfd_X}$ has rank $\sfd_Y$. Hence, it is possible to consider the right- (resp. left-) inverse of $W$ (resp. $W^\star$). Moreover, we have $ (W^\star)^{-1,l} W^{-1,r}= (W W^\star)^{-1}WW^\star(W W^\star)^{-1}=(WW^\star)^{-1}=V^{-1}$, form which we can derive
	\begin{align} \nonumber
		H_x^\star V^{-1} H_x&= H_x^\star (W^\star)^{-1,l} W^{-1,r} H_x=C^{-1/2} (H_x C^{1/2})^\star(W^\star)^{-1,l} W^{-1,r} (H_x C^{1/2}) C^{-1/2}
		\\ \label{E: HVH formula}
		&= C^{-1/2} W^\star (W^\star)^{-1,l} W^{-1,r} W C^{-1/2}=C^{-1/2} \widetilde{P} \widetilde{P} C^{-1/2},
	\end{align}
	where $ \widetilde{P}$ is the matrix of the orthogonal projection on  $(\ker(W))^\perp=(\ker(H_x C^{1/2}))^\perp$.
	We use this representation and we can write the trace in the L.H.S. of \eqref{E: comparison variance},
	\begin{align}\nonumber
		&\text{Tr}\big\{V^{-1}H_x(\partial_1C)H_x^\star V^{-1}
		H_x(\partial_1C)H_x^\star
		\} \\		\nonumber
	&	\qquad \qquad =
 	\text{Tr}\big\{H_x^\star V^{-1}H_x(\partial_1C)
	H_x^\star V^{-1}H_x(\partial_1C)\big\} \quad \text{(commutativity of the trace),}
		\\\nonumber
	&	\qquad \qquad = \text{Tr}\big\{
		C^{-1/2} \widetilde{P} C^{-1/2} (\partial_1 C) C^{-1/2} \widetilde{P}C^{-1/2} (\partial_1 C) \big\} \quad \text{by \eqref{E: HVH formula}, and $\widetilde{P}^2=\widetilde{P}$,}
			\\\nonumber
		&	\qquad \qquad = \text{Tr}\big\{
		 \widetilde{P} C^{-1/2} (\partial_1 C) C^{-1/2} \widetilde{P}C^{-1/2} (\partial_1 C) C^{-1/2}\big\} \quad \text{(commutativity of the trace),}
		 	\\ &  \label{E: trace withour Hx}
		 		\qquad \qquad = \text{Tr}\big\{
		 \widetilde{P} C^{-1/2} (\partial_1 C) C^{-1/2} \widetilde{P}\widetilde{P}C^{-1/2} (\partial_1 C) C^{-1/2}\widetilde{P}\big\} \qquad \text{(using $\widetilde{P}^2=\widetilde{P}$ and trace commutativity).}			 
\end{align}
Assume for simplicity of the exposition that the dimension of the parameter $\theta_1$ is $\sfp_1=1$. Then, $ R:= C^{-1/2} (\partial_1 C) C^{-1/2} \widetilde{P} 
	 \widetilde{P}C^{-1/2} (\partial_1 C) C^{-1/2}= (\widetilde{P}C^{-1/2} (\partial_1 C) C^{-1/2})^\star \widetilde{P}C^{-1/2} (\partial_1 C) C^{-1/2}$ is a nonnegative symmetric matrix, and we have
	 from \eqref{E: trace withour Hx},
	 \begin{equation*}
	 	\text{Tr}\big\{V^{-1}H_x(\partial_1C)H_x^\star V^{-1}
	 	H_x(\partial_1C)H_x^\star=
	 	\text{Tr} \big\{ \widetilde{P} R  \widetilde{P} 
	 	\big\}.
	 \end{equation*}
	Using that $R$ is nonnegative symmetric and that $P$ is a orthogonal projection, we can show that $	\text{Tr} \big\{ \widetilde{P} R  \widetilde{P} \big\} \le 	\text{Tr} \big\{ R \big\}$. It entails,
	\begin{align}\nonumber
			\text{Tr}\big\{V^{-1}H_x(\partial_1C)H_x^\star V^{-1}
		H_x(\partial_1C)H_x^\star \}& \le 	\text{Tr} \big\{ R \big\}
	\\ \nonumber
	&= \text{Tr}\{C^{-1/2} (\partial_1 C) C^{-1/2} \widetilde{P} 
	\widetilde{P}C^{-1/2} (\partial_1 C) C^{-1/2}  \}
		\\ &= \text{Tr}\{ 
	\widetilde{P}C^{-1/2} (\partial_1 C) C^{-1/2}C^{-1/2} (\partial_1 C) C^{-1/2} \widetilde{P}  \} \label{E: P remains to remove}
			\end{align}
	by commutativity of the trace operator.
	As $C^{-1/2} (\partial_1 C) C^{-1/2}C^{-1/2} (\partial_1 C) C^{-1/2}$ is a nonnegative symmetric matrix and $\tilde{P}$ an orthogonal projection, we get that the trace given by \eqref{E: P remains to remove} is increased if we remove the matrices $\widetilde{P} $ in the expression. It entails,
	\begin{equation}\label{E: scalar control Information} 
		\text{Tr}\big\{V^{-1}H_x(\partial_1C)H_x^\star V^{-1}
		H_x(\partial_1C)H_x^\star \} \le \text{Tr}\{ C^{-1/2} (\partial_1 C) C^{-1/2}C^{-1/2} (\partial_1 C) C^{-1/2}  \}
	\end{equation}
	and \eqref{E: comparison variance} follows.
	
	If $\sfp_1>1$, we get \eqref{E: comparison variance}  in a analogous way. Remarking that the matrix inequality \eqref{E: comparison variance} follows from the scalar  control 
	\begin{multline*}
		\sum_{1\le i ,j \le \sfp_1} u_i u_j 	\text{Tr}\big\{V^{-1}H_x(\partial_1C)_iH_x^\star V^{-1}
		H_x(\partial_1C)_jH_x^\star \}
		 \\ \le 
			\sum_{1\le i ,j \le \sfp_1} u_i u_j 
		\text{Tr}\{ C^{-1/2} (\partial_1 C)_i C^{-1/2}C^{-1/2} (\partial_1 C)_j C^{-1/2}  \}
	\end{multline*}
	for all $u \in \mathbb{R}^{\sfp_1}$. This scalar comparison, is equivalent to the expression
	\eqref{E: scalar control Information} if $\partial_1 C$ is replaced, in the latter, by the partial derivative along the vector $u$, that is $\sum_{1\le i \le \sfp_1} u_i (\partial_1 C)_i$. As \eqref{E: scalar control Information} holds true when $\partial_1 C$ denotes any directional derivative, we deduce the lemma for $\sfp_1 >1$.	
\qed\halflineskip
}

\section{{\vblue Basic estimations for $\cald_j(\theta_1,\theta_2,\theta_3)$}}\label{202001141618}
%
\begin{en-text}
{\tgreen 
We denote $U^{\otimes k}$ for $U\otimes \cdots\otimes U$ ($k$-times) for a tensor $U$. 
For tensors  
$S^1=(S^1_{i_{1,1},...,i_{1,d_1};j_{1,1},...,j_{1,k_1}})$, ..., 
$S^m=(S^m_{i_{m,1},...,i_{m,d_m};j_{m,1},...,j_{m,k_m}})$ 
and and a tensor $T=(T^{i_{1,1},...,i_{1,d_1},...,i_{m,1},...,i_{m,d_m}})$, we write 
\beas 
T[S^1,...,S^m]&=&T[S^1\otimes\cdots\otimes S^m]
\\&=&
\bigg(
\sum_{i_{1,1},...,i_{1,d_1},...,i_{m,1},...,i_{m,d_m}}
T^{i_{1,1},...,i_{1,d_1},...,i_{m,1},...,i_{m,d_m}}
S^1_{i_{1,1},...,i_{1,d_1};j_{1,1},...,j_{1,k_1}}
\\&&\hspace{3cm}\cdots S^m_{i_{m,1},...,i_{m,d_m};j_{m,1},...,j_{m,k_m}}\bigg)
_{j_{1,1},...,j_{1,k_1},...,j_{m,1},...,j_{m,k_m}}.
\eeas
This notation will be applied for a tensor-valued tensor $T$ as well. 
{\tgreen
Clearly, it has an advantage over the matrix notation since 
quite often the elements $S^1, ..., S^m$ have a long expression 
in the inference; the matrix notation repeats $S^i$s twice for the quadratic form,  
thrice for the cubic form, and so on. 
This notation was introduced in Yoshida \cite{yoshida2011polynomial} and 
adopted by many papers, e.g., 
Uchida \cite{uchida2010contrast}, 
Uchida and Yoshida \cite{uchida2012adaptive}, 
Yoshida \cite{yoshida2013martingale}, 
Masuda \cite{masuda2013convergence}, 
Nakakita and Uchida \cite{nakakita2019inference}, 
Kamatani and Uchida \cite{KamataniUchida2014}, 
Nualart and Yoshida \cite{nualart2019asymptotic}, 
Nakakita et al. \cite{nakakita2020quasi}
among many other papers. 
}
}
\end{en-text}
{\vblue 
\begin{lemma}\label{2206230342}
Suppose that $[A1]$ with $\bbi=(1,1,2,1,3,{\vredd1})$ and $[A2]$ $(i)$ are satisfied. 
Then 
\bea\label{2206230411} 
\sup_{(\theta_1,\theta_2,\theta_3)\in\overline{\Theta}_1\times\overline{\Theta}_2\times\overline{\Theta}_3}
\big|\cald_j(\theta_1,\theta_2,\theta_3)\big|
&\leq& 
h^{-1/2}r_{n,j}^{(\ref{2206230343})}
\eea
for 
some random variables $r_{n,j}^{(\ref{2206230343})}$ such that 
\bea\label{2206230343}
\sup_n\sup_j\big\|r_{n,j}^{(\ref{2206230343})}\big\|_p
&<&
\infty
\eea
for every $p>1$. 
\end{lemma}
\proof 
Use Lemmas 
\lsix\!\!, \lfive (b) and \lthree (b) of \cite{gloter2021adaptive}. 
\qed
\halflineskip

Lemma \lthree (b) and Lemma \lfive (b) of \cite{gloter2021adaptive} provides the following estimate. 
\begin{lemma}\label{2206241411}
Suppose that $[A1]$ with $\bbi=(1,1,2,1,3,0)$ and $[A2]$ $(i)$ are satisfied. 
Then 
\bea\label{0406241100}
\sup_{(\theta_1,\theta_2)\in\overline{\Theta}_1\times\overline{\Theta}_2}
\big|\cald_j(\theta_1,\theta_2,\theta_3^*)\big|
&=& 
O_{L^\inftym}(1)
\eea
as $n\to\infty$. 
\end{lemma}

}

\bibliographystyle{spmpsci}      
\bibliography{bibtex-20180426-20180615bis3++++}   

\end{document}